\pdfoutput=1
\documentclass[pdflatex,sn-mathphys-num]{sn-jnl}
\usepackage[utf8]{inputenc}
\usepackage[T1]{fontenc}
\usepackage{lmodern}
\usepackage{amsmath,amssymb,amsthm,mathtools}
\usepackage{booktabs,longtable,array,enumitem,xurl}
\usepackage[expansion=false]{microtype}
\theoremstyle{thmstyleone}
\newtheorem{theorem}{Theorem}[section]
\newtheorem{proposition}[theorem]{Proposition}
\newtheorem{lemma}[theorem]{Lemma}
\newtheorem{corollary}[theorem]{Corollary}
\theoremstyle{thmstylethree}

\newtheorem{example}[theorem]{Example}

\newtheorem{question}{Question}[section]
\newcommand{\D}{\mathbb D}\newcommand{\N}{\mathbb N}\newcommand{\Z}{\mathbb Z}\newcommand{\R}{\mathbb R}\newcommand{\F}{\mathbb F}
\newcommand{\dd}{\,d\mu}\newcommand{\eps}{\varepsilon}\newcommand{\ind}{\mathbf1}\newcommand{\wh}{\widehat}\newcommand{\norm}[1]{\left\|#1\right\|}\newcommand{\abs}[1]{\left|#1\right|}
\DeclareMathOperator{\supp}{supp}\DeclareMathOperator{\spec}{spec}\DeclareMathOperator{\Var}{Var}
\numberwithin{equation}{section}
\setlist{nosep}
\hypersetup{bookmarksdepth=2,pdftitle={Walsh Analysis in the Twenty-First Century: Structural Survey},pdfauthor={Sergo A. Episkoposian}}
\makeatletter
\let\ps@plain\ps@headings
\makeatother
\begin{document}
\title[Walsh Analysis: 2000--2026]{Walsh Analysis in the Twenty-First Century: Convergence, Summability, Nonlinear Approximation, Universality, and Digital Applications}
\author*[1,2]{\fnm{Sergo A.} \sur{Episkoposian}}\email{sergo.episkoposyan@polytechnic.am}
\affil*[1]{National Polytechnic University of Armenia}
\affil*[2]{Yerevan State University }
\abstract{Walsh analysis connects ordered Fourier convergence, summability, nonlinear approximation, and universal representation through a common character algebra. This survey examines developments from 2000--2026 against their classical foundations. Binary kernel geometry explains Hardy-space losses and summability thresholds; interval and packet methods control spectral selection; localized blocks produce correction and universal representations with distinct coefficient and exceptional-set requirements. Selected proofs and counterexamples compare these mechanisms and identify where an implication needs an additional hypothesis. The same distinctions organize digital integration, sparse measurements, Boolean spectra, and generalized Walsh bases. Quantitative thresholds, source conventions, and published theorem pointers are collected in comparative tables. The concluding problems separate source-listed questions, proposed refinements, and elementary obstructions. Supplementary derivations preserve normalization checks and finite covariance calculations while keeping the main exposition focused on the analytic mechanisms.}
\maketitle
\tableofcontents
\clearpage
\section{Scope, comparison principles, and main thresholds}
The Walsh characters support three different approximation problems: a prescribed order of Fourier partial sums, a selectable spectral support, and a representation chosen after the target. Their common algebra does not determine the same convergence behavior. Dyadic sums are conditional expectations, arbitrary sums have endpoint losses, and correction or coefficient selection can produce representations unavailable to an ordered Fourier expansion.

This survey concentrates on 2000--2026. Earlier results supply the foundations: Golubov--Efimov--Skvortsov \cite{GES}, Schipp--Wade--Simon--P\'al \cite{SWS}, and Weisz \cite{Weisz94,Weisz02}; Persson--Tephnadze--Weisz \cite{PTW22} provides the modern Hardy-space synthesis. The historical surveys \cite{Wade,WeiszSurvey04,Golubov18,GatSurvey16} and the accounts of greedy approximation, time--frequency analysis, and numerical integration \cite{GreedyReview25,Thiele06,DKS13,GodaSuzuki20} delimit the present perspective.

\subsection{The comparison framework}\label{sub:v2surveys}
Four data specify a theorem in this survey: the \emph{ordering} of the characters; the permitted \emph{selection} of frequencies, signs, or terms; the ambient \emph{topology}; and the \emph{quantifiers}, especially whether a set, weight, or coefficient sequence precedes the target. The organizing contribution is a comparison of mechanisms under these fixed data. Kernel localization and atoms expose sharp thresholds; time--frequency decompositions control moving spectral boundaries; constructive blocks and category address target-dependent representations.

Three frequency parameters recur: the largest active binary position, the numerical index, and Hamming weight. Numerical frequency and martingale scale have equivalent positive powers on each block, while Hamming weight has different multiplicities. Their exact comparisons, and the examples separating them, are developed below. Correction and universality receive detailed treatment as the author's research perspective, alongside independent developments in summability, interval analysis, and digital integration.

\subsection{Reading order, status, and conventions}\label{sub:readingroutes}
The progression is from characters and kernels to ordered convergence, summability, Hardy estimates, and smoothness; then to interval selection, multipliers, greedy approximation, and target-dependent representation. Coefficient recovery and multiparameter geometry precede the digital and noncompact applications. Section~\ref{sec:finite} develops finite transforms, Boolean spectra, sparse measurements, and generalized Walsh bases. Supplementary covariance calculations are in Appendix~\ref{app:quadratic-models}. The universality route uses Sections~\ref{sec:foundations}, \ref{sec:greedy}, and \ref{sec:universal}; maximal-operator questions use Sections~\ref{sec:convergence}--\ref{sec:hardy} and \ref{sec:timefreq}.

A named theorem is an external input unless its proof is supplied. Proofs of known results and elementary consequences are included when they explain a threshold, a construction, or an obstruction; their numbering carries no claim of novelty. The common-weight theorem is a combined expository formulation proved here. Paragraphs labeled ``proof mechanism'' identify an external proof strategy. The bibliography records published articles, books, proceedings, and standards. An author version may provide open access to an identified publication; it is not treated as a separate publication. Appendix~\ref{app:source-register} records normalization and source scope. Constants depend on the displayed parameters; no sharp dependence is claimed unless stated.

The selected-proof structure has three connected anchors: the binary-gap atom in Proposition~\ref{v2:gap}, the sharp Fej\'er boundary in Theorem~\ref{th:fejerhardy}, and the localized representation block in Lemma~\ref{lem:univblock}. A complete time--frequency proof of Walsh--Carleson is outside this selection; the relevant tree mechanism and original sources are identified in Section~\ref{sec:timefreq}.

\begin{longtable}{@{}p{.22\textwidth}p{.68\textwidth}@{}}
\caption{Principal notation and conventions}\label{tab:notation}\\
\toprule Symbol & Meaning\\\midrule\endfirsthead
\toprule Symbol & Meaning\\\midrule\endhead
$\D,\mu,x\oplus y$ & Compact binary group, probability Haar measure, and coordinatewise addition modulo two.\\
$I_m(x),I_m$ & Cylinder fixing the first $m$ coordinates; $I_m=I_m(0)$ and $\mu(I_m)=2^{-m}$.\\
$r_j,w_n,\kappa_n$ & Rademacher character, Walsh--Paley character, and Walsh--Kaczmarz function.\\
$\widehat f(n),S_N,D_N$ & Walsh coefficient, projection onto $[0,N)$, and its Dirichlet kernel.\\
$\sigma_N,K_N$ & Positive-index Fej\'er convention: $N^{-1}\sum_{k=1}^NS_k$ and its kernel.\\
$E_m,\Delta_m$ & Conditional expectation $S_{2^m}$ and difference $E_{m+1}-E_m$.\\
$u(n),\ell(n),d(n)$ & Highest and lowest active binary positions and $u(n)-\ell(n)$, for $n>0$.\\
$V(n),L_n$ & Binary variation and Lebesgue constant $\|D_n\|_1$.\\
$F^*,H^p$ & Maximal function of a martingale and its Hardy space; no integrable terminal function is presumed when $p\le1$.\\
$M_d,A_q^d$ & Positive dyadic maximal operator and spatial dyadic Muckenhoupt class.\\
$P_A,G_N,e_N$ & Fourier projection on $A$, thresholding selection, and best $N$-term approximation with free coefficients.\\
$E_N(f)_p$ & Best approximation by frequencies $[0,N)$; $e_m(f)_p$ in the dyadic notation uses $N=2^m$.\\
$T_a,\mathcal D,P_t$ & Walsh multiplier, numerical-frequency generator, and its semigroup; $P_t(x)$ also denotes the kernel when an argument is displayed.\\
$h(n),T_\rho$ & Hamming weight and Boolean noise multiplier $\rho^{h(n)}$.\\
$P^\perp,\mu_\beta(k)$ & Digital annihilator and sum of the largest $\beta$ nonzero digit positions; $\mu_\beta$ is not Haar measure.\\
$I_\alpha,J_\alpha$ & Scale potential and its positive probability-kernel normalization; $J_\alpha=E_0+2^{-\alpha}(I_\alpha-E_0)$.\\
$\operatorname{VP}_N$ & Low-frequency reproducing filter $2\sigma_{2N}-\sigma_N$.\\
$H_N,Q_N$ & Unnormalized Sylvester matrix and unitary matrix $N^{-1/2}H_N$.\\
$\lesssim_\theta,\asymp_\theta$ & Constants depend only on the listed fixed parameters; dependence on dimension or a weight characteristic is stated separately.\\\bottomrule
\end{longtable}

\paragraph{Normalization dictionary.}
The table above fixes Paley order unless another order is named. Kaczmarz order reverses lower binary digits within a scale; sequency uses the reflected Gray permutation in Proposition~\ref{v2:seq}. Our Fej\'er mean averages $S_1,\ldots,S_N$, so its multiplier at frequency $n$ is $(1-n/N)_+$. Reversed N\"orlund means use $q_{N-k}S_k$, whereas direct means use $q_{k-1}S_k$; they agree only under an explicit symmetry of the row. Haar measure is a probability on $\D$ and has $\mu(H)=1$ on the noncompact field. On $\F_2^m$ the average Fourier coefficient is $N^{-1}H_Nf$, distinct from the unitary vector $N^{-1/2}H_Nf$. These conventions determine every dimensional factor below.

\subsection{The question answered by the comparison}\label{sub:research-design}
The central analytical question is which part of a convergence or representation theorem is paid for by localization, which by orthogonality, and which by a choice made after the target. The comparisons below identify that step through a threshold, a norm loss, or an explicit obstruction.

For instance, complete-block averaging has norm one, while selecting equally many individual characters incurs the power in Theorem~\ref{th:greedyrate}. The correction block in Lemma~\ref{lem:correctionblock} pays instead with a small exceptional set. The representation block in Lemma~\ref{lem:univblock} spreads coefficients over a longer spectrum. These are three different costs for improving an approximation process.

\begin{longtable}{@{}p{.18\textwidth}p{.34\textwidth}p{.35\textwidth}@{}}
\caption{Research structure: problem, mechanism, and unresolved step}\label{tab:research-design}\\
\toprule Part & What the proved theory explains & Boundary carried forward\\\midrule\endfirsthead
\toprule Part & What the proved theory explains & Boundary carried forward\\\midrule\endhead
Ordered sums and means & Binary kernels, atoms, and matrix structure determine critical exponents and subsequence losses. & Full-sequence endpoints and geometry finer than a worst-case count; Section~\ref{sec:boundaries}.\\
Intervals and multipliers & Moving frequency boundaries need packet analysis; compactness needs operator-tail control. & Boundedness at one exponent versus control of tails at that exponent.\\
Greedy selection & Conditionality distinguishes coefficient thresholding from best sparse approximation. & Algorithm choice, quantitative loss, and what correction changes.\\
Representation & Localized blocks and common sets permit target-dependent selection. & Formal coefficients versus Fourier origin; a fixed alphabet versus ternary freedom.\\
Finite applications & Annihilators, local coherence, and Boolean spectra identify the relevant signal geometry. & A smoothness ball, a recovery model, or a specified component determines the guarantee.\\\bottomrule
\end{longtable}

The organizing contribution of the survey is this comparison and its explicit boundary tests. The selected proofs serve four purposes: show a sharp loss, identify a counterexample, expose a reusable construction, or connect two formulations. Standard auxiliary calculations are included when they fix a normalization needed later. The exposition does not require a new theorem to justify each section, and a newly displayed elementary consequence is not assigned research priority.

\paragraph{Relation to other surveys and selection of literature.}
The modern greedy survey \cite{GreedyReview25} organizes known characterizations around precise research questions, while the longer account \cite{TemlyakovActa08} develops approximation algorithms and their guarantees. The QMC survey \cite{DKS13} connects error criteria, weighted spaces, and constructive algorithms. The present article uses a broader Walsh-centered comparison, so its common structure is a mathematical problem and a mechanism rather than a prescribed page count. The bibliography is a selective expert synthesis: foundational monographs, source theorems needed for the stated comparisons, sharp counterexamples, and representative applications. It is not an exhaustive database census. Recent work is included through an identified published source; journal publication and independent proof verification remain separate questions. The problem register distinguishes a published open question from a question formulated here and from an already settled obstruction.

\subsection{Historical transitions and the modern comparison}\label{sub:historical-transitions}
The modern questions retain three features of the original system. Walsh's construction \cite{Walsh} provides bounded orthogonal functions taking only the values $1$ and $-1$. Paley's ordering \cite{Paley} makes complete initial blocks the functions of finitely many binary coordinates. Fine's analysis \cite{Fine,FineCesaro55} shows how the resulting kernels behave under ordinary summation and averaging. These are analytic features of a specified ordering, rather than consequences of orthogonality alone.

The pointwise line and the summability line then separate. Billard's $L^2$ theorem and Sj\"olin's refinement \cite{Billard,Sjolin} concern unaveraged Walsh--Fourier sums. Fej\'er convergence uses a different family of kernels, and the maximal and Hardy-space theories make its exceptional endpoints visible. Schipp's work and the monographs \cite{Schipp69,Schipp75,GES,SWS,Weisz94,Weisz02} connect explicit binary identities with martingale localization. The later sharp results in \cite{PTW22} develop that structure beyond a qualitative convergence theorem.

A second transition replaces convergence in a prescribed order by nonlinear choice. The framework of greedy bases \cite{KonyaginTemlyakov,Wojtaszczyk,GreedyReview25} distinguishes retention of large original coefficients from approximation with coefficients chosen anew. Temlyakov's survey \cite{TemlyakovRMS25} also treats algorithms that recompute approximants in the selected span. In Walsh analysis, these choices meet a highly localized family of Dirichlet polynomials and an independent-character family with different $L^p$ sizes. Correction and universal representation exploit further freedom: the function, the exceptional set, or the permitted terms may change. The classical representation background \cite{Talalyan60,Ulyanov64} and the source-specific constructions in Section~\ref{sec:universal} identify what is fixed before the target.

A third transition moves between exact finite character identities and an analytic guarantee. Digital nets use annihilators together with a function norm, while sparse measurements use restricted injectivity on a signal class. Both retain the finite Walsh table, but they constrain different parts of a measurement map. The numerical literature \cite{DickPillichshammer,DKS13,GodaSuzuki20} is therefore connected with the preceding Fourier theory by its coefficient geometry, not by an unconditional transfer of a pointwise convergence theorem.

\begin{longtable}{@{}p{.12\textwidth}p{.30\textwidth}p{.43\textwidth}@{}}
\caption{Selected published milestones and the transition they explain}\label{tab:chronology}\\
\toprule Date & Source or result & Role in the comparison\\\midrule\endfirsthead
\toprule Date & Source or result & Role in the comparison\\\midrule\endhead
1923 & Walsh \cite{Walsh} & Bounded orthogonal sign functions.\\
1932 & Paley \cite{Paley} & Ordering compatible with complete binary blocks.\\
1949--1955 & Fine \cite{Fine,FineCesaro55} & Character and kernel structure; Ces\`aro summability.\\
1967--1969 & Billard, Sj\"olin \cite{Billard,Sjolin} & Almost everywhere convergence of unaveraged sums.\\
1969--1975 & Schipp \cite{Schipp69,Schipp75} & Divergence and maximal summability mechanisms.\\
1987--1994 & \cite{GES,SWS,Weisz94} & Systematic Walsh and martingale foundations.\\
1999--2000 & \cite{KonyaginTemlyakov,Wojtaszczyk} & Greedy and quasi-greedy basis theory.\\
2000 & Thiele \cite{ThieleQuartile00} & Quartile analysis and a time--frequency proof architecture for Walsh convergence.\\
2002--2004 & \cite{MTTBiCarleson02,MTTBiest04} & Distinct Walsh models for bi-Carleson and trilinear biest operators.\\
2010--2015 & \cite{Plotnikov10,SkvortsovTulone15} & Coefficient recovery, exceptional sets, and generalized integration.\\
2012--2014 & \cite{DoLaceyLac,DoLaceyWeighted,DiPlinioLac,DiPlinioLerner14} & Lacunary endpoints, variation, and quantitative weights.\\
2016 & Osipov \cite{Osipov16} & Uniform synthesis over arbitrary disjoint spectral intervals.\\
2017 & \cite{GSY17} & Sharp-order high-smoothness digital integration under specified net conditions.\\
2022 & \cite{PTW22} & Modern Hardy-space and summability synthesis.\\
2025 & \cite{DFPR25,GreedyReview25,TemlyakovRMS25} & Sharp weighted interval estimates and modern approximation surveys.\\
2026 & \cite{HaoLiWeisz26,EpiVilenkin26,Getsadze26} & Generalized Hardy scales, common-set representation, and strong divergence criteria.\\\bottomrule
\end{longtable}
The chronology selects transitions used later; it assigns neither completeness nor priority from a publication date alone. Earlier independent results are cited in their own sections. Its function is to show how an ordered convergence problem becomes a question about selection, representation, or a finite measurement model.

\subsection{Main quantitative thresholds}\label{sub:theoremmap}
Table~\ref{tab:thresholds} records the principal scales. Its references distinguish external inputs from proofs reproduced below. Operator norms refer to the indicated spaces; a critical weak bound is listed separately from strong boundedness.
\begin{longtable}{@{}p{.22\textwidth}p{.29\textwidth}p{.36\textwidth}@{}}
\caption{Thresholds, losses, and sharpness mechanisms}\label{tab:thresholds}\\
\toprule Operator or problem & Range or scale & Source and boundary\\\midrule\endfirsthead
\toprule Operator or problem & Range or scale & Source and boundary\\\midrule\endhead
Individual $S_N$ & $L^1\to L^{1,\infty}$ uniformly & Proposition~\ref{prop:weakindividual}; no supremum in $N$.\\
Full Walsh--Carleson maximum & $L^p$, $p>1$; weak norm $O((p-1)^{-1})$ near one & \cite{Billard,Sjolin,DiPlinio}; full-sequence endpoint in Section~\ref{sub:endpoint-refinement}.\\
Hardy partial sums & $2^{d(n)(1/p-1)}$, $0<p<1$ & \cite{TephPartial,PTW22}; proof and block-atom lower bound in Proposition~\ref{v2:gap}.\\
Full Fej\'er maximum & $H^p\to L^p$, $p>1/2$ & Theorem~\ref{th:fejerhardy}; strong failure at $1/2$, separate weak endpoint.\\
Dyadic Fej\'er maximum & $H^p\to L^p$, $p>0$ & Theorem~\ref{th:dyadicfejer}; no strong terminal-$L^1$ bound.\\
Reversed N\"orlund means, $Q_{2^j}\asymp2^{\alpha j}$ & $p>1/(1+\alpha)$ within $1/2<p\le1$ & \cite{GogNorlund23}; requires the stated $L^\infty$ hypothesis.\\
Greedy approximation & $N^{|1/p-1/2|}e_N(f)_p$ & Theorem~\ref{th:greedyrate}; Dirichlet and independent-character examples give the exact power.\\
Disjoint interval synthesis & $1<p\le2$ & \cite{Osipov16}; proof in Section~\ref{sec:timefreq}, uniform in the interval family.\\
Walsh--Carleson weights & Upper power $\max\{p',2/(p-1)\}$ & \cite{DiPlinioLerner14}; no assertion of sharpness throughout the range.\\
High-order digital integration & $N^{-\alpha}(\log N)^{(s-1)/2}$ & \cite{GSY17}; RMS worst-case error, with fixed $t$, dimension, and positive product weights.\\\bottomrule
\end{longtable}

\section{Characters, orderings, and exact kernels}\label{sec:foundations}
The first question is what is preserved when the characters are reordered or represented by a finite matrix. The answers are exact finite-digit identities. They fix the normalization for later endpoint estimates and distinguish block preservation from preservation of every initial segment.

The algebraic model fixes which manipulations are exact. We first identify the continuous characters and the finite sampling convention, then compute kernels whose localization will be used in the maximal and representation arguments. Binary indices enter both through their highest scale and through the internal pattern of their digits.

\subsection{The compact group and the three orderings}
Let
\[
 \D=\prod_{j\ge0}\Z_2,\qquad x\oplus y=(x_j+y_j\pmod2)_{j\ge0},
\]
with normalized Haar measure $\mu$. Put $r_j(x)=(-1)^{x_j}$. If $n=\sum_{j\ge0}n_j2^j$, the Walsh--Paley character is
\[
 w_n(x)=\prod_{j\ge0}r_j(x)^{n_j}.
\]
Thus $w_nw_k=w_{n\mathbin{\mathrm{xor}}k}$, not generally $w_{n+k}$. We use the cylinders
\[
 I_m(x)=\{y:y_j=x_j\ (0\le j<m)\},\qquad I_m=I_m(0).
\]
Except at the usual double binary expansions, $x$ corresponds to $\sum_jx_j2^{-j-1}\in[0,1)$. It is safer to prove group identities on $\D$ first; dyadic addition is not ordinary addition modulo one.

\begin{proposition}\label{prop:characters}
The functions $w_n$ are exactly the continuous characters of $\D$. They form an orthonormal basis of $L^2(\D)$, and their finite linear combinations are dense in $C(\D)$ and in $L^p(\D)$, $1\le p<\infty$.
\end{proposition}
The finite-cylinder character theorem gives this description: continuity forces a character to be trivial on some cylinder subgroup, and the resulting finite quotient has exactly the characters $w_n$ below that level; see \cite{GES,SWS}.

The classical sources are \cite{Rademacher,Walsh,Paley,Fine,Vilenkin}.

In sequency ordering one arranges finite Walsh patterns by their numbers of sign changes; the finite relation with Paley ordering uses the binary reflected Gray code, together with the bit convention used for sampling. The endpoints of $[0,1)$ are not joined when counting sign changes. The Walsh--Kaczmarz system reverses the active lower digits within each dyadic block. More explicitly, if $n=2^m+\sum_{j<m}n_j2^j$, then
\[
 \kappa_n(x)=r_m(x)\prod_{j<m}r_{m-1-j}(x)^{n_j},\qquad \kappa_0=1.
\]
The block $\{\kappa_n:2^m\le n<2^{m+1}\}$ equals the corresponding Paley block as a set, but its initial segments differ. The historical role of Sneider's work \cite{Sneider48} and the distinct summability theory \cite{Gat98,GN11} should therefore be retained. Likewise, non-equivalence of rearranged Walsh and trigonometric systems in $L^p$, $p\ne2$, is a genuine obstruction to transferring results solely by rearranging bases \cite{HinrichsWenzel}.

\subsection{Dirichlet kernels and binary variation}
Our conventions throughout are
\begin{equation}\label{eq:conventions}\begin{gathered}
 \wh f(n)=\int f w_n\dd,\quad S_Nf=\sum_{n=0}^{N-1}\wh f(n)w_n,
 \quad D_N=\sum_{n=0}^{N-1}w_n,\\
 \quad \sigma_Nf=\frac1N\sum_{k=1}^NS_kf,
 \quad K_N=\frac1N\sum_{k=1}^ND_k.
\end{gathered}\end{equation}
Write $E_m=S_{2^m}$ and $\Delta_m=E_{m+1}-E_m$. Then
\begin{equation}\label{eq:dirichlet}
 D_{2^m}=\prod_{j<m}(1+r_j)=2^m\ind_{I_m},\qquad
 E_mf(x)=2^m\int_{I_m(x)}f\dd.
\end{equation}
For $0\le k\le2^m$,
\[
 D_{2^m+k}=D_{2^m}+r_mD_k.
\]
If $\tau_mx=(x_m,x_{m+1},\ldots)$, then a second exact identity is
\begin{equation}\label{eq:dilate}
 D_{2^mq}(x)=D_{2^m}(x)D_q(\tau_mx),\qquad \norm{D_{2^mq}}_1=\norm{D_q}_1.
\end{equation}
The identity comes from writing each integer below $2^mq$ as $2^ma+b$.

Set $L_N=\norm{D_N}_1$. Integrating the preceding recursion separately over $I_m$ and its complement gives
\begin{equation}\label{eq:lebesgue}
 L_{2^m+k}=1+L_k-\frac{k}{2^m}\quad(0\le k<2^m),\qquad
 L_{2^m-1}=2-2^{1-m}\quad(m\ge1).
\end{equation}
Indeed, on $I_m$, $D_k=k$ and the average over $x_m$ of $|2^m+r_mk|$ is $2^m$; off $I_m$ the first term vanishes. For
\[
 V(n)=n_0+\sum_{j\ge1}|n_j-n_{j-1}|,
\]
classical estimates give $L_n\asymp V(n)$; binary variation and related refinements remain active topics \cite{GES,SWS,AstashkinSemenov16}. In particular, $L_N$ is not bounded, despite $L_{2^m}=1$.

\subsection{Fej\'er kernels: positivity at blocks only}
\begin{proposition}\label{prop:fejerkernel}
For $m\ge1$, $K_{2^m}$ has the following values, determined by the first $m$ coordinates:
\[
 K_{2^m}(x)=
 \begin{cases}
 (2^m+1)/2,&x\in I_m,\\
 2^{j-1},&x_j=1\text{ and }x_i=0\text{ for }i<m,\ i\ne j,\\
 0,&\text{at least two of the first }m\text{ coordinates are }1.
 \end{cases}
\]
Consequently $K_{2^m}\ge0$ and $\norm{K_{2^m}}_1=1$.
\end{proposition}
\begin{proof}
Pair the terms $D_k$ and $D_{2^{m-1}+k}$ and use the Dirichlet recursion. This yields
\[
 2^mK_{2^m}=2^{m-1}(1+r_{m-1})K_{2^{m-1}}+2^{m-1}D_{2^{m-1}}.
\]
Induction gives the three stated cases, with $K_1=1$. Integration, or the constant Fourier coefficient, gives unit mass.
\end{proof}
For non-dyadic indices kernels can change sign: at $x_0=x_1=x_2=1$, $D_1,D_2,D_3,D_4,D_5$ equal $1,0,-1,0,-1$, so $K_5=-1/5$. Uniform $L^1$ boundedness of the full family is still valid; precise kernel estimates are developed in \cite{Toledo18}. This small example prevents an incorrect positivity proof for the full maximal theorem.

\subsection{Exact sequency and the finite-to-infinite identification}\label{sub:v2seq}
Let $\operatorname{seq}(w_n)$ count the internal changes of sign on $[0,1)$; the two endpoints are not joined. The binary convention in Section~\ref{sec:foundations} gives an exact formula, rather than only a qualitative relation with Gray coding.
\begin{proposition}\label{v2:seq}
If $b_j=n_j\oplus n_{j+1}\oplus\cdots$, then
\[
 \operatorname{seq}(w_n)=\sum_{j\ge0}b_j2^j=G^{-1}(n),
 \qquad G(s)=s\mathbin{\mathrm{xor}}\lfloor s/2\rfloor.
\]
The function of sequency $s$ is therefore $w_{G(s)}$.
\end{proposition}
\begin{proof}
For $n=0$ the character is constant and the formula gives zero. For $n>0$, suppose the highest active digit is $m$. The function is constant on $2^{m+1}$ consecutive intervals. A boundary number with exactly $t$ terminal zeros changes the fractional coordinates $x_{m-t},\ldots,x_m$ simultaneously. The character changes sign precisely when $n_{m-t}\oplus\cdots\oplus n_m=1$, and there are $2^{m-t}$ such boundaries. Summing over $t=0,\ldots,m$ gives the formula. The inverse identity is $n_j=b_j\oplus b_{j+1}$.
\end{proof}
\begin{center}
\begin{tabular}{c|rrrrrrrr}\toprule
$n$&0&1&2&3&4&5&6&7\\\midrule
$\operatorname{seq}(w_n)$&0&1&3&2&7&6&4&5\\
$\pi_K(n)$&0&1&2&3&4&6&5&7\\
$G(n)$&0&1&3&2&6&7&5&4\\\bottomrule
\end{tabular}\end{center}
Both permutations preserve the initial sets $[0,2^m)$, but not their intermediate initial segments. For example, at six terms the Paley and Kaczmarz frequency sets are $\{0,1,2,3,4,5\}$ and $\{0,1,2,3,4,6\}$. Bit reversal of an entire variable-length integer, without keeping the leading bit fixed, does not define the Kaczmarz permutation.

For the finite Sylvester matrix, use ordinary binary integer indices. Then
\[
 H_{2^m}[i,j]=(-1)^{\sum_{k<m}i_kj_k},\qquad
 w_i(j/2^m)=H_{2^m}[i,\operatorname{rev}_m(j)].
\]
The reversal appears because the $k$th fractional digit of $j/2^m$ is the integer digit $j_{m-1-k}$. On the abstract finite group, where $x$ is already a coordinate vector, the reversal is absent. This convention reconciles finite algorithms with continuous character calculations without changing either normalization.

\subsection{One algebra, several non-equivalent operators}\label{sub:v2geometry}
Numerical frequency and dyadic scale satisfy $2^{u(n)}\le n<2^{u(n)+1}$ for $n>0$. Their positive powers therefore define equivalent Hilbert graph norms. Hamming weight only satisfies $h(n)\le n$; there is no converse uniform bound on $n$ in terms of $h(n)$. Binary variation is different again: long strings of consecutive ones can have small variation even at very large frequency. These are exact comparisons, not a claim that the three scales have no inequalities between them.

The spectral counting function of the numerical generator is $\#\{n:n\le\lambda\}=\lfloor\lambda\rfloor+1$ for $\lambda\ge0$, and each eigenvalue has multiplicity one. The scale symbol $2^{u(n)}$ repeats over an entire dyadic block. The Boolean number operator has infinitely many frequencies of Hamming weight one on the infinite product group. Consequently its fixed positive-time kernel is not the same trace-class object as \eqref{eq:semigroup}; one cannot transfer a trace estimate from the numerical symbol by identifying the character bases.

\paragraph{Comparison and research boundary.}
The finite identities fix the comparison data for the next section.

\begin{figure}[tb]
\centering\includegraphics[width=\textwidth]{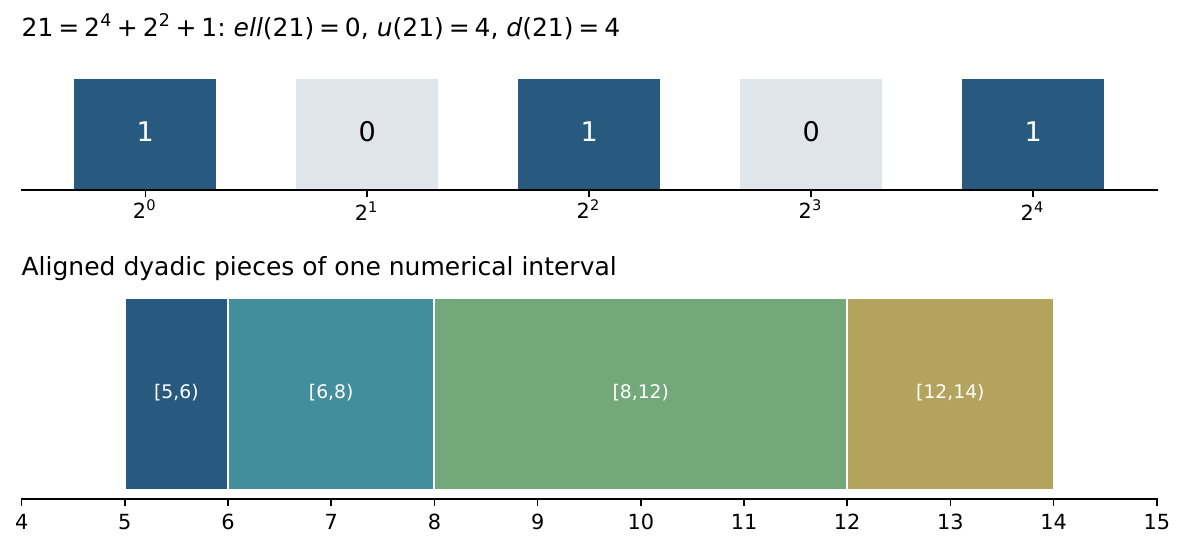}
\caption{Binary data for $21=2^4+2^2+1$ and a disjoint dyadically aligned decomposition of $[5,14)$. The active-bit gap and the lengths of constant binary strings encode different information. The interval pieces illustrate alignment; endpoint-XOR assembly is proved separately in Section~\ref{sec:timefreq}.}\label{fig:binary-geometry}
\end{figure}

\begin{figure}[tb]
\centering\includegraphics[width=\textwidth]{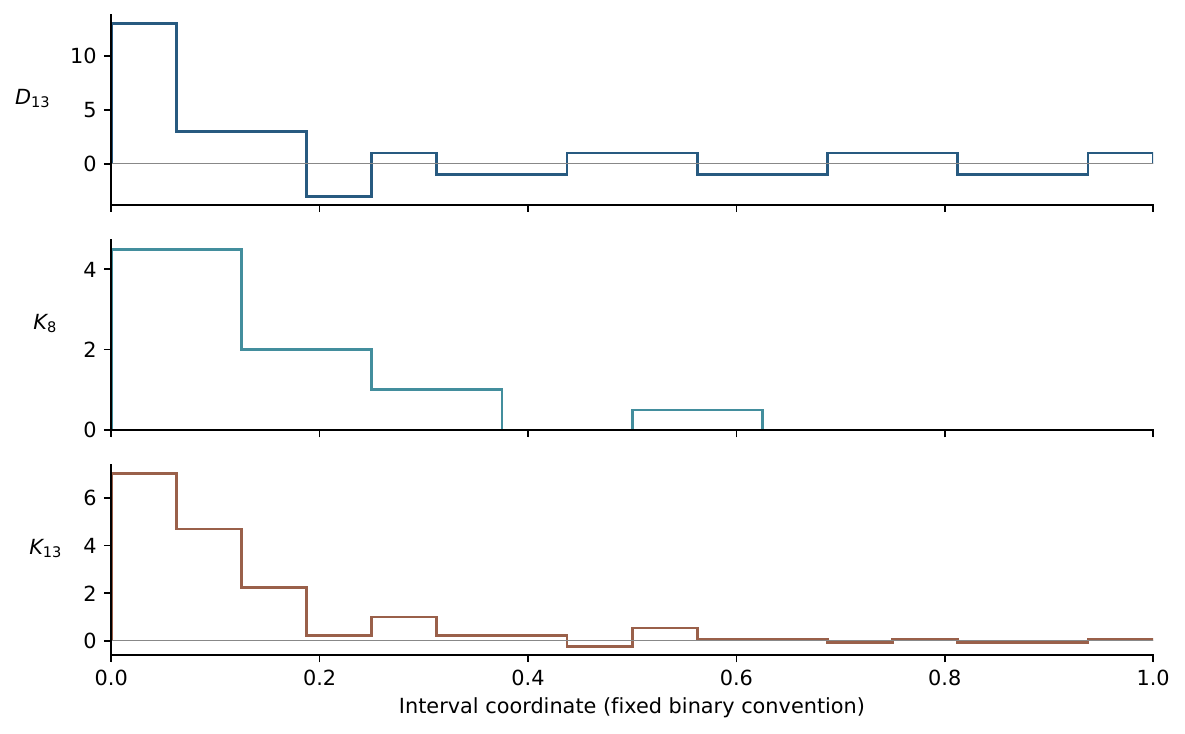}
\caption{Exact step values in the interval realization of the binary group, with positive-index Fej\'er normalization. $D_{13}$ oscillates; $K_8$ is nonnegative, whereas $K_{13}$ has negative values. The plots display cylinder values, not smooth interpolation.}\label{fig:kernel-geometry}
\end{figure}

Figures~\ref{fig:binary-geometry} and \ref{fig:kernel-geometry} display exact finite instances of the preceding identities.
 Preserving a complete dyadic block does not preserve its internal ordering or its averaged kernels, and a matrix implementation must retain the same bit convention before a rate can be transferred.

\section{Ordered convergence, divergence, and subsequences}\label{sec:convergence}
Here the function and Fourier order are fixed. We compare convergence of the complete sequence with individual or selected partial sums. The key boundary is between an estimate uniform in the index and an estimate for its supremum. Binary-gap examples quantify a cost that disappears at martingale indices.

Ordered convergence is a question about a moving endpoint of the spectrum. A fixed endpoint is a modulated martingale transform, giving uniform individual estimates. Controlling all endpoints at once requires a different argument. The section separates these mechanisms before comparing selected indices and averages of selected sums.

\subsection{An individual partial sum is not a maximal theorem}
The identity
\begin{equation}\label{eq:transform}
 S_Nf=w_N\sum_{j\ge0}n_j\Delta_j(fw_N)
\end{equation}
is an efficient explanation of the uniform weak estimate for individual sums. To verify it, observe that multiplication by $w_N$ replaces a frequency $k$ by $k\mathbin{\mathrm{xor}}N$. The condition $k<N$ means that, at the highest differing digit, $N$ has a one and $k$ a zero. The resulting XOR frequency belongs to the corresponding complete dyadic block.

\begin{proposition}\label{prop:weakindividual}
There is an absolute constant $C$ such that, for every $f\in L^1$, $N\ge1$, and $\lambda>0$,
\[
 \mu\{|S_Nf|>\lambda\}\le C\norm f_1/\lambda.
\]
Consequently $S_Nf\to f$ in measure for every $f\in L^1$.
\end{proposition}
\begin{proof}
It suffices by \eqref{eq:transform} to prove the estimate for a finite martingale transform with coefficients zero or one. If $\lambda>\norm f_1$, select maximal dyadic cylinders on which the average of $|f|$ exceeds $\lambda$. Dyadic differentiation shows that outside their union $|f|\le\lambda$ almost everywhere. Replace $f$ on each selected cylinder by its average. The resulting good part $g$ satisfies $\norm g_\infty\le2\lambda$ and $\norm g_1\le\norm f_1$. The bad part has mean zero on each cylinder, so its martingale transform is supported in the union. Orthogonality gives $\norm{Tg}_2^2\le\norm g_2^2\le2\lambda\norm f_1$. Chebyshev's inequality and the measure bound for the union complete the proof. If $\lambda\le\norm f_1$, the probability bound is automatic after increasing $C$.

For convergence in measure, approximate $f$ by a Walsh polynomial $P$. For all sufficiently large $N$, $S_NP=P$. Apply the uniform weak estimate to $f-P$ and also use $\mu\{|f-P|>\lambda\}\le\norm{f-P}_1/\lambda$. Let the approximation error tend to zero.
\end{proof}
The supremum $S^*f=\sup_N|S_Nf|$ cannot be inserted on the left. Uniform weak bounds on separate operators contain no control on a union of infinitely many exceptional sets.

\begin{theorem}[Classical pointwise convergence]\label{th:carleson}
For $1<p<\infty$, $S^*$ is bounded on $L^p(\D)$, and $S_Nf\to f$ almost everywhere for every $f\in L^p(\D)$.
\end{theorem}
Billard's $L^2$ theorem \cite{Billard} and Sj\"olin's work \cite{Sjolin} are essential historical sources. Sj\"olin also proved almost everywhere convergence under the stronger-than-$L\log L$ endpoint integrability condition
\[
 \int |f|\log^+|f|\log^+\log^+|f|\dd<\infty,
\]
with harmless modifications of the logarithms near zero. This should not be reduced in historical attribution to the statement for $p>1$. The Vilenkin treatment \cite{PSSTW22} includes the Walsh case. Modern explanations of the time--frequency mechanism appear in \cite{Thiele06,DemeterGuide16}; vector-valued extensions require geometric assumptions on the range space \cite{HytLacey,HLP}.

At $L^1$, divergence results belong to a different level from growth of operator norms. Stein's principle \cite{Stein}, Schipp's everywhere-divergent constructions \cite{Schipp69}, and Bochkarev's work for Walsh and multiplicative systems \cite{Bochkarev} must be read with their respective conclusions. Hunt--Taibleson's local-field framework \cite{HuntTaibleson71} is another part of the classical pointwise background.

\subsection{\texorpdfstring{Near-$L^1$}{Near-L1} convergence: the historical bound and its refinement}\label{sub:endpoint-refinement}
The integrability condition quoted above is a historical sufficient condition, not the modern frontier. Sj\"olin--Soria \cite{SjolinSoria03} extended Antonov's argument to kernels without the smoothness needed in earlier formulations and obtained almost everywhere Walsh convergence in $L\log L\log\log\log L$. With positive logarithmic regularizations, one may express the integrability class by
\[
 \int_\D\Phi_A(|f|)\dd<\infty,\qquad
 \Phi_A(t)=t\log(e+t)\log\log\log(e^{e^e}+t).
\]
An equivalent Young function may be used to define its Luxemburg norm. This is weaker than $t\log t\log\log t$ growth at infinity. The number and nesting of the logarithms are part of the theorem.

Di Plinio \cite{DiPlinio} established the quantitative weak estimate
\[
 \|S^*f\|_{L^{p,\infty}}\le \frac{C}{p-1}\|f\|_p,
 \qquad 1<p\le2.
\]
Appropriate extrapolation yields the Antonov-class endpoint; it does not give weak $(1,1)$ boundedness by setting $p=1$. The Fourier theorem of Di Plinio--Fragkos \cite{DiPlinioFragkos25} supplies a further quantitative multi-frequency method, discussed in Section~\ref{sub:quartile}.

\paragraph{The Orlicz gap and its quantifiers.}
Bochkarev's Theorem~1 in \cite{Bochkarev} gives an everywhere divergent Walsh--Fourier series of a function in $\Phi(L)$ whenever $\Phi(t)=t\psi(t)$, $\psi$ is nondecreasing, $\psi(0)=1$, and $\psi(t)=o(\sqrt{\log t})$ at infinity. Thus a sufficient convergence class and a divergence construction must be displayed together; neither determines the missing boundary. The word ``everywhere'' is stronger than the failure of convergence on a positive-measure set. The admissible growth assumptions in the source are retained when using its Orlicz notation.

\begin{longtable}{@{}p{.19\textwidth}p{.31\textwidth}p{.37\textwidth}@{}}
\caption{Endpoint classes: order, quantifiers, and what is established}\label{tab:orlicz-gap}\\
\toprule Model & Established assertion & Boundary and source\\\midrule\endfirsthead
\toprule Model & Established assertion & Boundary and source\\\midrule\endhead
Full Walsh order & Every $f\in L\log L\log_3 L$ has a.e. convergence. & Sufficient class, \cite{SjolinSoria03}; the $L\log L$ endpoint remains distinct.\\
Full Walsh order & For each admissible $\psi=o(\sqrt{\log})$, some $f\in(t\psi(t))(L)$ diverges everywhere. & Lower-side obstruction, \cite{Bochkarev}, Theorem~1; not divergence for every input.\\
Each prescribed lacunary Walsh sequence & A.e. convergence for every $f\in L\log_2 L\log_4 L$. & Sufficient class, \cite{DiPlinioLac}; constants and the chosen sequence enter.\\
Some lacunary Walsh sequences & Nested binary spectra give everywhere divergence in admissible classes $\Phi(t)=o(t\log_2t)$. & \cite{GogOniani21}, Theorem~2.6 and Remark~2.7; not a theorem for every lacunary sequence.\\
Lacunary trigonometric order & Lorentz-space classification under specified fundamental-function hypotheses. & \cite{LieLac17}; it is not automatically a Walsh or unrestricted Orlicz classification.\\\bottomrule
\end{longtable}
Here $\log_k$ means the $k$-fold logarithm with positive regularization. Full-sequence convergence and convergence along each lacunary sequence are different questions. The table records the cited conclusions and does not claim a complete classification of all rearrangement-invariant spaces.

\subsection{Norm convergence and the geometry of a selected sequence}
For convolution operators, testing against a translate of the sign of $D_N$ gives
\begin{equation}\label{eq:endpointnorm}
 \norm{S_N}_{C\to C}=\norm{S_N}_{L^1\to L^1}=L_N.
\end{equation}
The first sign is a cylinder function and hence continuous on $\D$; the second equality also follows by convolution with an approximate identity. For a sequence $n_k\to\infty$, polynomial density and uniform boundedness prove the exact criterion
\[
 S_{n_k}f\to f\text{ in }C\text{ for all }f\in C
 \quad\Longleftrightarrow\quad \sup_kL_{n_k}<\infty,
\]
and the analogous statement in $L^1$. If the supremum is infinite, the functions with unbounded selected norms form a residual set. This is a norm statement, not everywhere divergence at a prescribed set of points.

For $n>0$, put
\[
 u(n)=\max\{j:n_j=1\},\quad \ell(n)=\min\{j:n_j=1\},\quad d(n)=u(n)-\ell(n).
\]
In martingale Hardy spaces below $1$, the gap $d(n)$, rather than $V(n)$ alone, controls a different phenomenon. A basic sharp estimate is
\begin{equation}\label{eq:bitgap}
 \norm{S_nF}_p\le C_p2^{d(n)(1/p-1)}\norm F_{H^p},\qquad 0<p<1.
\end{equation}
Tephnadze's results \cite{TephPartial}, together with \cite{PTW22}, develop boundedness, approximation, and necessity statements at this scale. An $H^p\to L^p$ atom calculation is not by itself a proof of every $H^p\to H^p$ criterion in that literature.

The unit block atom in Proposition~\ref{v2:gap} realizes this loss at $n=2^M+1$, although $V(2^M+1)$ stays bounded. Thus binary variation and the Hardy-space gap measure different obstructions.

For a prescribed lacunary sequence $n_{j+1}/n_j\ge\lambda>1$, Do--Lacey \cite{DoLaceyLac} proved almost everywhere convergence in $L\log_2 L\log_3 L$. Di Plinio \cite{DiPlinioLac} improved this to $L\log_2 L\log_4 L$. Here $\log_k t$ denotes the $k$-fold logarithm, positively regularized near zero; the associated Young functions have the stated growth at infinity. These classes are much closer to $L^1$ than the full-sequence Antonov class.

The improvement comes from the lacunary weak-$L^p$ estimate $\|\sup_j|S_{n_j}f|\|_{p,\infty}\lesssim_\lambda\log(e+p')\|f\|_p$, $1<p\le2$, and extrapolation. This quantitative gain replaces the full-sequence $(p-1)^{-1}$ loss. Do--Lacey use multi-frequency decomposition and lacunary coefficient control; Di Plinio's Walsh argument avoids the Antonov approximation step. The normalization of the partial sum in \cite{DiPlinioLac} is inclusive, $\sum_{k=0}^n\widehat f(k)w_k=S_{n+1}f$ here. Passing to $(n_j-1)$, which is lacunary after discarding finitely many indices, gives the stated assertion for $S_{n_j}$.

Optimality must retain the model and function-space scale. Lie \cite{LieLac17} obtains a trigonometric Lorentz-space classification under assumptions on the fundamental function; it is not a Walsh classification. Goginava--Oniani \cite{GogOniani21} give Walsh divergence criteria tied to binary spectra and formulate the optimal universal lacunary class as an open problem in their account. Their result does not settle every fixed lacunary sequence. Selected-sum averaging \cite{Konyagin93,GatSub10} and strong summability \cite{Simon00,GatGog20} have different conclusions.

\subsection{The atomic proof of the binary-gap estimate}\label{sub:v2gap}
We supply the mechanism behind \eqref{eq:bitgap}; it is an estimate into $L^p$, distinct from stronger Hardy-norm assertions in \cite{TephPartial}.
\begin{proposition}\label{v2:gap}
For $0<p<1$ and $n\ge1$,
\[
 \norm{S_nF}_p\le C_p2^{d(n)(1/p-1)}\norm F_{H^p}.
\]
The order is attained at $n=2^M+1$ on suitable unit atoms.
\end{proposition}
\begin{proof}
Let $a$ be an atom on $I_M$. There is no contribution if $n\le2^M$. For $n>2^M$, the integral over $I_M$ is at most one, because
\[
 \int_{I_M}|S_na|^p\dd
 \le2^{-M(1-p/2)}\norm{S_na}_2^p
 \le2^{-M(1-p/2)}\norm a_2^p\le1.
\]
The exact kernel expansion
\[
 D_n=w_n\sum_jn_jr_jD_{2^j}
\]
follows by the same highest-differing-digit partition as \eqref{eq:transform}. For $x\notin I_M$, its terms with $j\ge M$ vanish after translation by $t\in I_M$. All remaining lower-coordinate factors are unchanged by $t$. Hence
\[
 S_na(x)=w_n(x)\Big(\int_{I_M}a(t)w_n(t)\dd(t)\Big)
                 \sum_{j<M}n_jr_j(x)D_{2^j}(x).
\]
The integral factor is at most $2^{M(1/p-1)}$. On a shell with first one at position $k<M$, the last sum vanishes for $k<\ell(n)$ and otherwise has modulus below $2^{k+1}$. Therefore
\[
 \int_{\D\setminus I_M}|S_na|^p\dd
 \lesssim_p2^{M(1-p)}\sum_{k=\ell(n)}^{M-1}2^{k(p-1)}
 \lesssim_p2^{(M-\ell(n))(1-p)}
 \le C_p2^{d(n)(1-p)}.
\]
If $\ell(n)\ge M$ this part is zero. To pass from atoms to a general martingale, fix $n$ and choose $m$ with $n<2^m$. Since $F_m$ is constant on cylinders of measure $2^{-m}$,
\[
 |\widehat F(n)|\le\|F_m\|_1
 \le2^{m(1/p-1)}\|F_m\|_p
 \le2^{m(1/p-1)}\|F\|_{H^p}.
\]
Thus each fixed coefficient is continuous on $H^p$. Atomic expansions converging in $H^p$ may consequently be passed through the finite operator $S_n$. Finite atomic sums satisfy the estimate by $p$-subadditivity, and passage to the limit proves the bound. The atom $2^{M(1/p-1)}r_MD_{2^M}$ has unit Hardy norm and a single nonzero partial coefficient at $n=2^M+1$, proving the claimed lower order.
\end{proof}
\begin{corollary}\label{v2:gapconvergence}
If $n_k\to\infty$ and $\sup_kd(n_k)<\infty$, then $S_{n_k}F$ converges in $L^p$ to the terminal martingale limit for every $F\in H^p$, $0<p<1$.
\end{corollary}
\begin{proof}
Approximate $F$ in $H^p$ by finite martingales using atomic decomposition and conditional expectations of finitely many bounded atoms. Partial sums stabilize on these polynomials. The proposition controls the remaining error uniformly in $k$. The martingale convergence theorem for $H^p$ gives an almost everywhere finite terminal limit \cite{Weisz94,PTW22}, and dominated convergence applies to $|F_m-F_\infty|^p\le2^p(F^*)^p$. This does not identify $F$ uniquely from $F_\infty$.
\end{proof}

\subsection{Selected sums and averages of selected sums}\label{sub:v2selected}
Konyagin \cite{Konyagin93} proved the existence of a sequence $(n_j)$ for which $\sup_jL_{n_j}=\infty$ and $S_{n_j}f\to f$ almost everywhere for every $f\in L^1$. Since $\|S_{n_j}\|_{L^1\to L^1}=L_{n_j}$ is unbounded, the uniform boundedness principle gives an $f\in L^1$ with $\sup_j\|S_{n_j}f\|_1=\infty$. The same sequence can thus have a universal almost everywhere property and fail a universal norm property. Its arithmetic construction, rather than only the fact that a subsequence exists, is the content of the external result.

The selection rule must be retained when comparing averaged sums. G\'at's published trigonometric divergence theorem \cite{GatArithmetic24} concerns Ces\`aro means of a selected sequence of trigonometric partial sums. It cannot be transferred to Walsh sums by replacing the basis name: the kernels and the arithmetic of the selected indices change. The positive Walsh statements below therefore retain their explicit index hypotheses.

G\'at's selected-sum theorems \cite{GatSub10,GatSurvey16} provide another layer. If $a(k+1)/a(k)\ge q>1$, then for every $f\in L^1$,
\[
 \frac1N\sum_{k=1}^NS_{a(k)}f\longrightarrow f\quad\text{almost everywhere}.
\]
For an increasing convex sequence of integer indices tending to infinity, the logarithmic averages
\[
 \frac1{H_N}\sum_{k=1}^N\frac{S_{a(k)}f}{k}\longrightarrow f,
 \qquad H_N=\sum_{k=1}^N\frac1k,
\]
converge almost everywhere. Normalization by $H_N$ is asymptotically equivalent to the $\log N$ convention and reproduces constants exactly. Neither assertion claims convergence of the unaveraged $S_{a(k)}f$. The average and the unaveraged sequence remain different approximation processes.

\begin{proposition}[Uniform interior norm bounds]\label{v2:partialnorms}
For $1<p<\infty$, $\sup_N\norm{S_N}_{L^p\to L^p}<\infty$, and $S_Nf\to f$ in $L^p$. This norm result follows from martingale transforms without proving the full Walsh--Carleson maximal theorem.
\end{proposition}
\begin{proof}
The proof of Proposition~\ref{prop:weakindividual} gives uniform weak $(1,1)$ bounds for the finite transforms in \eqref{eq:transform}; orthogonality gives their $L^2$ norm at most one. Interpolation proves the uniform range $1<p\le2$. Each $S_N$ is self-adjoint in the $L^2$ pairing, so duality gives $p\ge2$. Polynomial density then proves convergence in norm. Taking a supremum over all $N$ before integration is a different estimate, supplied by Theorem~\ref{th:carleson}.
\end{proof}

\paragraph{Comparison and research boundary.}
The distinction between an individual partial sum and its maximum remains decisive at the endpoint. The specified full-sequence endpoint question in Section~\ref{sec:boundaries} concerns a stronger distributional conclusion than a collection of individual weak bounds. Selecting indices changes the geometry and must be recorded as a separate operation.

\section{Summation methods and matrix criteria}\label{sec:summation}
Averaging changes the operator rather than the target function. The question is which restrictions on averaging rows give a regular approximation process, and when they control a maximum or an almost-everywhere limit. Kernel identities give elementary sufficient mechanisms; the matrix criteria below have additional hypotheses that must survive any comparison.

Averaging replaces the Dirichlet kernel by a row-dependent combination of kernels. The central issues are reproduction of fixed frequencies, uniform kernel mass, and control of a maximal family. Abel domination handles one class of rows; decreasing reversed weights and dispersed nonmonotone rows require the cited structural theorems.

\subsection{Fej\'er convergence and regular approximation}
\begin{theorem}[Schipp's maximal theorem]\label{th:fejerweak}
The full Fej\'er maximal operator $\sigma^*f=\sup_N|\sigma_Nf|$ is of weak type $(1,1)$. Hence $\sigma_Nf\to f$ almost everywhere for every $f\in L^1(\D)$.
\end{theorem}
This is the weak endpoint theorem \cite{Schipp75}; it is not an $L^1\to L^1$ estimate. The Hardy-space endpoint belongs to a different domain \cite{Fujii,LadhawalaPankratz}. Walsh--Kaczmarz summability has its own proofs and maximal estimates \cite{Gat98,GN11}.

The norm convergence part follows from a useful elementary principle. If $B_N$ are convolution kernels with $\sup_N\norm{B_N}_1<\infty$ and $\wh B_N(k)\to1$ for every fixed $k$, then $f*B_N\to f$ in $L^p$, $1\le p<\infty$, and in $C(\D)$. Prove this first for polynomials, then use uniform boundedness and density. Almost everywhere convergence additionally needs a suitable maximal estimate; regularity of the matrix alone does not provide one.

\subsection{N\"orlund means: direction and critical growth}
Let $q_j\ge0$, $q_0>0$, and $Q_N=\sum_{j=0}^{N-1}q_j$. We distinguish
\begin{equation}\label{eq:norlund}
 t_Nf=\frac1{Q_N}\sum_{k=1}^Nq_{N-k}S_kf,
 \qquad T_Nf=\frac1{Q_N}\sum_{k=1}^Nq_{k-1}S_kf.
\end{equation}
The first is the reversed, or usual N\"orlund, mean. Exchanging it with the direct mean changes which monotonicity assumption is useful.

Indeed, for $A_Nf=\sum_{k=1}^Na_{N,k}S_kf$, set $a_{N,N+1}=0$ by convention. Abel summation gives
\begin{equation}\label{eq:matrixabel}
 A_Nf=\sum_{k=1}^Nk(a_{N,k}-a_{N,k+1})\sigma_kf,
 \qquad |A_Nf|\le B_N\sigma^*f,
 \quad B_N=\sum_{k=1}^Nk|a_{N,k}-a_{N,k+1}|.
\end{equation}
Uniform control requires $\sup_NB_N<\infty$. If a row is nonnegative, nonincreasing, and has sum one, then $B_N=1$. Consequently nondecreasing $q_j$ give domination for $t_N$, and nonincreasing $q_j$ give domination for $T_N$. With the usual regularity conditions, these yield norm and almost everywhere convergence. This argument cannot settle the difficult reversed means with decreasing coefficients.

\begin{theorem}[Goginava's criterion \cite{GogNorlund23}, Theorem 3]\label{th:norlundcriterion}
Suppose $q_j>0$ is nonincreasing, $Q_N\to\infty$, $1/2<p\le1$, and $t^*=\sup_N|t_N|$ is bounded on $L^\infty$. Then $t^*:H^p\to L^p$ is bounded if and only if
\begin{equation}\label{eq:gammap}
 \sup_{M\ge1}\frac{2^{M(1-p)}}{Q_{2^M}^{p}}
 \sum_{j=1}^MQ_{2^j}^{p}2^{j(p-1)}<\infty.
\end{equation}
\end{theorem}
The $L^\infty$ hypothesis is part of the theorem. It cannot silently be removed on the ground that each matrix row is normalized.

The criterion has an instructive model. If $Q_{2^j}\asymp2^{\alpha j}$, then the summand in \eqref{eq:gammap} is comparable to $2^{j[p(1+\alpha)-1]}$. Writing $b=p(1+\alpha)-1$, the outer factor is comparable to $2^{-Mb}$. For $b>0$ the geometric sum cancels it; for $b=0$ the expression grows like $M$; for $b<0$ it grows exponentially. Thus the strong threshold, under the stated hypotheses, is $p>1/(1+\alpha)$. Logarithmic factors can decide finer boundaries and should be substituted into the actual sum. In particular, reversed logarithmic means $q_j=(j+1)^{-1}$ behave very differently from the direct means, which are covered by \eqref{eq:matrixabel}.

The analysis of individual N\"orlund means and their norms \cite{GogNagy23}, subsequences \cite{LogNorlund22,SubNorlund23}, and critical weak estimates \cite{GMM25} gives a more detailed map than one universal ``N\"orlund theorem.'' Recent norm and almost everywhere convergence results for Vilenkin subsequences appear in \cite{Norlund26}. The weak critical result for power growth is separate from the strict strong range in \eqref{eq:gammap}; interpolation does not recover a missing strong endpoint.

\subsection{Nonmonotone rows and structural convergence theorems}
Matrix methods with no monotonicity of individual rows are a significant modern development \cite{GoginavaNagy22,GogMuk24,BlahotaGat25}. One representative result is as follows.
\begin{theorem}[Blahota--G\'at \cite{BlahotaGat25}, sufficient norm condition]\label{th:bg}
Let $d_j\to\infty$, $t_{k,d_j}\ge0$, and $\sum_{k=1}^{d_j}t_{k,d_j}=1$. If
\[
 \sum_{k=1}^{d_j}t_{k,d_j}^2\le C/d_j,
\]
then the matrix means $\sum_{k=1}^{d_j}t_{k,d_j}S_kf$ converge to $f$ in $L^1$ for every $f\in L^1$, and uniformly for every $f\in C(\D)$.
\end{theorem}
This is a sufficient condition, not an equivalence. For the almost everywhere conclusion in the corresponding theorem one additionally imposes lacunary growth $d_{j+1}\ge\lambda d_j$, $\lambda>1$. The negative results in \cite{BlahotaGat25} explain why arbitrary concentration of row mass or removal of structural conditions cannot be ignored. They are not consequences of the elementary Abel bound.

The approach of Goginava--Mukhamedov \cite{GogMuk24} connects uniform $H^1\to L^1$ bounds with pointwise convergence for a specified family of Walsh matrix operators. The structure of that family is indispensable: uniform bounds on a collection of arbitrary operators do not imply almost everywhere convergence. The classical spatial $A_p$ theory in Section~\ref{sec:timefreq} concerns a different, explicitly specified operator and weight class.

\subsection{Uniform kernel bounds and the operator convergence principle}\label{sub:v2kernelbounds}
\begin{proposition}\label{v2:kernelL1}
With the positive-index convention \eqref{eq:conventions}, $\sup_N\norm{K_N}_1\le2$.
\end{proposition}
\begin{proof}
For $N=2^m+l$, $0\le l<2^m$, sum the Dirichlet recursion to obtain
\[
 NK_N=2^mK_{2^m}+lD_{2^m}+r_mlK_l.
\]
Dyadic kernels have norm one by Proposition~\ref{prop:fejerkernel}. Induction on $N$ gives
$N\norm{K_N}_1\le2^m+l+2l\le2(2^m+l)$.
The case $l=0$ is already known; $K_1=1$ starts the induction.
\end{proof}
\begin{proposition}\label{v2:regularconvolution}
Suppose $\sup_N\norm{B_N}_1<\infty$ and $\wh B_N(k)\to1$ at every fixed frequency. Then $f*B_N\to f$ in $L^p$, $1\le p<\infty$, and uniformly in $C(\D)$. If additionally $B^*f=\sup_N|f*B_N|$ has weak type $(1,1)$, then convergence is almost everywhere for every $L^1$ function.
\end{proposition}
Polynomial density and the uniform convolution norm give the norm assertions. The weak maximal approximation argument, applied to $f-P$ and then letting $\|f-P\|_1\to0$, gives the almost everywhere assertion; see \cite{GES,SWS,Weisz02}.
For general matrix rows, \eqref{eq:matrixabel} proves a uniform maximal estimate whenever $\sup_NB_N<\infty$. If their sums tend to one and each fixed coefficient tends to zero, their action on a polynomial converges to that polynomial: the finitely many exceptional low indices vanish. Density and the weak maximum then finish the argument. The convention $a_{N,N+1}=0$ merely expresses the boundary term in Abel summation.

The nonmonotone condition in Theorem~\ref{th:bg} has sharpness information beyond norm convergence. Theorems 3 and 5 of \cite{BlahotaGat25} permit an arbitrarily slowly diverging relaxation factor in the squared row mass and construct divergent examples, in norm or almost everywhere in the respective statements. The factor is evaluated at the paper's specified index, and the almost everywhere example uses dyadic row lengths. Thus this is an optimality statement for a universal sufficient scale; it is not a claim that every row failing the bound is divergent.

More explicitly, for every positive sequence $\lambda_n\to\infty$, the norm-divergence statement allows normalized nonnegative rows at indices $m_j\uparrow\infty$ with
\[
 \sum_{k\le m_j}t_{k,m_j}^2=O(\lambda_{m_j}/m_j)
\]
and an $f\in L^1$ whose mean norms tend to infinity. In the almost everywhere counterpart, for nondecreasing $\lambda_n\to\infty$ one can choose dyadic lengths $2^{N_j}$ and rows satisfying
\[
 \sum_{k\le2^{N_j}}t_{k,2^{N_j}}^2=O(\lambda_{N_j}/2^{N_j})
\]
with an integrable function whose mean moduli tend to infinity almost everywhere. The different placement of $\lambda_{m_j}$ and $\lambda_{N_j}$ belongs to the statements. These examples complement the positive theorem without asserting necessity of its condition for an individual matrix.

\subsection{Fixed and varying Ces\`aro orders}\label{sub:v2cesaro}
For $\alpha>0$, put $A_0^\alpha=1$ and
$A_n^\alpha=\prod_{j=1}^n(1+\alpha/j)$. The normalized positive-index mean is
\[
 \sigma_n^\alpha f=\frac1{A_{n-1}^\alpha}
       \sum_{k=1}^nA_{n-k}^{\alpha-1}S_kf,
 \qquad\sum_{j=0}^{n-1}A_j^{\alpha-1}=A_{n-1}^\alpha.
\]
For $0<\alpha\le1$, the reversed coefficients are decreasing and their accumulated sums have power growth $n^\alpha$. Their uniformly bounded kernel norms supply the $L^\infty$ hypothesis; Theorem~\ref{th:norlundcriterion} then gives the strong maximal range $p>1/(1+\alpha)$ within its stated domain. The critical weak estimate and sharper weighted statements require separate theorems \cite{WeiszCesaro01,BTCesaro14}. Fine's earlier summability theorem \cite{FineCesaro55} belongs to the classical foundation and should be distinguished from later maximal and Hardy-space results.

For $0<\alpha<1$ and $0<p<1/(1+\alpha)$, Blahota--Tephnadze--Toledo \cite{BTTCesaro15}, Theorems 1 and 2, prove
\[
 \norm{\sup_{n\ge1}\frac{|\sigma_n^\alpha F|}{(n+1)^{1/p-1-\alpha}}}_p
 \le C_{\alpha,p}\norm F_{H^p},
\]
and the strong summability estimate
\[
 \sum_{m\ge1}\frac{\norm{\sigma_m^\alpha F}_{H^p}^p}{m^{2-(1+\alpha)p}}
 \le C_{\alpha,p}\norm F_{H^p}^p.
\]
If a nondecreasing $\varphi\ge1$ satisfies $n^{1/p-1-\alpha}/\varphi(n)\to\infty$, the maximum normalized by $\varphi(n)$ is not bounded $H^p\to L^p$. The mean normalization in that source differs, if $A_n^\alpha$ is used, by $A_n^\alpha/A_{n-1}^\alpha=(n+\alpha)/n$; its boundedness preserves these estimates after changing constants. The strong summability formula controls a weighted sum of norms and is not boundedness of the unweighted full maximum below its critical index.

When $\alpha=\alpha_n\to0$, even a positive order at every index no longer supplies a uniform fixed-order theorem. The trigonometric generalized-order work of Akhobadze \cite{Akhobadze07} is methodological background; the Walsh statements are due to the corresponding Walsh papers \cite{AbuJoudehGat18,GatGogVary22}. In the latter, write
\[
 V_\alpha(n)=\sum_{i\ge0}|n_i-n_{i+1}|2^{i\alpha_n},\quad
 \mathcal V_K=\{n:V_\alpha(n)\le K n^{\alpha_n}\}.
\]
For $0<\alpha_n<1$, Theorem 1.1 of \cite{GatGogVary22} gives almost everywhere convergence for every $L^1$ input along $n\in\mathcal V_K$, for fixed finite $K$. Theorem 1.8 gives an unbounded almost everywhere divergent example if, eventually,
\[
 \frac{c_1}{\log_2n}\le\alpha_n\le
 \frac{c_0\log_2\log_2n}{\log_2n},\qquad c_1>0,\quad0<c_0<\tfrac12.
\]
Here $V_\alpha$ is the source's weighted boundary-digit expression, not our $V(n)$ with its additional lowest-digit term. Thus the order sequence and binary geometry interact, rather than only the scalar limit $\alpha_n\to0$ deciding convergence.

The later negative theorem is stronger. Goginava \cite{GogVary24} proves that for every positive order sequence $\alpha_n\to0$, however slowly, there exists $f\in L^1$ whose full sequence of $(C,\alpha_n)$ means diverges almost everywhere. This is the quantified assertion of the published paper. Thus there is no decay rate of the scalar order alone that guarantees full-sequence a.e. convergence for every integrable input. The positive conclusion on $\mathcal V_K$ concerns a restricted index set and remains compatible with this theorem. The harmless normalization factor $(n+\alpha_n)/n\to1$ preserves convergence or its failure between the two Ces\`aro conventions. No growth rate of the divergence is asserted from the abstract-level statement.

\begin{longtable}{@{}p{.20\textwidth}p{.31\textwidth}p{.36\textwidth}@{}}
\caption{Summation map: sufficient ranges, exact criteria, and obstructions}\label{tab:summation-map}\\
\toprule Method & Established range or hypothesis & Logical status and location\\\midrule\endfirsthead
\toprule Method & Established range or hypothesis & Logical status and location\\\midrule\endhead
Full Fej\'er maximum & $H^p\to L^p$ for $p>1/2$; critical weak estimate is separate. & Sharp strong threshold, Theorem~\ref{th:fejerhardy}.\\
Fixed Ces\`aro order $0<\alpha\le1$ & Strong maximal range $p>1/(1+\alpha)$ within the stated Hardy domain. & Under the kernel and coefficient hypotheses above; critical weak input is separately cited.\\
Reversed N\"orlund & Positive nonincreasing coefficients and the $L^\infty$ maximal hypothesis. & Exact criterion \eqref{eq:gammap}, Theorem~\ref{th:norlundcriterion}; not a statement for every normalized row.\\
General nonnegative matrix rows & Unit row mass and $\sum_k t_{k,d_j}^2\lesssim d_j^{-1}$. & Sufficient norm condition, Theorem~\ref{th:bg}; a.e. convergence additionally uses the source's geometric conditions.\\
Varying $\alpha_n\to0$ & Convergence on each fixed $\mathcal V_K$ for all $L^1$ inputs. & Restricted positive theorem \cite{GatGogVary22}; full-sequence negative theorem \cite{GogVary24}.\\
Arithmetic means of selected sums & A prescribed lacunary sequence gives a positive a.e. theorem for every $L^1$ input. & \cite{GatSub10}; the selection hypothesis belongs to the statement.\\\bottomrule
\end{longtable}

\subsection{Model rows and what their hypotheses actually control}\label{sub:matrixmodels}
The squared row-mass condition in Theorem~4.3 describes dispersion of mass. For normalized nonnegative rows, Cauchy--Schwarz gives
\[
 \sum_{k=1}^{d}t_{k,d}^2\ge \frac1d,
\]
with equality precisely for the uniform row. Consequently the sufficient condition bounds the squared mass by a constant multiple of its smallest possible value. It does not require any ordering of the entries. A row uniform on an arbitrary set of at least $\theta d$ positions satisfies it with constant $\theta^{-1}$; that set can consist of scattered indices. A row concentrated at $k=d$ has squared mass one. It fails the condition but can still give convergent means on a separately controlled sequence, such as dyadic conditional expectations. The sufficient criterion is therefore a universal statement about a class of matrices, not a necessity test for one row sequence.

Regularity can also be checked directly from this condition. For fixed $K$,
\[
 \sum_{k\le K}t_{k,d}
 \le \sqrt K\Big(\sum_{k\le d}t_{k,d}^2\Big)^{1/2}
 \le (CK/d)^{1/2}\longrightarrow0.
\]
Every fixed low index loses its mass, and the row sum remains one. This proves stabilization on any fixed polynomial. Uniform kernel bounds supply the additional norm control in the cited theorem. The lacunary restriction needed for the corresponding almost everywhere statement is not a consequence of this finite estimate.

\paragraph{A direct-row test.}
If $q_j$ is nonincreasing and nonnegative, the direct row $a_{N,k}=q_{k-1}/Q_N$ has
\[
 \sum_{k=1}^{N}k|a_{N,k}-a_{N,k+1}|=1,
 \qquad a_{N,N+1}=0.
\]
This is a telescoping equality. Abel's formula then gives $|T_Nf|\le\sigma^*f$ in the positive-index convention. For regular approximation one must also have $q_j/Q_N\to0$ for every fixed $j$, for example $Q_N\to\infty$. If $Q_N$ instead converges to a finite number, a nonzero proportion of mass can remain at low partial sums and the limit need not reproduce all polynomials. Maximal control and regularity are independent checks.

\paragraph{Two accumulated-growth substitutions.}
The criterion for reversed N\"orlund means concerns a different row. Under its hypotheses, $Q_{2^j}\asymp2^{\alpha j}$ gives the strong threshold $p>1/(1+\alpha)$. At the formally critical value, replacing this growth by
$Q_{2^j}\asymp2^{\alpha j}(j+1)^\beta$ makes the criterion comparable to
\[
 (M+1)^{-\beta p}\sum_{j=1}^{M}(j+1)^{\beta p}.
\]
It is unbounded for every real $\beta$: if $\beta p>-1$ it grows like $M$; at equality it grows like $M\log M$; and below $-1$ it grows like $M^{-\beta p}$. Thus a pure logarithmic correction to this particular accumulated-growth model does not restore strong critical boundedness. This calculation uses both-sided comparability and the actual finite sum; it should not be replaced by a general statement that a logarithm always resolves the endpoint.

For logarithmic accumulated growth $Q_{2^j}\asymp j+1$, the criterion at $p<1$ contains an exponentially growing factor $2^{M(1-p)}/(M+1)^p$ multiplying a sum bounded below by a positive constant. At $p=1$ it is comparable to $(M+1)^{-1}\sum_{j\le M}(j+1)\asymp M$. This supplies an explicit failure test for that model whenever the other hypotheses have been checked. These examples are consequences of the criterion, not additional N\"orlund theorems.

\paragraph{Comparison and research boundary.}
The averaging results show where a kernel estimate is sufficient and where an external structural criterion enters. At the critical Hardy exponents, the next issue is localization on atoms: neither regular approximation on polynomials nor a single uniform kernel norm determines the full low-exponent maximum.

\section{Hardy spaces and sharp atomic mechanisms}\label{sec:hardy}
Below the Banach range, a terminal integrable function is no longer the default model. Martingale atoms make the support geometry visible and locate the precise summability threshold. We keep the full maximum, restricted indices, scalar normalizing weights, and spatial measure weights as different problems.

Below the integrable endpoint, martingale atoms expose the scale and cancellation responsible for sharp bounds. The Fej\'er calculation identifies a critical exponent by the number of active coordinates. The later stopping-time tools fix the duality and embedding conventions used in generalized spaces and interval analysis.

\subsection{Martingales, atoms, and coefficient estimates}
Let $\mathcal F_m$ be the sigma-algebra generated by rank-$m$ cylinders. For a dyadic martingale $F=(F_m)$ define $F^*=\sup_m|F_m|$ and
\[
 \norm F_{H^p}=\norm{F^*}_p,\qquad 0<p<\infty.
\]
If $p>1$, this is equivalent to the norm of the terminal $L^p$ function. For $p\le1$, the martingale definition is essential. For example, $F_m=D_{2^m}$ tends to zero almost everywhere but, if $0<p<1$,
\[
 \norm F_{H^p}^p=\frac12\sum_{j\ge0}2^{-j(1-p)}<\infty,
\]
and its stabilized Walsh coefficients all equal one. The terminal zero function therefore does not encode this martingale.

A $p$-atom, $0<p\le1$, is a bounded function supported on a cylinder $I$, with integral zero and $\norm a_\infty\le\mu(I)^{-1/p}$. The standard atomic theorem \cite{Weisz94,PTW22} represents $H^p$ martingales by constants and such atoms, with coefficient sum comparable to the $p$th power of the norm. For a sublinear spectral operator bounded on $L^\infty$, it is therefore sufficient to bound
\begin{equation}\label{eq:atomcriterion}
 \int_{\D\setminus I}|Ta|^p\dd\le C_p
\end{equation}
uniformly in atoms. The integral inside $I$ is bounded by the $L^\infty$ estimate. Apply the finite atomic theorem to finite martingales first, obtain bounds independent of the terminal level and of finite maximal truncations, then pass by completion and monotone convergence. This order avoids an unproved use of an integrable terminal function.

The coefficient theory illustrates how atoms expose a natural scale. One has
\begin{equation}\label{eq:hardylittlewood}
 \sum_{n\ge1}\frac{|\wh F(n)|^p}{n^{2-p}}\le C_p\norm F_{H^p}^p,
 \qquad 0<p\le1.
\end{equation}
For an atom on $I_M$, the coefficients below $2^M$ vanish. For $p<1$, $|\wh a(n)|\le2^{M(1/p-1)}$, and summing the remaining powers proves the assertion on one atom. At $p=1$ use Cauchy--Schwarz:
\[
 \sum_{n\ge2^M}\frac{|\wh a(n)|}{n}
 \le\norm a_2\Big(\sum_{n\ge2^M}n^{-2}\Big)^{1/2}\le C.
\]
Atomic decomposition gives \eqref{eq:hardylittlewood}. The corresponding pointwise coefficient growth $|\wh F(n)|\lesssim_p n^{1/p-1}\norm F_{H^p}$ is consistent with the bit-gap example in Section~\ref{sec:convergence}.

For a general $F=(F_m)\in H^p$, its Fourier coefficient is the stabilized value
\[
 \widehat F(n)=\int_\D F_mw_n\dd\qquad(m\ge u(n)+1),
\]
with $m\ge0$ for $n=0$. Conditional expectation makes this independent of $m$. Thus all finite spectral operators below are defined without assuming an $L^1$ terminal function.

\subsection{Full Fej\'er maximal operator: the critical exponent}
The exact kernel formula in Proposition~\ref{prop:fejerkernel} leads to a proof of the strong range, not just a citation of its value. Repeated use of the Dirichlet recursion gives the positive majorant
\begin{equation}\label{eq:kernelmajorant}
 N|K_N(x)|\le3\sum_{j=0}^{u(N)}2^jK_{2^j}(x).
\end{equation}
Indeed, if $N=2^L+l$, then
$N|K_N|\le2^LK_{2^L}+lD_{2^L}+l|K_l|$ and $D_{2^L}\le2K_{2^L}$.

For $k<l<M$, let $E_{k,l}$ be the set whose first two nonzero coordinates occur at $k,l$. Let $E_{k,M}=I_M(e_k)$, where $e_k$ has a single one at position $k$. These disjoint sets partition $\D\setminus I_M$ and have measures $2^{-l-1}$ and $2^{-M}$, respectively.
\begin{lemma}\label{lem:offatom}
If $N>2^M$, then
\[
 \int_{I_M}|K_N(x\oplus t)|\dd(t)\lesssim\frac{2^{k+l}}{N2^M}
 \quad(x\in E_{k,l}),
\]
and
\[
 \sup_{N>2^M}\int_{I_M}|K_N(x\oplus t)|\dd(t)\lesssim2^{k-M}
 \quad(x\in E_{k,M}).
\]
\end{lemma}
\begin{proof}
On $E_{k,l}$, $K_{2^j}=0$ for $j>l$. For $j\le k$, $K_{2^j}\le2^j$, hence $2^jK_{2^j}\le2^{2j}$; for $k<j\le l$, $K_{2^j}=2^{k-1}$. Thus the sum in that majorant is $O(2^{k+l})$; integration multiplies by $2^{-M}$.

On $E_{k,M}$ the contribution from $j<M$ has integral $O(2^k)$. For $j\ge M$, a nonzero kernel requires the next $j-M$ coordinates to cancel, a relative proportion $2^{-(j-M)}$ of $I_M$. Hence
\[
 \int_{I_M}2^jK_{2^j}(x\oplus t)\dd(t)=2^{k-1}.
\]
If $2^L\le N<2^{L+1}$, division by $N$ gives $O(2^k(L-M+2)/N)\lesssim2^{k-M}$.
\end{proof}
\begin{theorem}\label{th:fejerhardy}
For $p>1/2$, $\norm{\sigma^*F}_p\le C_p\norm F_{H^p}$. For every $0<p\le1/2$, the strong estimate fails. At $p=1/2$ a weak-type $H^{1/2}\to L^{1/2,\infty}$ estimate is known.
\end{theorem}
\begin{proof}[Proof of the strong range and its sharpness]
Let $1/2<p\le1$ and let $a$ be an atom on $I_M$. For $N\le2^M$, $\sigma_Na=0$. Lemma~\ref{lem:offatom} yields
\[
 \sigma^*a(x)\lesssim2^{M(1/p-2)}2^{k+l}\quad(x\in E_{k,l}).
\]
Consequently
\begin{align*}
 \sum_{l<M}\sum_{k<l}\int_{E_{k,l}}|\sigma^*a|^p\dd
 &\lesssim_p2^{M(1-2p)}\sum_{l<M}2^{l(p-1)}\sum_{k<l}2^{kp}\\
 &\lesssim_p2^{M(1-2p)}\sum_{l<M}2^{l(2p-1)}\lesssim_p1.
\end{align*}
On the remaining sets the sum is bounded by $C_p\sum_{k<M}2^{(k-M)p}$. Inside $I_M$ use the uniform $L^1$ bound for the kernels. The atomic theorem applies. For $p>1$, interpolate the weak $(1,1)$ and strong $(\infty,\infty)$ estimates.

For sharpness put $A_M=2^{M(1/p-1)}$ and
\[
 a_M=A_M(D_{2^{M+1}}-D_{2^M})=A_Mr_MD_{2^M}.
\]
Its $H^p$ norm is one. At $N=2^M+1$, $\sigma_Na_M=A_Mr_M/(2^M+1)$, whose $L^p$ norm tends to infinity if $p<1/2$.

For $p=1/2$, take $N_l=2^M+2^l$, $1\le l<M$. Since $S_{2^M+j}a_M=A_Mr_MD_j$ for $1\le j\le2^l$, summing and dividing gives
\[
 \sigma_{N_l}a_M=\frac{A_M2^l}{2^M+2^l}r_MK_{2^l}.
\]
On $E_{k,l}$ this has absolute value at least $2^{k+l-2}$. These sets are disjoint, so
\[
 \int|\sigma^*a_M|^{1/2}\dd
 \gtrsim\sum_{l=1}^{M-1}2^{-l}\sum_{k<l}2^{(k+l)/2}\gtrsim M-1.
\]
This proves failure of the strong endpoint. The weak endpoint is a separate theorem of the Hardy-space theory \cite{Weisz02,PTW22}, not a consequence of the divergent geometric sum above.
\end{proof}
The threshold records two active coordinates of a kernel. It is not a formal consequence of an analogy with trigonometric Ces\`aro means.

\subsection{Restricted indices and two meanings of weight}
\begin{theorem}\label{th:dyadicfejer}
For $\sigma^\#F=\sup_m|\sigma_{2^m}F|$ and every $p>0$,
$\norm{\sigma^\#F}_p\le C_p\norm F_{H^p}$.
\end{theorem}
\begin{proof}
For an atom on $I_M$, all means with $m\le M$ vanish. If $m>M$, the kernel vanishes off $I_M$ whenever the first $M$ digits contain two ones. On $E_{k,M}$, its integral over $I_M$ is $2^{k-1}2^{-m}$. Therefore
\[
 \int_{\D\setminus I_M}|\sigma^\#a|^p\dd
 \lesssim_p2^{-M}2^{M(1-p)}\sum_{k<M}2^{kp}\lesssim_p1.
\]
The integral on the support is bounded by positivity and unit mass of the dyadic kernels. This proves the atomic range $p\le1$; the range $p>1$ follows from the full maximal theorem.
\end{proof}
Even this restricted maximum is not strongly bounded on $L^1$. For $f_M=D_{2^M}$, $\norm{f_M}_1=1$. On the shell whose first one is at $j<M$, $\sigma_{2^{j+1}}f_M=K_{2^{j+1}}=2^{j-1}$. Summing the integrals over these shells gives $\norm{\sigma^\#f_M}_1\ge M/4$.

At the critical exponent, the geometry of a family of indices matters, not only that of one integer. For $\Lambda_s=\Lambda\cap[2^s,2^{s+1})$, collect all endpoints of the maximal strings of ones in the binary expansions of its members into a set $A_s$. Uniform boundedness of $|A_s|$ gives a strong restricted maximal estimate $H^{1/2}\to L^{1/2}$; the corresponding failure and weighted variants are developed in \cite{BPTRestricted23,ABPTWeighted23,BaramidzeTeph24}. A representative normalization divides each $\sigma_n$ by $|A_{u(n)}|^2$. This is an index-dependent normalization, not a spatial measure weight. The family $n=2^s+2^k$, $k<s$, has individually bounded $V(n)$ but an increasing collection of endpoints, illustrating the distinction. Divergence on null sets is another separate issue \cite{AreshidzeZero23}.

By contrast, \cite{WeiszWeighted19} studies spatial weights. For a dyadic Muckenhoupt weight $w\in A_{2p}^d$, $p>1/2$, the Fej\'er maximal estimate is formulated in the corresponding weighted martingale Hardy space. The measure is then $w\dd$, and its condition cannot be replaced by a normalization in $n$. The often-used convention $N^{-1}\sum_{k=0}^{N-1}S_k$ differs from ours by a shifted index and a uniformly bounded factor; exact identities must still be adjusted.

The norms and the standard duality and embedding inputs are recorded in Section~\ref{sub:v2bmo}. There the BMO seminorm, the constant component of $H^1$, and the Carleson normalization are fixed together.

\subsection{Coefficient growth and the endpoint summation argument}\label{sub:v2coeff}
\begin{proposition}\label{v2:coeffgrowth}
For $0<p\le1$, $n\ge1$, and $F\in H^p$,
$|\wh F(n)|\le C_p n^{1/p-1}\norm F_{H^p}$. For $p<1$, the exponent cannot be lowered uniformly.
\end{proposition}
\begin{proof}
For an atom on $I_M$, coefficients vanish below $2^M$, while above that point their modulus is at most $\norm a_1\le2^{M(1/p-1)}\le n^{1/p-1}$. If $F=\sum\lambda_ja_j$, then $\sum|\lambda_j|\le(\sum|\lambda_j|^p)^{1/p}$; for each fixed coefficient this ensures absolute convergence of the atomic sum. The constant atom contributes nothing at positive frequencies. Take the infimum over decompositions. The block atom from Proposition~\ref{v2:gap} attains the stated power.
\end{proof}
The individual bound does not imply \eqref{eq:hardylittlewood}: substituting it gives a divergent harmonic majorant. The atom-by-atom summation proof is stronger. At $p=1$ its use of Parseval, rather than the crude $L^1$ bound on each coefficient, is indispensable.

\subsection{Martingale tools and the extension to \texorpdfstring{$H^p$}{Hp}}\label{sub:v2bmo}
We use the following standard binary-filtration facts; see \cite{Weisz94,Weisz02,PTW22}. They are recorded to fix norms and dependencies. With $b_I=\mu(I)^{-1}\int_I b\dd$, put
\[
 \|b\|_{\mathrm{BMO}_{d,q}}
 =\sup_I\left(\mu(I)^{-1}\int_I|b-b_I|^q\dd\right)^{1/q}.
\]
\begin{proposition}[Doob]\label{v2:doob}
For $1<p<\infty$, $\|M_df\|_p\le p'\|f\|_p$.
\end{proposition}
\begin{theorem}[Dyadic John--Nirenberg]\label{v2:JN}
For every cylinder $I$ and $b\in\mathrm{BMO}_d$,
\[
 \mu\{x\in I:|b-b_I|>t\}
 \le C\mu(I)\exp(-ct/\|b\|_{\mathrm{BMO}_{d,1}}).
\]
In particular the oscillation seminorms are equivalent for $1\le q<\infty$. If the seminorm is zero, $b$ is constant almost everywhere and the distribution bound is read with zero exceptional measure.
\end{theorem}
\begin{theorem}[Carleson embedding]\label{v2:carleson}
If $\alpha_I\ge0$ and $A=\sup_J\mu(J)^{-1}\sum_{I\subseteq J}\alpha_I$, then
\[
 \sum_I\alpha_I|f_I|^2\le4A\|f\|_2^2.
\]
Conversely an embedding constant $B$ implies $A\le B$. For normalized Haar functions $h_I$, the paraproduct $\Pi_b f=\sum_I f_I\langle b,h_I\rangle h_I$ consequently satisfies
\[
 \|b\|_{\mathrm{BMO}_{d,2}}\le\|\Pi_b\|_{2\to2}
 \le2\|b\|_{\mathrm{BMO}_{d,2}}.
\]
\end{theorem}
For duality one either pairs the initial constant separately or works with zero-mean $H^1$ and BMO modulo constants. The important issue for the maximal operators here is their extension from finite martingales: when $p\le1$, an integrable terminal function need not exist.

\begin{lemma}[Atomic extension principle]\label{v2:atomicextension}
Let $0<p\le1$. Suppose a sublinear operator on finite martingales has bounds $\norm{T1}_p\le C$ and $\norm{Ta}_p\le C$ on every $p$-atom, uniformly over all finite spectral truncations under consideration. Assume the finite atomic decomposition has coefficient sum bounded by $C_p\norm F_{H^p}^p$, and that $|TF-TG|\le T(F-G)$. Then $T$ extends boundedly $H^p\to L^p$.
\end{lemma}
\begin{proof}
Include the constant as $a_0=1$. For $F=\sum_{j=0}^s\lambda_ja_j$, $p$-subadditivity gives
$\norm{TF}_p^p\le\sum_j|\lambda_j|^p\norm{Ta_j}_p^p\le C^p\sum_j|\lambda_j|^p$.
The finite atomic theorem proves the norm estimate, and its difference form makes images Cauchy along polynomial approximants. Completeness gives the extension. For a countable maximal family, first apply the argument to finite suprema with a common constant and then use monotone convergence. This avoids defining the operator on a nonexistent integrable terminal function.
\end{proof}

\paragraph{Finite approximation and the initial constant.}\label{sub:hardy-map}
The binary filtration is regular: each child has half its parent's measure. Its $H^p$ atomic theorem therefore permits bounded mean-zero atoms on single cylinders, with the initial constant separate \cite{Weisz94,PTW22}. This is the cylinder form used in Lemma~\ref{v2:atomicextension}.
If $a$ is such an atom, $E_La\to a$ in $H^p$: maximal errors tend to zero almost everywhere and are bounded by $2\|a\|_\infty$. Dominated convergence applies. Finite atomic sums, followed by finite conditional expectations, prove polynomial density. The constant is controlled separately by $T1$. No integrable terminal function is needed for this construction.

\paragraph{Comparison and research boundary.}
The full Fej\'er threshold and its endpoint obstruction are settled by the preceding atom calculations. For restricted families, the sharp worst-case endpoint-count growth is recorded in Section~\ref{sec:boundaries}; a finer question must specify which feature of endpoint placement improves on that worst case.

\section{Approximation, smoothness, and generalized spaces}\label{sec:approximation}
The question now concerns rates and regularity: how do approximation tails describe a function, and which norm or modular measures them? Fixed powers, variable exponents, and slowly varying factors behave differently at critical indices. The elementary examples below test the actual hypotheses used by the generalized-space theorems.

Approximation rates depend on what regularity the norm measures. Cylinder tails, dyadic translation moduli, numerical-frequency differentiation, and variable local normalization are compared through explicit estimates. The section also distinguishes sufficient absolute-convergence scales from regularity assumptions that fail to control the coefficient sum.

\subsection{Dyadic approximation and absolute convergence}
For $1\le p\le\infty$, define the dyadic modulus
\[
 \omega_p(f,2^{-m})=\sup_{h\in I_m}\norm{f(\cdot\oplus h)-f}_p.
\]
For $p=\infty$ in approximation by continuous functions we work in $C(\D)$; an $L^\infty$ version uses essential norms and an explicit equivalence-class convention. Averaging translations over $I_m$ proves
\begin{equation}\label{eq:modulus}
 \norm{f-E_mf}_p\le\omega_p(f,2^{-m})\le2\norm{f-E_mf}_p.
\end{equation}
For finite $p$, use Minkowski's inequality; for $p=\infty$, estimate the supremum directly. If $e_m(f)_p$ denotes the best approximation by frequencies below $2^m$, contractivity of $E_m$ gives
\[
 e_m(f)_p\le\norm{f-E_mf}_p\le2e_m(f)_p.
\]
These elementary estimates are the foundation of dyadic Jackson-type theorems \cite{Yano,Watari,GES}, but ordinary translation smoothness on $[0,1)$ is a separate property.

The martingale scales
\[
 \norm f_{B^s_{p,q}}\asymp\norm{E_0f}_p+
 \Big(\sum_{m\ge0}(2^{ms}\norm{\Delta_mf}_p)^q\Big)^{1/q}
\]
and the corresponding $F^s_{p,q}$ norm with the $\ell^q$ sum inside the $L^p$ norm describe approximation by dyadic blocks. Endpoint and quasi-Banach cases require their own definitions. A convenient sufficient condition for absolute Walsh convergence is
\begin{equation}\label{eq:absolute}
 \sum_{m\ge0}2^{m/2}\norm{\Delta_mf}_2<\infty.
\end{equation}
Indeed, Cauchy--Schwarz on each block bounds the sum of absolute coefficients by its number of terms to the power $1/2$ times its $L^2$ norm. Thus $B^{1/2}_{2,1}$ embeds into the Walsh Wiener algebra. This is a precise sufficient result, not a complete classification of absolute convergence. Work on variation and integrable sums of absolute blocks \cite{MTKh15} addresses more refined hypotheses.

Affine Walsh systems and dyadic wavelets extend the approximation viewpoint. Protasov--Farkov \cite{ProtasovFarkov06} study refinable functions and wavelets on the half-line; Terekhin \cite{Terekhin09} and Astashkin--Terekhin \cite{AstashkinTerekhin18} investigate synthesis and basis properties in Banach and symmetric spaces. A wavelet packet library admits best-basis selection \cite{CoifmanWickerhauser}, but choosing an orthogonal basis by an additive cost is not the thresholding greedy algorithm in one fixed Walsh system.

\subsection{Saturation and filters that reproduce low frequencies}
The Fej\'er multiplier is $(1-n/N)_+$ under \eqref{eq:conventions}. Hence if $f\in L^2$ satisfies $\sum n^2|\wh f(n)|^2<\infty$, then
\begin{equation}\label{eq:saturation}
 N^2\norm{\sigma_Nf-f}_2^2\longrightarrow\sum_{n\ge1}n^2|\wh f(n)|^2.
\end{equation}
For $n<N$, the squared multiplier error multiplied by $N^2$ is $n^2$; for $n\ge N$ it is $N^2\le n^2$. Dominated convergence proves the limit. Conversely, if $\norm{\sigma_Nf-f}_2=o(N^{-1})$, each fixed nonconstant coefficient must vanish. This is the elementary Hilbert-space saturation mechanism.

The comparison with a low-frequency reproducing filter is made in Section~\ref{sub:v2rates}, where its multiplier and best-approximation constants are computed together. This keeps saturation, approximation order, and reproduction in the same normalization.

\subsection{Variable, Musielak--Orlicz, and Lorentz--Karamata scales}\label{sub:v2generalspaces}
The generalized-space literature is not captured by replacing $p$ with $p(x)$ in a constant-exponent formula. For a measurable exponent $0<p_-\le p(x)\le p_+<\infty$, put
\[
 \rho_{p(\cdot)}(f)=\int |f(x)|^{p(x)}\dd(x),\qquad
 \norm f_{p(\cdot)}=\inf\{\lambda>0:\rho_{p(\cdot)}(f/\lambda)\le1\}.
\]
For $0<\norm f_{p(\cdot)}<\infty$, modular inequalities compare $\rho(f)$ with the $p_-$ and $p_+$ powers of the norm, in the order determined by whether the norm is below or above one. The zero case follows from $f=0$ almost everywhere; the infinite case is interpreted in the extended sense, rather than dividing by an infinite norm.

For a cylinder $I$ with measure at most one,
\begin{equation}\label{eq:indicatorvariable}
 \mu(I)^{1/p_-(I)}\le\norm{\ind_I}_{p(\cdot)}\le\mu(I)^{1/p_+(I)}.
\end{equation}
These bounds follow by testing the two constants in the modular. Their ratio is
\[\mu(I)^{1/p_+(I)-1/p_-(I)},\]which can become arbitrarily large without oscillation control. This simple calculation explains why an atom cannot always be normalized by one fixed power of its support. The martingale framework \cite{JiaoVariable}, Walsh--Vilenkin applications \cite{WeiszVariable}, and variable Hardy--Lorentz characterizations \cite{WeiszLorentz23} impose regularity hypotheses tailored to such estimates.

For a growth function $\Phi(x,t)$, the Musielak--Orlicz norm instead uses $\int\Phi(x,|f|/\lambda)\dd$. Lower and upper type, spatial weight conditions, and regularity of the filtration are distinct hypotheses. The Musielak--Orlicz--Lorentz theory of Jiao--Weisz--Xie--Yang \cite{JiaoMO21} and the product martingale theory \cite{MOProduct25} provide atomic and maximal tools beyond fixed powers. An assertion about one of these scales does not automatically apply to all others.

For a nonincreasing rearrangement $f^*$, a typical Lorentz--Karamata expression is
\[
 \norm f_{p,q;b}=\Big(\int_0^1[t^{1/p}b(t)f^*(t)]^q\frac{dt}{t}\Big)^{1/q},
\]
with the supremum modification when $q=\infty$ and a slowly varying factor $b$. The Orlicz--Lorentz--Karamata martingale theory and Fourier applications are developed in \cite{HaoLiLongWeisz24,HaoLiWeisz24,HaoWeisz24}. These are separate papers with different author lists. Variable-exponent martingale inequalities \cite{HaoVariable24} and their Walsh applications \cite{HaoLiWeisz26} are distinct published contributions. The former develops the relationships among five Hardy-space definitions; the latter supplies the Fourier-analytic applications used below. Their author lists and theorem scopes are kept separate.

An endpoint test makes the logarithms concrete. Set, for sufficiently small $t$,
\[
 f^*(t)=t^{-1/p}(\log(e/t))^{-\gamma},\qquad b(t)=(\log(e/t))^\beta,
\]
and continue $f^*$ boundedly and nonincreasingly on the remaining interval. Then for finite $q$ the integral near zero is finite exactly when $(\gamma-\beta)q>1$. A formula declared on all of $(0,1)$ need not be decreasing near one; the continuation is part of the example. The endpoint depends on the secondary parameter $q$ and the slowly varying factor, not only on $p$.

Local Morrey and Campanato formulations require equally explicit normalizations. For example,
\[
 \sup_I\mu(I)^{-\alpha}\Big(\frac1{\mu(I)}\int_I|f-f_I|^q\dd\Big)^{1/q}
\]
is a homogeneous Campanato seminorm. At $\alpha=0$ its equivalence with BMO, for appropriate finite $q$, uses John--Nirenberg; it is not a definition-free identity for every generalized norm. Precise hypotheses are particularly important when importing Fej\'er and fractional-integral estimates into these local scales.

A concrete recent theorem shows the amount of structure needed for an actual generalized-space application. In the notation of Hao--Li--Weisz \cite{HaoLiWeisz26}, let $b$ be slowly varying on $[1,\infty)$ and $\gamma_b(t)=b(\max\{t,t^{-1}\})$. Define
\[
 \norm f_{p(\cdot),q,b}=
 \norm{t\,\norm{\ind_{\{|f|>t\}}}_{p(\cdot)}
 \gamma_b(\norm{\ind_{\{|f|>t\}}}_{p(\cdot)})}_{L^q((0,\infty),dt/t)},
\]
with the usual supremum at $q=\infty$ and zero contribution for empty level sets. The conditional square-function Hardy model $H^s_{p(\cdot),q,b}$ uses
$s(F)=(\sum_n\mathbb E_{n-1}|d_nF|^2)^{1/2}$, including the initial constant component consistently. Their Theorem 1.4 gives
\[
 \norm{\sigma^*F}_{p(\cdot),q,b}\lesssim\norm F_{H^s_{p(\cdot),q,b}}
\]
under $0<q\le\infty$, $1/2<p_-\le p_+<\infty$,
$1/p_- -1/p_+<1$, and the filtration-atom oscillation condition
\[
 \sup_{I\text{ an atom}}\mu(I)^{p_-(I)-p_+(I)}<\infty.
\]
This is a theorem about the conditional square-function model; its identification with a maximal Hardy model uses further regular-filtration results. The atom condition is not implied by global exponent bounds: put $p=a$ on one half and $p=b>a$ on the other half of disjoint cylinders with ranks $m_j\to\infty$. The displayed quantity on those cylinders is $2^{m_j(b-a)}\to\infty$. This example tests a hypothesis, rather than claiming unboundedness of the operator for every exponent failing it.

\subsection{Tail norms, absolute convergence, and a sharp block boundary}\label{sub:v2approx}
Write $\mathcal P_N=\operatorname{span}\{w_n:n<N\}$ and $E_N(f)_p=\inf_{P\in\mathcal P_N}\norm{f-P}_p$. The notation $e_N$ in Section~\ref{sec:greedy} permits arbitrary supports and coefficients; these are different approximation problems.
\begin{proposition}\label{v2:Jackson}
For $1\le p\le\infty$,
\[
 E_{2^m}(f)_p\le\norm{f-E_mf}_p\le2E_{2^m}(f)_p,
 \quad\norm{f-S_Nf}_p\le(1+L_N)E_N(f)_p.
\]
Together with \eqref{eq:modulus}, each of the dyadic rates $E_{2^m}(f)_p=O(2^{-m\alpha})$, $\norm{f-E_mf}_p=O(2^{-m\alpha})$, and $\omega_p(f,2^{-m})=O(2^{-m\alpha})$ is equivalent to the others for every $\alpha>0$.
\end{proposition}
\begin{proof}
Both projectors reproduce the corresponding polynomial spaces. Subtract an arbitrary approximant and use $\norm{E_m}_{p\to p}=1$ or $\norm{S_N}_{p\to p}\le L_N$. The lower best-approximation bound follows because $E_mf$ lies in that space. The modulus comparison proves the last assertion.
\end{proof}
No upper bound $\alpha\le1$ is needed: cylinder polynomials satisfy every such dyadic rate. On this exponent-two group, $T_h^2=I$ implies $(T_h-I)^k=(-2)^{k-1}(T_h-I)$; repeating the same difference does not create Euclidean higher-order smoothness.

\begin{proposition}\label{v2:Besovtails}
For $s>0$, $1\le p,q\le\infty$, the martingale Besov norm is equivalent to
\[
 \norm f_p+\norm{(2^{ms}\omega_p(f,2^{-m}))_{m\ge0}}_{\ell^q}.
\]
For $1\le p<\infty$ the mixed-norm inequalities also give
$B^s_{p,\min(p,q)}\subset F^s_{p,q}\subset B^s_{p,\max(p,q)}$.
\end{proposition}
\begin{proof}
Since $\Delta_jE_jf=0$, $\norm{\Delta_jf}_p\le2\norm{f-E_jf}_p$. In the other direction put $u_j=2^{js}\norm{\Delta_jf}_p$. If $(u_j)\in\ell^q$, its boundedness and $s>0$ give absolute convergence of the block expansion in $L^p$; at $p=\infty$ this supplies a uniformly convergent representative after the low constant is included. Consequently
\[
 2^{ms}\norm{f-E_mf}_p\le\sum_{j\ge m}2^{-(j-m)s}u_j.
\]
The geometric sequence is in $\ell^1$, so discrete convolution is bounded on every indicated $\ell^q$. Combine this with \eqref{eq:modulus}. The final embeddings follow from Minkowski's inequality, with the order of the $L^p$ and sequence norms made explicit.
\end{proof}
The next proposition is an elementary consequence of blockwise Cauchy--Schwarz. Its counterexample establishes the boundary within the displayed Besov scale; it carries no independent novelty claim.

\begin{proposition}[Block criterion for absolute convergence]\label{v2:absolute}
Let $\mathcal A(\D)=\{f:\sum_n|\wh f(n)|<\infty\}$. Then
$B^{1/2}_{2,1}\subset\mathcal A$, whereas $B^{1/2}_{2,q}\subset\mathcal A$ fails for every $q>1$, including $q=\infty$. For $s>1/2$, $B^s_{2,q}\subset\mathcal A$ for $1\le q\le\infty$.
\end{proposition}
The detailed derivation is recorded in Appendix~\ref{app:auxiliary-proofs}, preserving the normalization used here.
For a polynomial $P\in\mathcal P_{2^m}$ there is also the sharp finite-dimensional estimate
$\norm P_\infty\le2^{m/p}\norm P_p$, $0<p<\infty$: the value is constant on every cylinder of measure $2^{-m}$, and a single cylinder indicator attains equality. This fact is useful in approximation arguments, but is a consequence of the partition geometry rather than an additional deep theorem.

\subsection{Campanato regularity in an ultrametric geometry}\label{sub:v2campanato}
Use $\rho(x,y)=2^{-k}$ when $k$ is the first differing coordinate. For $\alpha>0$ and $1\le p<\infty$, put
\[
 [f]_{\mathcal C_p^\alpha}=\sup_I\mu(I)^{-\alpha}
       \Big(\frac1{\mu(I)}\int_I|f-f_I|^p\dd\Big)^{1/p}.
\]
\begin{proposition}\label{v2:campanato}
This seminorm is finite exactly when $f$ has a continuous representative $\widetilde f$ satisfying
$|\widetilde f(x)-\widetilde f(y)|\le C\rho(x,y)^\alpha$.
The equivalence holds for every $\alpha>0$, with constants depending on $p,\alpha$.
\end{proposition}
\begin{proof}
If $J$ is a child of $I$ and the seminorm is $A$, then $|f_J-f_I|\le2A\mu(I)^\alpha$. Thus $E_mf$ is uniformly Cauchy, since the differences from level $m$ onward sum to at most $2A\sum_{j\ge m}2^{-j\alpha}$. Its continuous uniform limit is a representative of $f$. Compare two points with the mean on their last common cylinder to obtain the stated bound. Conversely a H\"older estimate bounds $|f-f_I|$ on $I$ by $C\mu(I)^\alpha$, proving the seminorm bound. There is no Euclidean argument forcing constancy for $\alpha>1$.
\end{proof}
At $\alpha=0$ the convergent geometric sum in this proof is lost. Equivalence of finite-$p$ oscillation norms with BMO then uses Theorem~\ref{v2:JN}, not the positive-$\alpha$ proof. Local Morrey norms and scale summation can be combined in either order to define further spaces; the order is part of the definition and matters for operator estimates.

\subsection{Rates of Fej\'er approximation and logarithmic saturation}\label{sub:v2rates}
\begin{proposition}\label{v2:Fejerrate}
For $1\le p\le\infty$ and an absolutely convergent block expansion,
\[
 \norm{f-\sigma_{2^M}f}_p
 \le2^{1-M}\sum_{j<M}2^j\norm{\Delta_jf}_p
       +\sum_{j\ge M}\norm{\Delta_jf}_p.
\]
If $\norm{\Delta_jf}_p\le A2^{-j\alpha}$, its upper order is $A2^{-M\alpha}$ for $0<\alpha<1$, $A(M+1)2^{-M}$ for $\alpha=1$, and $A2^{-M}$ for $\alpha>1$.
\end{proposition}
\begin{proof}
On the first $j+1$ digits,
\[
 \mathcal D=\sum_{k=0}^j2^{k-1}(I-T_{e_k}),\qquad
 \norm{\mathcal Dg}_p\le(2^{j+1}-1)\norm g_p.
\]
For $j<M$, $(I-\sigma_{2^M})\Delta_j=2^{-M}\mathcal D\Delta_j$, whereas higher blocks are removed by the mean. Triangle inequalities give the bound first on polynomials; absolute block convergence gives the extension. Sum the three geometric-series regimes to get the rates.
\end{proof}
The logarithm at $\alpha=1$ is necessary for the uniform norm; the Hilbert estimate below and Section~\ref{sub:approx-models} have a smaller critical factor. For the continuous function $f=\sum_{j\ge0}2^{-j}r_j$, the error at zero is
\[
 (f-\sigma_{2^M}f)(0)=M2^{-M}+\sum_{j\ge M}2^{-j}=(M+2)2^{-M}.
\]
In $L^2$ the low-block contribution is only $\sqrt M\,2^{-M}$ by orthogonality. This difference is why a Hilbert computation cannot justify the uniform-norm rate. The example is outside the domain of the $L^2$ derivative in \eqref{eq:saturation}.

For the reproducing method $\operatorname{VP}_N=2\sigma_{2N}-\sigma_N$, its multiplier is one on $[0,N)$, $2-n/N$ on $[N,2N)$, and zero thereafter. The kernel norm is at most six if $\sup_N\norm{K_N}_1\le2$. Reproduction therefore gives the explicit estimate
\[
 E_{2N}(f)_p\le\norm{f-\operatorname{VP}_Nf}_p\le7E_N(f)_p.
\]
Its signed matrix coefficients explain why bounded row variation is a useful condition beyond stochastic rows. The contrast with Fej\'er saturation concerns the operator's reproduction property, rather than only the smoothness class.

\begin{example}[An $L^1$ smoothness counterexample]\label{v2:L1absolute}
For $0<\alpha\le1$, let
\[
 f_\alpha=\sum_{j\ge1}2^{-j\alpha}r_jD_{2^j}.
\]
The series converges in $L^1$, since each $r_jD_{2^j}$ has norm one. Its coefficients on block $j$ all equal $2^{-j\alpha}$, so
\[
 \sum_n|\wh f_\alpha(n)|=\sum_j2^{j(1-\alpha)}=\infty.
\]
Yet its $L^1$ martingale tail is $O(2^{-m\alpha})$, and hence so is $\omega_1(f_\alpha,2^{-m})$. Even exponent one in this $L^1$ modulus does not give absolute convergence. At $\alpha=1$ its $L^2$ block norm is $2^{-j/2}$, also giving a critical $L^2$ example. Both the smoothness exponent and the norm in a sufficient theorem are indispensable.
\end{example}

For clarity, the fixed-variable-exponent results in \cite{WeiszVariable} also include natural partial-sum norm convergence in bounded Vilenkin systems when $p_->1$ and the exponent satisfies the source's log-H\"older regularity. Its Fej\'er maximal result uses $p_->1/2$ and $1/p_- -1/p_+<1$, with that regularity retained. On $[0,1)$ the condition is
\[
 |p(x)-p(y)|\le\frac C{\log(e+|x-y|^{-1})}.
\]
The filtration-atom oscillation condition used in the later Lorentz--Karamata statement controls only cylinder oscillation; replacing the Euclidean hypothesis in an earlier theorem by that weaker condition needs an explicit version supporting the replacement. Hardy--Lorentz extensions likewise specify whether their Hardy norm is maximal or square-functional. These are substantive hypotheses, rather than alternate names for one space.

There is a further endpoint distinction in the Hardy--Lorentz convergence conclusions: when $q=\infty$, the closure of polynomials in the weak Hardy--Lorentz space replaces the whole space for the corresponding norm-convergence assertion in \cite{WeiszVariable}. The maximal inequality alone does not prove density of polynomials in a weak norm. This is a second reason to separate boundedness from convergence, even after all exponent conditions have been met.

\paragraph{The Hilbert critical rate and weak-norm density.}\label{sub:approx-models}
The exact identity
\[
 \|f-\sigma_{2^M}f\|_2^2
 =2^{-2M}\sum_{n<2^M}n^2|\widehat f(n)|^2
 +\sum_{n\ge2^M}|\widehat f(n)|^2
\]
gives the critical rate $(M+1)^{1/2}2^{-M}$ when $\|\Delta_jf\|_2\lesssim2^{-j}$. The Rademacher example above attains it and the larger uniform-norm rate simultaneously.
In weak $L^p(0,1)$, $x^{-1/p}$ has distance at least $1/2$ from every bounded function, as $\{f>2\lambda\}\subset\{|f-P|>\lambda\}$ for $\lambda>\|P\|_\infty$. Thus the polynomial-closure qualification in weak Hardy--Lorentz convergence is essential.

\paragraph{Comparison and research boundary.}
These examples connect regularity with approximation only inside the stated norm or modular model. They do not turn a constant-exponent estimate into a variable-space theorem. The next comparison changes the geometry of the selected spectrum rather than the regularity of the target.

\section{Spectral intervals, time--frequency analysis, and weights}\label{sec:timefreq}
A selected spectral interval need not begin at zero or align with a martingale block. Moving boundaries require interval decompositions and time--frequency tools. The comparison proceeds from disjoint intervals to oscillation, bilinear interaction, and spatially weighted estimates, where orthogonality alone no longer supplies the desired norm inequality.

The spectral boundary can move within a block, or several disjoint boundaries can be present simultaneously. Endpoint XOR decompositions connect interval estimates to martingale transforms, while time--frequency trees control a varying cutoff. Spatial weights and sparse bounds quantify these conclusions beyond the unweighted norm.

\subsection{The interval theorem of Osipov and its extensions}
For a finite interval $I\subset\N_0$, let $P_I$ be its Fourier projection. The interval Littlewood--Paley--Rubio de Francia theorem is a substantial strengthening of estimates for complete dyadic blocks.
\begin{theorem}[Osipov \cite{Osipov16}, Theorem 1]\label{th:rubio}
If $I_j$ are pairwise disjoint finite intervals and $\spec f_j\subset I_j$, then
\[
 \norm{\sum_jf_j}_p\le C_p\norm{(\sum_j|f_j|^2)^{1/2}}_p,\qquad1<p\le2.
\]
Equivalently, for $2\le q<\infty$,
\[
 \norm{(\sum_j|P_{I_j}f|^2)^{1/2}}_q\le C_q\norm f_q.
\]
Constants are independent of the number, positions, and lengths of the intervals.
\end{theorem}
The assembly lemma and binary decomposition are proved below. They yield a complete finite-family proof; countable versions follow by the corresponding norm limits. Set $\delta_{-1}=\{0\}$, $\delta_k=[2^k,2^{k+1})$, and $\mathsf d_{-1}=E_0$, $\mathsf d_k=\Delta_k$.
\begin{lemma}[Assembly operator]\label{lem:rubioassembly}
Suppose the sets $a_\nu\mathbin{\mathrm{xor}}\delta_{k_\nu}$ are pairwise disjoint. Then
\[
 Gh=\sum_\nu w_{a_\nu}\mathsf d_{k_\nu}h_\nu
\]
has weak type $(L^1(\ell^2),L^{1,\infty})$ and is bounded $L^p(\ell^2)\to L^p$ for $1<p\le2$.
\end{lemma}
\begin{proof}
For $\lambda\le\|h\|_{L^1(\ell^2)}$, the weak estimate follows from $\mu(\D)=1$. For the decomposition below assume $\lambda>\|h\|_{L^1(\ell^2)}$, so every maximal bad cylinder has a parent. Disjoint spectra give $\norm{Gh}_2^2=\sum_\nu\norm{\mathsf d_{k_\nu}h_\nu}_2^2\le\norm h_{L^2(\ell^2)}^2$. Apply the dyadic good--bad decomposition at level $\lambda$ to $\norm{h(x)}_{\ell^2}$. Replace $h$ by its vector average on the selected cylinders. The good vector $g$ satisfies $\norm g_{L^\infty(\ell^2)}\le2\lambda$ and $\norm g_{L^1(\ell^2)}\le\norm h_{L^1(\ell^2)}$. Each bad vector has coordinatewise mean zero on its cylinder. Its low martingale differences vanish, and the remaining differences remain supported on that cylinder; modulation does not enlarge support. Outside the bad union, $Gh=Gg$. Chebyshev and the $L^2$ estimate give measure at most $3\norm h_{L^1(\ell^2)}/\lambda$. Interpolation proves the strong range.
\end{proof}
The same argument, followed by self-adjoint duality, proves uniform boundedness on $L^p(\ell^2)$ of coordinatewise martingale transforms with fixed coefficients in $\{0,1\}$. Randomization and Hilbert-valued Khintchine then give
\begin{equation}\label{eq:vectorlp}
 \norm{(\sum_{j,k}|\mathsf d_kg_j|^2)^{1/2}}_p
 \lesssim_p\norm{(\sum_j|g_j|^2)^{1/2}}_p,
 \qquad1<p<\infty.
\end{equation}

The other ingredient is a binary decomposition using both interval endpoints. For $a<b$, let $h$ be their largest differing digit and $c=2^h\lfloor b/2^h\rfloor$. Then $[a,b)=L\cup R$ with $L=[a,c)$ and $R=[c,b)$, and
\[
 a\mathbin{\mathrm{xor}}L=\delta_{-1}\cup\bigcup_{\substack{k<h\\a_k=0}}\delta_k,
 \qquad
 b\mathbin{\mathrm{xor}}R=\bigcup_{\substack{k<h\\b_k=1}}\delta_k.
\]
For the first identity compare $n\ge a$ with $a$ at the highest differing lower digit; that digit is zero in $a$ and one in $n$. The second identity reverses this comparison for $n<b$.

\begin{proof}[Completion of Theorem~\ref{th:rubio}]
First work with finitely many polynomials. For each interval $I_j=[a_j,b_j)$ use the left and right pieces $L_j,R_j$ above. On multiplying $f_j$ by $w_{a_j}$, projection onto $a_j\mathbin{\mathrm{xor}}L_j$ is a martingale transform selecting the blocks $\delta_{-1}$ and the indicated $\delta_k$. Define
\[
 g_j^L=w_{a_j}P_{L_j}f_j,
 \qquad g_j^R=w_{b_j}P_{R_j}f_j.
\]
The coordinatewise transform estimate gives
\[
 \|(g_j^L)_j\|_{L^p(\ell^2)}+
 \|(g_j^R)_j\|_{L^p(\ell^2)}
 \le C_p\|(f_j)_j\|_{L^p(\ell^2)}.
\]
Indeed the input modulation is an isometry of $L^p(\ell^2)$, each coordinate transform has deterministic zero--one coefficients, and the output has no frequencies outside its selected blocks. Apply \eqref{eq:vectorlp} separately to the two vectors. The square sums of their block differences are bounded by the displayed right side.

For every nonzero block of $g_j^L$ its remodulated spectrum is $a_j\mathbin{\mathrm{xor}}\delta_k\subset L_j$; the corresponding statement holds for $R_j$. These sets are pairwise disjoint within either assembly, since the original $I_j$ are disjoint and the binary pieces partition each interval. Lemma~\ref{lem:rubioassembly} therefore gives
\[
 \Big\|\sum_jP_{L_j}f_j\Big\|_p+
 \Big\|\sum_jP_{R_j}f_j\Big\|_p
 \le C_p\|(f_j)_j\|_{L^p(\ell^2)}.
\]
Adding establishes synthesis.

For completeness, the duality requires a vector projection estimate. Let $(h_j)$ be arbitrary. Formula~\eqref{eq:transform}, applied coordinatewise with the possibly different $a_j$ or $b_j$, and the vector martingale-transform bound yield
\[
 \|(P_{I_j}h_j)_j\|_{L^p(\ell^2)}
 \le C_p\|(h_j)_j\|_{L^p(\ell^2)}.
\]
Given $f\in L^q$, $q=p'$, test the vector $(P_{I_j}f)_j$ against $(h_j)$ in the unit ball of $L^p(\ell^2)$. Self-adjointness of the finite Fourier projections and synthesis for $P_{I_j}h_j$ bound the pairing by
$C_p\|f\|_q\|(h_j)_j\|_{L^p(\ell^2)}$. Taking the supremum proves the analysis estimate. Density handles $L^p$ inputs. For a countable family in the analysis inequality, apply the finite result and monotone convergence to the square function. A synthesis series is asserted only when its right-hand square function lies in $L^p$; the same estimate applied to tails then gives its norm convergence.
\end{proof}

Borovitskiy's two-parameter theorem \cite{Borovitskiy20} concerns pairwise disjoint spectral rectangles, whose individual coordinate projections need not be disjoint. It uses product martingale atomic theory; simply iterating the one-dimensional theorem would miss this geometry. Tselishchev extends interval inequalities to bounded Vilenkin systems \cite{Tselishchev21} and studies vector-valued Walsh inequalities \cite{TselishchevVector21}. Multiparameter and unbounded-base advances appear in \cite{BorovitskiyMulti24,TselishchevUnbounded24}. This chain is an essential modern contribution of the Russian school, rather than a historical footnote.

\subsection{Packets, trees, and quantitative oscillation}
For a dyadic interval $I=[a,a+|I|)$, the normalized packet
\[
 p_{I,n}(x)=|I|^{-1/2}w_n((x-a)/|I|)\ind_I(x)
\]
has unit $L^2$ norm. Associate it to the tile $I\times[n/|I|,(n+1)/|I|)$ of area one. The packets for $n=0,1$ are the normalized indicator and Haar function on $I$. Nested dyadic intervals and frequency nesting organize tiles into trees. Tree estimates separate energy from density; a stopping-time decomposition then sums the trees. This is the basic architecture of Walsh time--frequency analysis \cite{Thiele06,DemeterGuide16}, and it controls a moving spectral boundary unavailable to one fixed Doob maximal theorem.

For $0<r<\infty$ and a scalar sequence $a_N$, define
\[
 V^r(a_N)=\sup_{N_0<\cdots<N_J}
 \Big(\sum_{j=1}^J|a_{N_j}-a_{N_{j-1}}|^r\Big)^{1/r}.
\]
Variation controls accumulated oscillation; its deterministic jump consequence is the following.

\begin{proposition}\label{v2:jumps}
Let $N_\lambda(z)$ count successively disjoint ordered pairs $a_j<b_j\le a_{j+1}$ with $|z_{b_j}-z_{a_j}|>\lambda$. Then $\lambda N_\lambda(z)^{1/r}\le V^r(z)$. Finite $r$-variation of a scalar sequence implies that it has a finite limit.
\end{proposition}
\begin{proof}
Insert all endpoints of any $K$ such pairs into one increasing partition. Its sum of $r$th powers is at least $K\lambda^r$. Take the supremum. A sequence that is not Cauchy supplies infinitely many successively disjoint pairs separated by one positive $\lambda$, contradicting finite variation. Boundedness follows already by comparing each term with the first, so the Cauchy limit is finite.
\end{proof}
In particular $\sup_{\lambda>0}\lambda N_\lambda^{1/r}\le V^r$ deterministically. An estimate with exponent $1/2$ in place of $1/r$ when $r>2$ is a genuinely stronger jump question; it does not follow from that inequality.

A precise weighted range from Do--Lacey \cite{DoLaceyWeighted} is
\begin{equation}\label{eq:dolacey}
 \norm{V^r(S_Nf)}_{L^p(w)}\lesssim\norm f_{L^p(w)},
 \quad w\in A_a^d,\quad1\le a<p,\quad r>2a,\quad
 \frac1r<\frac1a-\frac1p.
\end{equation}
Here $[w]_{A_a^d}$ is formed over dyadic cylinders (with the essential-infimum definition at $a=1$). For $w=1$, taking $a=1$ gives $r>\max\{2,p'\}$, equivalently $r>2$ and $p>r'$. The scalar specialization of \cite{HLP} has this same range; the condition $p>r'$ must not be dropped.

Indeed $1/r<1-1/p$ is exactly $p>r'=r/(r-1)$ for $r>1$. Finite variation controls both the number of large jumps and the existence of a limit; the estimate therefore measures oscillation in addition to the size controlled by a maximal theorem. Vector-valued bounds additionally depend on tile type or the range-space geometry \cite{HLP,HytLacey}.

The outer-$L^p$ approach of Amenta--Uraltsev \cite{AmentaUraltsev20} packages tree sizes and Carleson embeddings into a framework that also supports vector-valued and sparse conclusions. Their model is the \emph{ternary} Walsh model. The role of modulation and tree geometry transfers conceptually; its statements should not be quoted as literal theorems for the binary Paley system without translating the model.

\subsection{Bilinear interaction: the quartile model}\label{sub:quartile}
Linear interval projections select pieces of one spectrum. A bilinear model must additionally control interactions between two inputs and an output frequency. The quartile operator is the classical Walsh model of the bilinear Hilbert transform. This paragraph uses the standard dyadic phase plane $\mathbb R\times\mathbb R_+$ rather than the compact probability normalization of the preceding sections.

A quartile $s=I_s\times\omega_s$ has area four. Its frequency interval has four dyadic grandchildren; let $w_{s,j}$ be the $L^2$-normalized wave packet for the $j$th child, with the standard consistent numbering. For a finite collection $\mathcal S$, the model is
\[
 \mathcal V_{\mathcal S}(f_1,f_2)
 =\sum_{s\in\mathcal S}|I_s|^{-1/2}
   \langle f_1,w_{s,1}\rangle\langle f_2,w_{s,2}\rangle w_{s,3}.
\]
Its uniform strong range includes
\[
 \|\mathcal V_{\mathcal S}(f_1,f_2)\|_{p_3}
 \lesssim_{p_1,p_2}\|f_1\|_{p_1}\|f_2\|_{p_2},
 \quad 1<p_1,p_2<\infty,\quad
 \frac1{p_3}=\frac1{p_1}+\frac1{p_2},\quad p_3>\frac23;
\]
see \cite{Thiele06,OberlinThiele11}. The output can be quasi-Banach. At $p_3=2/3$, Demeter--Di Plinio \cite{DemeterDiPlinio13} obtain restricted-input weak estimates and, for $p_1,p_2>1$ with $1/p_1+1/p_2=3/2$, the consequence
\[
 \|\mathcal V_{\mathcal S}(f_1,f_2)\|_{2/3,\infty}
 \lesssim_{p_1,p_2}\|f_1\|_{p_1}\|f_2\|_{L^{p_2,2/3}}.
\]
This does not assert the unrestricted strong endpoint.

\paragraph{Why the extra geometry matters.}
The interval assembly argument above uses disjoint output spectra. Quartile trees instead organize several correlated packet families; a size estimate must be paired with energy control for the other inputs. Multi-frequency Calder\'on--Zygmund decomposition preserves cancellation against a selected packet family rather than just one constant on each bad cylinder. This distinction connects the nonlinear model to modern endpoint methods without identifying it with an ordinary martingale transform. The weak-type Carleson theorem of Di Plinio--Fragkos \cite{DiPlinioFragkos25} provides a further development of quantitative multi-frequency methods for ordinary Fourier sums. Its Fourier endpoint statements are not invoked here as new Walsh endpoint theorems.

\subsection{Quantitative weights and sparse multiplier estimates}
For $1<q<\infty$,
\[
 [w]_{A_q^d}=\sup_I\Big(\frac1{\mu(I)}\int_Iw\dd\Big)
 \Big(\frac1{\mu(I)}\int_Iw^{-1/(q-1)}\dd\Big)^{q-1}.
\]
The weighted theory distinguishes membership in $A_q^d$, dependence on its characteristic, variation, and endpoint Orlicz control. Di Plinio--Lerner \cite{DiPlinioLerner14}, together with \cite{DiPlinio}, provide quantitative bounds for the Walsh--Carleson operator; one resulting $A_p$ upper exponent is $\max\{p',2/(p-1)\}$. An upper exponent is not automatically sharp throughout the full parameter range. The variational theorem \eqref{eq:dolacey} is a stronger oscillation statement with a different range restriction.

A family $\mathcal S$ is sparse if each $I\in\mathcal S$ contains a measurable $E_I$ with $\mu(E_I)\ge\eta\mu(I)$ and the $E_I$ are disjoint. For a Young function $\Psi$, write
\[
 \langle f\rangle_{I,\Psi}=\inf\{\lambda>0:\mu(I)^{-1}\int_I\Psi(|f|/\lambda)\dd\le1\}.
\]
Sparse domination reduces bilinear estimates to sums $\sum_{I\in\mathcal S}\mu(I)\langle f\rangle_{I,\Psi}\langle g\rangle_{I,q}$. It is not merely a different name for boundedness on an unweighted $L^p$ space.

\paragraph{A weighted endpoint for arbitrary Walsh intervals.}
For $w>0$, write $w(I)=\int_Iw\dd$ and use the dyadic Fujii--Wilson characteristic
\[
 [w]_{A_\infty^d}
 =\sup_I\frac1{w(I)}\int_I M_d(w\ind_I)\dd,
 \qquad [w]_{A_1^d}=\|M_dw/w\|_\infty.
\]
\begin{theorem}[Di Plinio--Fl\'orez-Amatriain--Parissis--Roncal \cite{DFPR25}, Theorem D]\label{th:weighted-rubio}
Let $\Omega$ be a family of pairwise disjoint finite intervals of nonnegative integers and $\mathcal R_\Omega f=(\sum_{\omega\in\Omega}|P_\omega f|^2)^{1/2}$. Uniformly in $\Omega$,
\[
 \|\mathcal R_\Omega f\|_{L^2(w)}
 \lesssim [w]_{A_1^d}^{1/2}[w]_{A_\infty^d}^{1/2}
          \|f\|_{L^2(w)},\qquad w\in A_1^d.
\]
For $2<p<\infty$ and $w\in A_{p/2}^d$, one also has
\[
 \|\mathcal R_\Omega f\|_{L^p(w)}
 \lesssim_p [w]_{A_{p/2}^d}\|f\|_{L^p(w)}.
\]
The exponent one in the latter characteristic is sharp.
\end{theorem}
\paragraph{Proof mechanism and structural consequence.}
The external proof compares the interval square function with an intrinsic packet square function through a weighted Chang--Wilson--Wolff estimate. The packet estimate has the form $\int|\mathcal W_\Omega f|^2w\lesssim\int|f|^2M_dw$; $A_1^d$ control and extrapolation then enter. Thus weighted orthogonality requires an additional localization theorem, not just Parseval. This sharpens the unweighted interval map and shows a setting where Walsh localization resolves a weighted endpoint whose real Fourier counterpart is more delicate. We use the dyadic characteristics justified in the source's Walsh proof. For the $p>2$ display, the source's proof, Section 5.3, supplies the $L^p(w)$ norm on both sides.

For Walsh Marcinkiewicz multipliers, bounded variation on each dyadic frequency block defines the classical regular symbol class used in Section~\ref{sec:multipliers}. The endpoint coefficient geometry is already more delicate than Parseval: independent characters have the $L(\log L)^{1/2}$ coefficient estimate in Appendix~\ref{app:orthogonal}, and the indicator test there proves the sharp logarithmic power for that estimate. Related lacunary and rearrangement-invariant results are treated in \cite{SagherZhou,Carrillo}. A coefficient inequality is not itself a sparse bound for the full multiplier. Sparse or weighted operator conclusions require an additional localization theorem, such as the published interval theorem above, with its actual operator and weight hypotheses.

Permutations and Weyl multipliers form a related but distinct classical-modern line. Rademacher--Menshov estimates concern maximal partial sums of an orthogonal sequence with logarithmic losses; arbitrary permutations are not controlled by natural-order Carleson convergence. Recent work on Calder\'on--Zygmund maximal operators and Weyl multipliers \cite{KaragulyanLaceyNavoyan25} places these questions in a broader framework. Its conclusions must be specialized with the relevant operator hypotheses rather than inferred solely from orthogonality of Walsh characters.

\paragraph{Comparison and research boundary.}
Interval and packet estimates control moving spectral boundaries; their weighted constants depend on additional localization. This is why an arbitrary bounded diagonal symbol cannot be treated as an automatic consequence of Fourier orthogonality. The next section separates symbol boundedness from operator-tail control.

\subsection{What localization, interval assembly, and trees each control}\label{sub:mechanism-comparison}
The preceding methods can be compared through the location of the cancellation they preserve. For an atom on $I_M$, the off-support kernel estimate in Lemma~\ref{lem:offatom} records the first two nonzero spatial digits. The summation in its strong-range proof contains
\[
 2^{M(1-2p)}\sum_{l<M}2^{l(2p-1)}.
\]
For $p>1/2$ the geometric sum compensates the prefactor; at $p=1/2$ it becomes a count of scales. The lower atom realizes that count on disjoint sets. Thus this proof identifies both the threshold and the mechanism of failure: the shell contributions cease to be summable. A positive complete-block kernel alone would not detect the obstruction.

Interval synthesis uses another cancellation. Modulating an interval endpoint converts the adjacent binary pieces into martingale differences. Orthogonality then sums disjoint output pieces in $L^2$, while the vector good--bad decomposition keeps the bad part supported on its original cylinders. The bound survives arbitrarily many intervals because disjointness is used before the sum is estimated. If output pieces overlap without control, the $L^2$ step must carry their multiplicity, and the family-uniform conclusion is lost. The vector input is equally important: a sum of separate scalar bounds would charge for the number of intervals.

A moving boundary in a Carleson maximum has neither a fixed spatial atom nor one disjoint interval family. Linearization selects a boundary depending on the spatial point, and the associated tiles overlap in the phase plane. Tree methods separate two tasks: an energy bound controls coefficients inside a tree, and a density or size stopping rule controls how trees may accumulate. Their contribution is the packing estimate needed after spatial selection. This explains why the uniform weak bound for individual $S_N$ does not prove a weak bound for $\sup_N|S_N|$. Even for a finite set of $K$ indices, the direct union bound gives only
\[
 \mu\{\max_{j\le K}|S_{N_j}f|>\lambda\}
 \le CK\|f\|_1/\lambda.
\]
Removing the factor $K$ is an additional operator assertion, not a change of notation. The time--frequency sources cited above provide the needed organization for their stated exponent ranges.

\begin{longtable}{@{}p{.19\textwidth}p{.25\textwidth}p{.23\textwidth}p{.20\textwidth}@{}}
\caption{Comparison of cancellation mechanisms}\label{tab:mechanisms}\\
\toprule Mechanism & Information retained & What pays for the bound & Failure test\\\midrule\endfirsthead
\toprule Mechanism & Information retained & What pays for the bound & Failure test\\\midrule\endhead
Atom and shell estimate & Support, mean zero, and the first active spatial digits. & Summable off-support shells plus an on-support norm bound. & Critical Fej\'er shells accumulate at $p=1/2$.\\
Endpoint interval assembly & Endpoint XOR decomposition and disjoint output spectra. & $L^2$ orthogonality and a vector good--bad decomposition. & Uncontrolled spectral overlap invalidates the family-uniform $L^2$ step.\\
Tree organization & Nested spatial intervals and a selected frequency boundary. & Energy inside trees and packing across their tops. & Separate weak bounds leave a factor equal to the number of selected indices.\\\bottomrule
\end{longtable}
The useful transfer question is therefore concrete: which support, disjointness, or packing hypothesis survives in a new system? The binary kernel formulas answer it for atoms, and the endpoint-XOR formulas answer it for intervals. A further weighted or vector-valued result must supply the corresponding weighted packing or range-space estimate. This comparison identifies the added input instead of inferring a stronger operator theorem from an unchanged scalar identity.

\subsection{Walsh models and what a Fourier transfer requires}\label{sub:model-transfer}
The word ``model'' covers several different operators. Thiele's quartile exposition \cite{ThieleQuartile00} links Walsh pointwise convergence with a dyadic bilinear model. The bi-Carleson model of Muscalu--Tao--Thiele \cite{MTTBiCarleson02} combines bilinear interaction with a selected truncation, whereas their biest model \cite{MTTBiest04} is trilinear. They are not alternative names for one theorem. The scalar Carleson maximum, the quartile form, and the higher-order models have different linearizations, admissible exponent regions, and endpoint questions.

For the scalar maximum, first restrict to finitely many stopping indices and choose measurably an index $N(x)$ and a phase so that its linearization is paired with a testing function. The estimate must be uniform over this spatial selection. Complete-block martingale estimates control each fixed block but do not control this choice. In the Walsh phase plane, packet support in a dyadic spatial interval is exact, and its frequency interval is exact in the Walsh packet algebra. A tree organizes simultaneous spatial and frequency nesting. Energy estimates control the coefficients of a single tree; stopping and packing estimates control the aggregate contribution of many trees. The linearized bound, uniform in $N(x)$, is what permits the maximum. This outlines the role of the method in the cited proofs, rather than replacing their selection and tree lemmas.

In the real Fourier setting, a nonzero packet cannot simultaneously have compact spatial and Fourier support. Tail estimates and interactions between translated or dilated packet families must replace exact Walsh support statements. Thus the useful implication is methodological: one isolates an energy, size, or packing principle and proves its analogue with the necessary tails. A Walsh bound by itself is not a transference theorem for the corresponding real operator. At an endpoint, both the estimate and the manner of selecting the exceptional set must be checked again.

\begin{longtable}{@{}p{.24\textwidth}p{.30\textwidth}p{.33\textwidth}@{}}
\caption{Model comparison: a shared method and the additional Fourier step}\label{tab:model-transfer}\\
\toprule Problem & Walsh mechanism & Additional issue in a Fourier argument\\\midrule\endfirsthead
\toprule Problem & Walsh mechanism & Additional issue in a Fourier argument\\\midrule\endhead
Scalar selected cutoff & Exact packet support, tree energy, and packing for a linearized maximum. & Smooth packet localization and summable off-support interactions.\\
Disjoint interval square function & Endpoint XOR decomposition and disjoint output spectra. & Frequency modulation uses ordinary addition; a weighted endpoint needs its own localization estimate.\\
Quartile interaction & Three packets in a dyadic frequency configuration. & The singular bilinear multiplier requires its own discretization and exponent analysis.\\
Bi-Carleson or biest & A selected bilinear cutoff or a trilinear interaction, respectively. & Different operators and selection structures must remain separate.\\\bottomrule
\end{longtable}
The comparison also separates three notions of uniformity: independence of the number of intervals, independence of a spatially selected cutoff, and independence of a parameter in a bilinear model. One theorem may supply one of them without supplying the others. The interval assembly proved here gives the first; the cited Carleson and model theories address their explicitly stated selections and parameters.

\section{Dyadic calculus, multipliers, and semigroups}\label{sec:multipliers}
Diagonal coefficients determine Hilbert-space norms exactly, while a conditional Banach basis also has operator geometry. This section compares numerical frequency, martingale scale, and Hamming degree through scalar spectra, finite-rank tails, regular symbol classes, and cross-space smoothing. Positive semigroups are constructed from their character kernels and closed quadratic forms; a Banach-space conclusion retains the additional boundedness hypothesis it needs.

\subsection{Three diagonal calculi and their domains}
A Walsh polynomial admits several natural diagonal operators. The numerical-frequency derivative has symbol $n$, the martingale-scale derivative has symbol comparable to $2^{u(n)}$, and the Boolean number operator has symbol $h(n)=\sum_jn_j$. They encode different notions of regularity. The Gibbs--Butzer tradition \cite{ButzerWagner,Golubov18,DyadicBook15} concerns a dyadic differentiation built from Walsh characters. Its Chinese development, including Weiyi Su's work documented in \cite{DyadicBook15}, and the martingale Hardy inequalities of Nie--Li--Lou \cite{NieLiLou05}, belong to this history alongside the Russian and European contributions.

For the numerical symbol, define on polynomials
\[
 \mathcal Df=\sum_{n\ge0}n\wh f(n)w_n.
\]
In $L^2$ its natural closed domain is exactly $\sum n^2|\wh f(n)|^2<\infty$. This identifies the saturation limit \eqref{eq:saturation}. A spatial difference formula must be proved to represent this operator; it does not follow simply from naming it a derivative. Likewise, on one martingale block the scale multiplier $2^{m\alpha}$ is exact by definition, whereas a Bernstein estimate for a different symbol uses multiplier bounds inside the block.

For general $a=(a_n)$ write $T_af=\sum a_n\wh f(n)w_n$ initially on polynomials. A classical sufficient Marcinkiewicz condition is
\begin{equation}\label{eq:marcink}
 \sup_n|a_n|+\sup_m\sum_{n=2^m}^{2^{m+1}-1}|a_{n+1}-a_n|<\infty.
\end{equation}
It implies $L^p$ boundedness for $1<p<\infty$; Hardy-space variants on locally compact Vilenkin groups appear in \cite{OnneweerQuek}. Condition \eqref{eq:marcink} is block variation, not automatically a Euclidean H\"ormander--Mikhlin condition. Episkoposian's published multiplier survey \cite{EpiReview26} concerns boundedness, compactness, and spectral questions. The statements below separate the completely diagonal Hilbert case from the Banach-space problem.

\subsection{Compactness: exact Hilbert facts and a general tail criterion}
\begin{theorem}\label{th:hilbertmult}
On $L^2(\D)$, $T_a$ is bounded exactly when $a\in\ell^\infty$, with norm $\norm a_\infty$. It is compact exactly when $a\in c_0$. For bounded $a$,
\[
 \spec(T_a)=\overline{\{a_n:n\ge0\}},\qquad
 \norm{T_a}_{\mathrm{ess}}=\limsup_n|a_n|.
\]
For $0<r<\infty$,
\[
 T_a\in\mathcal S_r\quad\Longleftrightarrow\quad\sum_n|a_n|^r<\infty,
 \qquad \norm{T_a}_{\mathcal S_r}^{r}=\sum_n|a_n|^r.
\]
For a compact diagonal, each nonzero eigenvalue has multiplicity $|\{n:a_n=\lambda\}|<\infty$; the zero eigenspace can be infinite-dimensional.
\end{theorem}
\begin{proof}
Parseval gives the norm and the norm of every diagonal tail. Finite diagonals are finite rank. If $K$ is compact, then $Kw_n\to0$ since the characters are weakly null; this gives the lower bound for the essential norm. Outside the closure of the diagonal entries, $(a_n-z)^{-1}$ defines a bounded inverse. Inside it, characters give eigenvectors or approximate eigenvectors. The singular values are the moduli of diagonal entries, with multiplicities, giving the Schatten formula.
\end{proof}
The closure in the spectral formula and the power $r$ in the Schatten formula are essential.

\begin{proposition}[Tail criterion]\label{prop:compacttail}
Let $1\le p<\infty$, and suppose $T_a$ is already bounded on $L^p$. Then
\[
 T_a\text{ is compact}\quad\Longleftrightarrow\quad
 \norm{(I-E_m)T_a}_{L^p\to L^p}\longrightarrow0.
\]
\end{proposition}
The detailed derivation is recorded in Appendix~\ref{app:auxiliary-proofs}, preserving the normalization used here.
The tail criterion is a general finite-rank approximation principle. Compactness implies $a_n\to0$ also for $p=1$: $w_n$ is weakly null in $L^1$, since every dual test belongs to $L^\infty\subset L^2$ and has Walsh coefficients tending to zero. For $p>1$, polynomial density in the dual space gives the same conclusion. Scalar decay alone does not prove boundedness or tail-norm decay for an arbitrary symbol.

\begin{proposition}[Consequence of Marcinkiewicz boundedness and interpolation]\label{prop:marccompact}
Assume \eqref{eq:marcink}. For every $1<p<\infty$,
\[
 T_a:L^p(\D)\to L^p(\D)\text{ is compact}
 \quad\Longleftrightarrow\quad a_n\longrightarrow0.
\]
No vanishing condition on the block variation is required.
\end{proposition}
\begin{proof}
Necessity follows from weak nullity of the characters. For sufficiency put $R_M=(I-E_M)T_a$ and $\alpha_M=\sup_{n\ge2^M}|a_n|$. Parseval gives $\|R_M\|_{2\to2}=\alpha_M\to0$. For $p<2$ choose $1<r<p$; for $p>2$ choose $r>p$. The Marcinkiewicz theorem bounds $T_a$ on this second space and
\[
 \|R_M\|_{r\to r}\le2\|T_a\|_{r\to r}\le C_r\mathfrak M(a),
\]
where $\mathfrak M(a)$ is the left side of \eqref{eq:marcink}. If $1/p=(1-\theta)/r+\theta/2$, then $0<\theta<1$ and interpolation yields
\begin{equation}\label{eq:compact-interpolation}
 \|R_M\|_{p\to p}
 \le(C_r\mathfrak M(a))^{1-\theta}\alpha_M^\theta\longrightarrow0.
\end{equation}
The finite-rank operators $E_MT_a$ converge in norm. The case $p=2$ is immediate.
\end{proof}
This is a direct consequence of the classical multiplier theorem and compact interpolation in the sense of Krasnosel'skii \cite{Krasnoselskii60}. The explicit tail proof above uses Riesz--Thorin on finite-rank remainders and also records the rate; the general compact-interpolation theorem supplies the qualitative conclusion. The regularity hypothesis supplies boundedness on an exponent beyond the desired $p$. Boundedness only at that $p$ would not supply the second endpoint in this proof.

For example, $a_n=(-1)^n/(n+1)$ satisfies \eqref{eq:marcink} and tends to zero, while its block variation tends to $2\log2$. Its multiplier is compact on every interior $L^p$. Thus requiring variation itself to vanish would discard a simple compact example. The statement does not assert that arbitrary null sequences are $L^p$ multipliers.

An often tempting coefficient shortcut is unavailable: a uniform synthesis map $\ell^{p'}\to L^p$ cannot hold for $1<p<2$. With $M$ unit coefficients on independent Rademacher frequencies, the synthesized norm is comparable to $\sqrt M$, while the coefficient norm is $M^{1/p'}=o(\sqrt M)$. Appendix~\ref{app:orthogonal} gives this test explicitly. Thus scalar decay and a formal coefficient estimate cannot replace operator-tail control. The Hilbert theorem, the tail criterion, and the Marcinkiewicz consequence above have independent proofs.

\subsection{Positive semigroups and fractional potentials}
For $t>0$, the numerical-frequency heat-type multiplier $e^{-tn}$ has kernel
\begin{equation}\label{eq:semigroup}
 P_t(x)=\prod_{j\ge0}(1+e^{-t2^j}r_j(x)).
\end{equation}
The product converges uniformly for each positive $t$, is nonnegative, has integral one, and has Fourier coefficients $e^{-tn}$. Thus convolution is a positive contraction on every $L^p$, $1\le p\le\infty$, satisfies the semigroup law, and converges strongly to the identity for finite $p$. On $L^2$ its infinitesimal generator is $-\mathcal D$ on the domain above; we call $\mathcal D$ the nonnegative generator, so that $P_t=e^{-t\mathcal D}$. For an arbitrary nonnegative symbol $\lambda_n$, the operators $e^{-t\lambda_n}$ are $L^2$ contractions; positivity additionally requires positive definiteness on the dual group and cannot be inferred from $\lambda_n\ge0$ alone.

A transparent dyadic fractional potential is
\[
 I_\alpha f=E_0f+\sum_{m\ge0}2^{-\alpha m}\Delta_mf,
 \qquad0<\alpha<1.
\]
Its convolution kernel is bounded in modulus by $C_\alpha\rho(x,0)^{\alpha-1}$ away from zero, where $\rho(x,0)=2^{-\min\{j:x_j=1\}}$. The positive normalized kernel in Section~\ref{sub:v2positivepotential} gives the bound for this normalization by an affine adjustment of the constant component.

\begin{proposition}\label{prop:fractional}
If $1<p<1/\alpha$ and $1/q=1/p-\alpha$, then $I_\alpha:L^p\to L^q$ is bounded. At $p=1$, $I_\alpha:L^1\to L^{1/(1-\alpha),\infty}$ is bounded.
\end{proposition}
\begin{proof}
Split the kernel integral at a dyadic radius $r$. The near integral is at most $C r^\alpha M_df(x)$ by a shell sum. For $p>1$, H\"older bounds the far integral by $Cr^{\alpha-1/p}\norm f_p$, up to a harmless constant component. Choose the nearest dyadic radius to $\min\{1,(\norm f_p/M_df(x))^p\}$, treating $M_df(x)=0$ separately. The estimate becomes
\[
 |I_\alpha f(x)|\lesssim\norm f_p^{\alpha p}(M_df(x))^{1-\alpha p}+\norm f_p.
\]
The strong $L^p$ bound of $M_d$ gives the assertion. At $p=1$, use the far bound $Cr^{\alpha-1}\norm f_1$, giving the same inequality with $\alpha p=\alpha$. Now the weak $(1,1)$ bound of $M_d$ yields the stated weak exponent. This explicit kernel step is needed; the fractional weak endpoint does not follow from weak $(1,1)$ by an unexplained substitution.
\end{proof}
Numerical-frequency fractional powers, dyadic-scale potentials, and Boolean noise operators need not have the same kernels or domains. Tauberian and Wiener-type questions for dyadic transforms \cite{GolubovWiener03} provide a further connection to noncompact analysis, considered in Section~\ref{sec:noncompact}.

\subsection{Closed domains and a positive difference form}\label{sub:v2calculus}
\begin{proposition}\label{v2:derivative}
The closure of $\mathcal D$ on Walsh polynomials in $L^2$ is positive self-adjoint, with domain $\sum n^2|\widehat f(n)|^2<\infty$. Its inverse on the mean-zero subspace is compact. For every $s>0$, the graph norm of $\mathcal D^s$ is equivalent to the martingale $B^s_{2,2}$ norm.
\end{proposition}
These are direct diagonal spectral facts: Fourier transformation gives multiplication by $n$ on $\ell^2$, polynomials are a graph core, and $n\asymp2^j$ on the $j$th block. The comparison does not assert boundedness of arbitrary within-block symbols in every Banach norm.

\begin{proposition}\label{v2:form}
On polynomials and then on the closed form domain,
\[
 \langle\mathcal Df,f\rangle
 =\sum_{n\ge1}n|\widehat f(n)|^2
 =\sum_{j\ge0}2^{j-2}\|f-T_{e_j}f\|_2^2.
\]
\end{proposition}
Indeed $|1-(-1)^{n_j}|^2=4n_j$ and $n=\sum_j2^jn_j$; Parseval and monotone convergence give the identity. Although $\mathcal D=\sum_j2^{j-1}(I-T_{e_j})$ on polynomials, it is not a derivation: $w_n^2=1$ contradicts a Leibniz rule for $n>0$. The positive form is the useful replacement.

The numerical generator has compact resolvent and the exact heat trace
\[
 \operatorname{Tr}(P_t)=\sum_{n\ge0}e^{-tn}
 =\frac1{1-e^{-t}}=\prod_{j\ge0}(1+e^{-t2^j})=P_t(0).
\]
In contrast, the Boolean number operator has infinitely many eigenvectors $r_j$ of eigenvalue one. Its resolvent is not compact on the infinite product. This multiplicity difference explains why dimension-free hypercontractivity and compact numerical smoothing answer different questions even though the eigenfunctions coincide.

\subsection{Compactness between different Lebesgue norms}\label{sub:v2cross}
\begin{proposition}\label{v2:crosscompact}
If $2\le p<\infty$, then $T_a:L^p\to L^2$ is compact exactly when $a\in c_0$. If $1\le q\le2$, the same equivalence holds for $T_a:L^2\to L^q$.
\end{proposition}
\begin{proof}
On the probability space, $\norm f_2\le\norm f_p$ and $\norm g_q\le\norm g_2$. Hence the norm of either truncated tail is at most $\sup_{n\ge N}|a_n|$; this gives finite-rank approximation for $a\in c_0$. Necessity for an $L^2$ target follows from orthogonality of $a_nw_n$. For an $L^q$ target, suppose $|a_n|\ge\delta$ on infinitely many indices. For distinct $n,k$, the characters coincide on half of the group and are opposite on the other half, so
\[
 \norm{a_nw_n-a_kw_k}_q^q
 =\tfrac12(|a_n-a_k|^q+|a_n+a_k|^q)\ge c_q\delta^q.
\]
Their images have no Cauchy subsequence, contradicting compactness.
\end{proof}
\begin{proposition}\label{v2:HYcompact}
Let $1<p\le2\le q<\infty$, $p<q$, and $1/r=1/p-1/q$. If $a\in\ell^r$, then
\[
 \norm{T_af}_q\le\norm a_r\norm f_p,
\]
and $T_a:L^p\to L^q$ is compact.
\end{proposition}
\begin{proof}
Hausdorff--Young gives $\norm{\wh f}_{p'}\le\norm f_p$. Sequence H\"older gives $\norm{a\wh f}_{q'}\le\norm a_r\norm{\wh f}_{p'}$, since $1/q'=1/r+1/p'$. The valid synthesis inequality is $\ell^{q'}\to L^q$ for $q\ge2$. Applying it proves the estimate. The $\ell^r$ norms of symbol tails tend to zero, so finite-rank truncations converge in operator norm.
\end{proof}
This is a useful positive counterpart to the invalid $\ell^{p'}\to L^p$ shortcut below two. It supplies a checkable cross-space estimate without deciding a general same-space compactness claim.

\subsection{A positive fractional kernel and its normalization}\label{sub:v2positivepotential}
For $0<\alpha<1$, let
\[
 \mathcal K_\alpha=(1-2^{-\alpha})\sum_{m\ge0}2^{-m\alpha}D_{2^m},
 \qquad J_\alpha f=\mathcal K_\alpha*f.
\]
All summands are nonnegative and have unit integral before their scalar weights. Thus $\mathcal K_\alpha$ is an $L^1$ probability kernel. Its symbol is one at zero and $2^{-\alpha(j+1)}$ on $[2^j,2^{j+1})$. With the earlier normalization,
$J_\alpha=E_0+2^{-\alpha}(I_\alpha-E_0)$, equivalently $I_\alpha=E_0+2^\alpha(J_\alpha-E_0)$.
If the first nonzero digit of $x$ is $k$, then
\[
 \mathcal K_\alpha(x)=(1-2^{-\alpha})\sum_{m=0}^k2^{m(1-\alpha)}
 \asymp_\alpha\rho(x,0)^{\alpha-1}.
\]
This proves the kernel majorant used in Proposition~\ref{prop:fractional} for both normalizations. At $\alpha=1$ the growth is logarithmic, and at $\alpha>1$ the kernel is bounded; the same exponent formula is not extended past those changes of regime.

\begin{proposition}[A directly summable variation class]\label{v2:BV}
Let $1<p<\infty$ and $\sum_{n\ge0}|a_{n+1}-a_n|<\infty$, with $a_n\to L$. Then $T_a$ is bounded and $T_a-LI$ is compact. In particular $T_a$ is compact if and only if $L=0$.
\end{proposition}
The detailed derivation is recorded in Appendix~\ref{app:auxiliary-proofs}, preserving the normalization used here.
This global bounded-variation condition is stronger than a uniform variation bound on separate dyadic blocks. For the blockwise class, Proposition~\ref{prop:marccompact} proves that scalar decay alone is enough for compactness. Global variation gives the stronger explicit absolutely convergent operator series in Appendix~\ref{app:auxiliary-proofs}, subsection C.3.

\subsection{Beyond a regular symbol class}\label{sub:multiplier-models}
The same tail argument identifies the limit of an interpolation repair. If $T_a$ is known to be bounded only on $L^{p_0}$ and $a\in c_0$, then $R_M$ is uniformly bounded on $L^{p_0}$ and tends to zero in operator norm on $L^2$. Consequently it tends to zero on every strictly intermediate $L^q$. The endpoint $q=p_0$ remains outside this argument. Proposition~\ref{prop:marccompact} reaches it because \eqref{eq:marcink} supplies a farther exponent.

The general tail criterion also yields, for a bounded $T_a$ on $1\le p<\infty$,
\[
 \|T_a\|_{\rm ess}\le\liminf_M\|R_M\|_{p\to p}
 \le\limsup_M\|R_M\|_{p\to p}\le2\|T_a\|_{\rm ess}.
\]
The first inequality uses finite rank. For the last subtract any compact $K$, use uniform convergence of $E_M$ on the compact image of the unit ball under $K$, and then minimize $2\|T_a-K\|$. For $1\le p<\infty$, $\limsup|a_n|\le\|T_a\|_{\rm ess}$ by weak nullity. At $p=1$, each dual test belongs to $L^\infty\subset L^2$, so its Walsh coefficients tend to zero; hence $w_n\rightharpoonup0$ in $L^1$ as well. Equality with this scalar limsup has been proved here only in the Hilbert case.

There is an analogous distinction for spectra. Scalar inversion gives the $L^2$ spectral formula; on $L^p$ one must also show that the reciprocal symbol is a bounded multiplier. Spectral mapping for multipliers on general compact abelian groups is a separate problem \cite{Ricker05}; this reference is context, not a claimed Walsh counterexample. The symbol, its multiplier norm, and the ambient operator algebra must all be fixed.

\paragraph{Comparison and research boundary.}
The exact Hilbert theory and the regular-symbol compactness consequence leave a precise operator question: what extra information, if any, is needed at a single non-Hilbert exponent? A proof using a second interpolation endpoint cannot answer that single-exponent question without supplying the endpoint.

\section{Greedy approximation and correction}\label{sec:greedy}
The spectral support is now selected from the original coefficients. The main comparison is between thresholding and free best sparse approximation, followed by the change produced by correcting the input function. Equal-size character families exhibit the sharp loss and explain the exceptional role of the Hilbert norm.

Allowing a selectable spectral support changes approximation even before one changes the function. The Dirichlet and independent-character tests quantify the resulting Lebesgue loss. Correction then changes the input on a small set to impose convergence and coefficient order; its common-set variants have a stronger quantifier pattern.

\subsection{Thresholding, conditionality, and the exact loss}
For $f\in L^p$, a greedy set $A$ of size $N$ consists of indices of $N$ largest coefficients in modulus. Fix increasing-frequency tie breaking. Put
\[
 G_Nf=P_Af=\sum_{n\in A}\wh f(n)w_n,
 \qquad e_N(f)_p=\inf_{|B|\le N}\inf_{(b_n)}\norm{f-\sum_{n\in B}b_nw_n}_p.
\]
The best approximant may change the coefficients, whereas $G_N$ retains Fourier coefficients. For $p=2$, Parseval makes $G_N$ optimal. For $p\ne2$, the Walsh system is not quasi-greedy. In $L^1$ it is not even a Schauder basis in its natural order, so statements about thresholding operators are more precise than indiscriminate basis terminology.

The general framework was developed around \cite{KonyaginTemlyakov,Wojtaszczyk}, with extensive subsequent theory \cite{Temlyakov11,GreedyReview25}. The constant-coefficient unconditionality forced by quasi-greediness \cite{Wojtaszczyk} gives a direct test for a uniformly bounded orthonormal system. In the Walsh setting, the contrast
\[
 \norm{D_N}_p=N^{1-1/p}\quad(N=2^m),\qquad
 \norm{r_1+\cdots+r_N}_p\asymp_p\sqrt N
\]
exhibits the underlying norm geometry. Complete martingale blocks still form an unconditional decomposition for $1<p<\infty$. Selecting whole blocks and selecting individual coefficients are different algorithms, with different costs.

\begin{theorem}\label{th:greedyrate}
For $1\le p<\infty$,
\[
 \norm{f-G_Nf}_p\le C_pN^{|1/p-1/2|}e_N(f)_p.
\]
The power is optimal uniformly in $f,N$; at $p=2$ the constant is one.
\end{theorem}
\begin{proof}
For every finite $A$, $|A|\le N$,
\begin{equation}\label{eq:projection}
 \norm{P_Af}_p\le N^{|1/p-1/2|}\norm f_p.
\end{equation}
For $p\ge2$, interpolate $\norm P_\infty\le\sqrt N\norm P_2$ and use $\norm{P_Af}_2\le\norm f_p$. For $p\le2$, use Hausdorff--Young followed by finite-dimensional H\"older to bound the coefficient $\ell^2$ norm by $N^{1/p-1/2}\norm f_p$.

Given an approximant $g$ with at most $N$ frequencies, enlarge its support set $B$ to cardinality $N$ by zero coefficients if necessary, and write $h=f-g$, $C=B\setminus A$, $D=A\setminus B$. Greediness and $|C|=|D|$ give
\[
 \norm{(\wh f(n))_{C}}_2\le\norm{(\wh h(n))_{D}}_2.
\]
The preceding coefficient and synthesis estimates imply $\norm{P_Cf}_p\le N^{|1/p-1/2|}\norm h_p$. Now
\[
 f-P_Af=h-P_Ah-P_Ch+P_Cf.
\]
Use \eqref{eq:projection} and take the infimum over $g$.

For the lower bound let $N=2^m$ and $R_N=\sum_{j=0}^{N-1}w_{2^{m+j}}$. Its spectrum is disjoint from that of $D_N$. For $p>2$, use $f=D_N+(1+\eps)R_N$: thresholding keeps $R_N$, whereas choosing $D_N$ gives an admissible error of order $\sqrt N$. Their ratio is $\gtrsim_p N^{1/2-1/p}$. For $p<2$, use $(1+\eps)D_N+R_N$ and interchange the roles. Positive $\eps$ eliminates a tie between the two groups.
\end{proof}
For $L^2$ coefficient sparsity, if $c_k^*\le Ak^{-1/r}$ with $r<2$, then $\norm{f-G_Nf}_2\lesssim_r A N^{1/2-1/r}$. This Hilbert-space fact does not describe all $L^p$ errors, since equal coefficients can create either a localized Dirichlet kernel or a randomized polynomial with a very different spatial norm.

\subsection{Where Walsh thresholding sits in greedy basis theory}\label{sub:greedy-comparison}
The approximation problem here is measured against free $N$-term error $e_N(f)_p$, whereas coordinate error restricts the approximant to the original coefficients. These two benchmarks lead to different notions of greedy and almost greedy bases. For a normalized Schauder basis of a Banach space, the classical characterizations recalled in \cite{GreedyReview25} are
\[\begin{aligned}
 \text{greedy}&\ \Longleftrightarrow\ \text{unconditional and democratic},\\
 \text{almost greedy}&\ \Longleftrightarrow\ \text{quasi-greedy and democratic}.
\end{aligned}\]
Democracy compares norms of sums of equally many normalized basis vectors; quasi-greediness controls greedy projections uniformly. Its weaker condition does not imply democracy for an arbitrary basis. These distinctions explain why the direct uniformly bounded-orthonormal argument in Proposition~\ref{v2:notquasi} is needed.

\begin{longtable}{@{}p{.25\textwidth}p{.29\textwidth}p{.33\textwidth}@{}}
\caption{Coefficient thresholding, locality, and the comparison benchmark}\label{tab:greedy-comparison}\\
\toprule Setting & Established conclusion & Mathematical reason\\\midrule\endfirsthead
\toprule Setting & Established conclusion & Mathematical reason\\\midrule\endhead
Walsh in $L^2$ & Thresholding is best $N$-term approximation. & Orthogonality makes error the discarded square sum.\\
Walsh in $L^p$, $1<p<\infty$, $p\ne2$ & TGA is not uniformly bounded; Theorem~\ref{th:greedyrate} gives the sharp power loss against $e_N$. & Equal coefficient magnitudes can produce localized kernels or Rademacher sums with different spatial norms.\\
$L^p$-normalized one-dimensional Haar, $1<p<\infty$ & The Haar basis is greedy; see \cite{GreedyReview25,TemlyakovActa08}. & Its unconditional and democratic structure differs from the Walsh character basis.\\
Walsh after correction & Convergence can hold for the corrected input. & The input function, and sometimes its agreement set, changes.\\
Adaptive approximation with reoptimized coefficients & A different algorithm and a different guarantee must be specified. & Selection uses residual information and is not simply retention of the largest original coefficients.\\\bottomrule
\end{longtable}

For example, the Weak Chebyshev Greedy Algorithm chooses an atom correlated with a norming functional of the current residual and then recomputes a best approximant in the selected span; \cite{TemlyakovActa08,GreedyReview25} give the precise Banach-space definitions. Its performance cannot be read from a negative TGA theorem, and results for a trigonometric dictionary do not become Walsh results by replacing its name. The useful comparison is between algorithms on a fixed dictionary, or between fixed algorithms on different dictionaries, with the same error benchmark and stated parameter dependence.

\paragraph{Local research boundary.}
The sharp TGA power and the Hilbert boundary below are settled comparisons. A proposed improvement must change an identified assumption, dictionary, input class, or algorithm. Correction and universality are then connected to greedy theory by their changed quantifiers rather than presented as a contradiction to an unconditional-basis obstruction.

\begin{longtable}{@{}p{.19\textwidth}p{.30\textwidth}p{.38\textwidth}@{}}
\caption{Which structural comparisons survive a change of system}\label{tab:systems-transfer}\\
\toprule System & Retained mechanism & Boundary of transfer\\\midrule\endfirsthead
\toprule System & Retained mechanism & Boundary of transfer\\\midrule\endhead
Binary Walsh & Characters, XOR modulation, regular binary martingale blocks. & Natural-order $L^p$ basis behavior does not imply quasi-greediness; Proposition~\ref{v2:notquasi}.\\
$L^p$-normalized one-dimensional Haar & Local support and unconditional, democratic basis structure for $1<p<\infty$. & Greediness uses this normalization; Haar functions are not uniformly bounded characters.\\
Trigonometric & Character modulation, orthogonality, and time--frequency methods. & Frequency addition replaces XOR; dyadic kernel identities and Walsh endpoint proofs do not transfer verbatim.\\
Bounded Vilenkin & Character algebra and uniformly regular finite-branching martingales. & Constants may depend on the digit bound; binary index formulas must be replaced by mixed-radix formulas.\\
Unbounded Vilenkin & Character orthogonality and complete scale-block projections. & Uniform branching estimates and a bounded overshoot cost are unavailable without further hypotheses; see Sections~\ref{sec:universal} and \ref{sub:ordinary-vilenkin}.\\\bottomrule
\end{longtable}

\subsection{One integrable function with divergent greedy sums}\label{subsec:badgreedy}
An explicit construction explains more than unbounded Lebesgue constants. It combines small full blocks, large intermediate sums, and a tiny perturbation fixing their greedy order. This is the mechanism in Episkoposian's divergence work \cite{Epi07,EpiGeneral07}; it is included as an expository binary realization, not a new result.

Put $v_q=1+4+\cdots+4^{q-1}$. The recursion \eqref{eq:lebesgue} gives $L_{v_q}\ge1+2(q-1)/3$. Choose strictly increasing $q_j\ge j2^j$ and set
\[
 B_j=2^{2q_j},\qquad a_j=2^{-j},\qquad \delta_j=2^{-2j-3},
\]
\[
 P_j=w_{B_j}\Big(a_jD_{B_j}+\delta_j\sum_{i=0}^{B_j-1}2^{-i-1}w_i\Big).
\]
The spectrum is $[B_j,2B_j)$; these blocks are disjoint. Its nonzero coefficients, in frequency order, are $a_j+\delta_j2^{-i-1}$, strictly decreasing within each block and across successive blocks. Also $\norm{P_j}_1\le a_j+\delta_j$. Therefore $f=\sum_jP_j$ converges in $L^1$, and continuity of each Fourier coefficient functional identifies precisely these coefficients.

Set $m_j=v_{q_j}<B_j$ and $M_j=\sum_{h<j}B_h+m_j$. All earlier blocks, followed by the first $m_j$ coefficients of block $j$, are exactly the $M_j$ largest nonzero coefficients of $f$. Thus
\[
 \norm{G_{M_j}f}_1
 \ge a_jL_{m_j}-\delta_j-\sum_{h<j}\norm{P_h}_1\longrightarrow\infty.
\]
Indeed, the first term grows at least linearly in $j$ and the last sum is bounded. This proves divergence for one function with strictly decreasing nonzero coefficients. The same construction works in the explicitly chosen union-of-blocks subsystem. It is not a theorem for every subsystem: independent Rademacher characters are equivalent to the $\ell^2$ basis in every finite $L^p$. Special positive subsystem results belong to \cite{GrigZink06}.

\subsection{A local corrector and the global assembly}
Correction allows one to change a function on a set of arbitrarily small measure. The task is to obtain convergence and coefficient order simultaneously; Lusin approximation alone does not do this.

Let $0<\eps<1$, let $\Delta$ be a rank-$m$ cylinder, and let $\gamma\ne0$. Choose $q$ with $2^{-q}<\eps$, $2^q\le2/\eps$, and $L\ge m$ above a prescribed spectral threshold. Let $B$ be the event $x_L=\cdots=x_{L+q-1}=0$. The polynomial
\begin{equation}\label{eq:localcorrector}
 Q=\gamma\ind_\Delta(1-2^q\ind_B)
\end{equation}
equals $\gamma$ on $\Delta\setminus B$, whose relative measure exceeds $1-\eps$, and vanishes outside $\Delta$. Its spectrum lies above $2^L-1$, because
\[
 2^q\ind_B=\sum_{\nu=0}^{2^q-1}w_{2^L\nu}
\]
and the constant component cancels. Multiplication by $\ind_\Delta$ only changes the lower $m$ digits. All nonzero coefficients have modulus $|\gamma|\mu(\Delta)$, and
\[
 \norm Q_1=2|\gamma|\mu(\Delta)(1-2^{-q}),\qquad
 \norm Q_2^2=|\gamma|^2\mu(\Delta)(2^q-1).
\]
Hence $\sup_N\norm{S_NQ}_1\le\norm Q_2\le|\gamma|\sqrt{2\mu(\Delta)/\eps}$. Notice also $\norm{Q-\gamma\ind_\Delta}_1=|\gamma|\mu(\Delta)$: equality on a large set does not imply small unweighted error.

The local construction explains concentration, but organizing different heights into decreasing coefficients needs a further lemma. The following is the block result used from \cite{Epi12}, Lemma 2.2 of its author version arXiv:1407.1496.
\begin{lemma}[Correction block]\label{lem:correctionblock}
Given $h\in L^1$, $h\ne0$, a starting frequency, and $\eta>0$, there are a function $g$, a set $E_h$ of measure greater than $1-\eta$, and a polynomial $P$ supported above that frequency such that
\[
 g=h\text{ on }E_h,\quad\norm g_1\le3\norm h_1,\quad\norm{P-g}_1<\eta.
\]
The nonzero coefficients of $P$ are nonincreasing in modulus and smaller than $\eta$, and every initial partial sum of $P$ has $L^1$ norm at most $3\norm h_1$.
\end{lemma}
The published lemma is stated for generalized Walsh systems of a fixed integer base. We use its binary specialization. Its proof refines cylinders, applies localized correctors, and matches coefficient sizes. Formula \eqref{eq:localcorrector} alone does not prove the full ordering assertion.

\begin{theorem}[Episkoposian's correction theorem \cite{Epi12}]\label{th:correction}
For every $f\in L^1$ and $\eps>0$, one can find $F\in L^1$ equal to $f$ outside a set of measure less than $\eps$, whose Fourier--Walsh series converges in $L^1$ and whose nonzero coefficient moduli are nonincreasing. Its greedy sums, with the stated tie convention, converge in $L^1$.
\end{theorem}
\begin{proof}[Assembly from Lemma~\ref{lem:correctionblock}]
Choose cylinder functions $f_j$ with $\norm{f_j-f}_1\le2^{-j}$ and put $f_0=r_0=0$. At step $j$ apply the lemma to
\[
 h_j=f_j-f_{j-1}+r_{j-1},\qquad r_j=g_j-P_j.
\]
Place the new spectrum after every preceding block. Choose $\eta_j\le2^{-2j}$ and $\eta_j\le\eps2^{-j-1}$, also below the smallest nonzero coefficient modulus of the last nonempty block. For $h_j=0$, set $g_j=P_j=r_j=0$ and $E_{h_j}=\D$. This ensures global coefficient order and $\norm{r_j}_1<\eta_j$.

Indeed,
\[
 \sum_j\|h_j\|_1\le\sum_j\|f_j-f_{j-1}\|_1+\sum_j\|r_{j-1}\|_1<\infty,
 \qquad \|P_j\|_1\le3\|h_j\|_1+\eta_j.
\]
Thus the series $\sum_jP_j$ converges absolutely in $L^1$ to a function $F$. On $E=\bigcap_jE_{h_j}$, of measure greater than $1-\eps$, the identity
\[
 P_j=f_j-f_{j-1}+r_{j-1}-r_j
\]
telescopes to
\[
 \sum_{i=1}^jP_i=f_j-r_j.
\]
Passing to the limit in $L^1(E)$ proves $F=f$ there. At block endpoints the Fourier series converges in $L^1$. The norm of any intermediate contribution from block $j$ is at most $3\norm{h_j}_1\to0$, so all partial sums converge. Disjoint spectra and continuity of coefficient functionals identify the constructed coefficients with those of $F$. Nonincreasing nonzero moduli now identify greedy sums with initial sums on that support.
\end{proof}
The proof establishes $\forall f\,\exists E_f,F_f$. The stronger variants in the source that select a common large set before all target functions require a separate universal construction.

\subsection{Grigoryan's line: common sets, profiles, and different topologies}
The correction programme includes uniform convergence \cite{GriUniform10,EpiUniform13}, $L^p$ convergence for generalized Walsh systems \cite{EpiGri12}, and nonlinear approximation extensions \cite{EpiTigran25,EpiGreedy26}. A particularly informative strengthening concerns a common modulus profile.
\begin{theorem}[Grigoryan--Navasardyan \cite{GriNava16}]\label{th:commonprofile}
For every $\eps>0$, there exist a measurable $E\subset[0,1)$, $|E|>1-\eps$, and $g\in L^1$ with positive decreasing Walsh coefficients such that for each $f\in L^1$ there is $\widetilde f\in L^1$ equal to $f$ on $E$, whose Fourier series converges in $L^1$ and satisfies
\[
 |\wh{\widetilde f}(n)|=\wh g(n)\quad\text{for every }n.
\]
\end{theorem}
The set and modulus profile precede the target. The signs may depend on it; imposing equality of the signed coefficients would force a single Fourier object and destroy the conclusion. Monotone-coefficient representation results \cite{GevNava99,Nava07} form part of this background.

These statements explain how correction coexists with the negative greedy theorem. Correction changes the input function. It does not assert uniform boundedness of $G_N$ on the original $L^p$ space, nor does it contradict the divergent function in Section~\ref{subsec:badgreedy}. Recent greedy-universal functions \cite{EGS22,ESG26} combine selection and approximation under additional definitions; the adjective ``universal'' must be accompanied by the allowed targets and convergence mode.

\subsection{Quasi-greediness and the Hilbert-space boundary}\label{sub:v2quasigreedy}
The negative conclusion needs more than observing a failure of democracy: in general a quasi-greedy basis need not be democratic. The uniformly bounded orthonormal setting supplies the extra argument.
\begin{proposition}\label{v2:notquasi}
For $1\le p<\infty$, $p\ne2$, the Walsh thresholding operators do not form a uniformly bounded family on $L^p$. For $p=2$ they give best $N$-term approximations.
\end{proposition}
\begin{proof}
Suppose their norms were bounded by $C$. For a finite set $A$, first give the coefficients on a subset $B\subset A$ modulus $1+\eta$ and all the others modulus one. Thresholding then selects $B$; let $\eta\downarrow0$ to obtain uniformly bounded coordinate suppression for polynomials having coefficient moduli one. A signed polynomial $\sum_{n\in A}\eps_nw_n$ consequently has norm comparable, uniformly in $A$ and the signs, to $\sum_{n\in A}w_n$: split into the positive and negative sets for one direction and reverse their signs for the other. Randomize $\eps_n$. Pointwise Khintchine, followed by integration, gives
\[
 \mathbb E_\eps\norm{\sum_{n\in A}\eps_nw_n}_p^p\asymp_p|A|^{p/2}.
\]
Thus every unsigned polynomial with $N$ coefficients one would have norm comparable to $\sqrt N$. But for $A=[0,2^m)$ it is $D_{2^m}$, with norm $N^{1-1/p}$. These powers agree only for $p=2$. At $p=2$, Parseval makes deletion of the smallest coefficients optimal and ensures convergence.
\end{proof}
This proof applies the constant-coefficient unconditionality forced by quasi-greediness, as developed in \cite{Wojtaszczyk}. The two equal-size families $\{w_n:n<N\}$ and $\{r_{m+j}:0\le j<N\}$ separately show non-democracy for $p\ne2$. Replacing the latter by a sum of merely $m$ characters would not compare sets of cardinality $N=2^m$.

Uniform correction results \cite{GriUniform10,EpiUniform13} require their own target classes and topology; they are not obtained by replacing $\norm{\cdot}_1$ in Lemma~\ref{lem:correctionblock} with a supremum. Similarly, the fixed-set universal-profile theorem \ref{th:commonprofile} makes the large set independent of the target, whereas the inductive proof of Theorem~\ref{th:correction} does not. These distinctions preserve the strength of the published results without manufacturing an unproved simultaneous strengthening.

\paragraph{Comparison and research boundary.}
Greedy convergence after correction changes the input, while a universal representation changes choices after the target. The sharp negative theorem for unchanged Walsh thresholding is therefore compatible with both mechanisms. The next section compares their quantifiers rather than placing them on one undifferentiated scale of strength.

\section{Correction, representation, and universality}\label{sec:universal}
Representation permits a choice after the target. Its strength depends on whether one selects stopping indices, deletes terms, changes signs, or uses every term in a new order. This section follows those operations separately, adding Fourier origin, common exceptional sets, coefficient restrictions, and weighted topology only when the construction supports them.

Universality requires the representation to be chosen after the target. The choice may be a subsequence of ordinary partial sums, signs, deletion, or a permutation, and the target topology determines which coefficient obstructions survive. Constructive blocks and category arguments are developed separately before their quantifiers are compared.

\subsection{The allowed operation is part of the theorem}
For a fixed coefficient sequence $c$, several notions of universality must be separated. A subseries uses $\delta_n\in\{0,1\}$ in natural order. A signed series uses $\delta_n\in\{-1,1\}$. Ternary selection uses $\delta_n\in\{0,\pm1\}$. A rearrangement uses every term exactly once. Pointwise universality of a function instead selects a subsequence of its ordinary Fourier partial sums. A weighted universal series need not be the Fourier series of an unweighted integrable function.

These distinctions are central to the Armenian representation programme and its international extensions \cite{GriEpi01,EpiWalsh99,EpiUniversal16,Epi13Double,GGS,Gri20,Gri24,UniversalPair23,Sargsyan25}. Related trigonometric constructions \cite{Epi06} are relevant as methods, but are not Walsh theorems by analogy. More abstract conditional universality for minimal systems is developed by Avetisyan--Grigoryan--Ruzhansky \cite{Avetisyan26}; the phase modification and the ambient sequential convergence space belong to its formulation.

The continuous-coefficient obstruction is recorded in Section~\ref{sub:universal-comparison}; weighted topologies and null-set convergence evade it by different mechanisms.

\subsection{Episkoposian--M\"uller pointwise universality}\label{subsec:muller}
\begin{theorem}[Episkoposian--M\"uller \cite{EpiMuller14}, Theorem 1.1]\label{th:muller}
For every fixed countable $E\subset\D$, the functions $f\in C(\D)$ such that, for every $h:E\to\R$, some subsequence $S_{n_j}f$ converges pointwise on $E$ to $h$, form a residual subset of the real space $C(\D)$.
\end{theorem}
The order is $E$ fixed, then a residual class of $f$, then arbitrary $h$, then a subsequence depending on $h$. There is no assertion of universality on all countable sets simultaneously.

\begin{lemma}[Finite control]\label{lem:finitecontrol}
Let $x_1,\ldots,x_s$ be distinct, $P$ a Walsh polynomial, $b_1,\ldots,b_s\in\R$, and $\eps>0$. For arbitrarily large $n$ there is a Walsh polynomial $Q$ with $\norm{Q-P}_\infty<\eps$ and $S_nQ(x_i)=b_i$ for every $i$.
\end{lemma}
\begin{proof}
Choose $m$ so that the cylinders $I_m(x_i)$ are disjoint. Set $n=2^mq$, with $q$ sufficiently large and $L_q>\max_i|b_i-P(x_i)|/\eps$. By \eqref{eq:dilate}, $D_n$ is supported on $I_m$ and has $L^1$ norm $L_q$. Put
\[
 H_i(t)=\frac{b_i-P(x_i)}{L_q}\ind_{I_m(x_i)}(t)
       \operatorname{sgn}D_n(x_i\oplus t).
\]
These are cylinder functions and hence Walsh polynomials. Their supports are disjoint, so $Q=P+\sum_iH_i$ is within $\eps$ of $P$. Convolution gives $S_nH_i(x_k)=0$ for $i\ne k$ and $b_i-P(x_i)$ for $i=k$. Take $n$ above the spectrum of $P$.
\end{proof}
\begin{proof}[Proof of Theorem~\ref{th:muller}]
If $E$ is empty the assertion is vacuous. If $E$ is infinite, enumerate $E=\{x_1,x_2,\ldots\}$; if $E$ is finite use its whole list and keep $s=|E|$ in the argument below. For a rational vector $b\in\mathbb Q^s$ and integers $k,N\ge1$, the set
\[
 \bigcup_{n\ge N}\{f:\max_{i\le s}|S_nf(x_i)-b_i|<1/k\}
\]
is open, since the finite coefficient maps are continuous, and dense by Lemma~\ref{lem:finitecontrol} and polynomial density. Baire's theorem gives a residual intersection over all these countably many choices. Approximate successively more values of an arbitrary $h$ by rationals (all values when $E$ is finite) and choose increasing indices with error tending to zero. This gives the required pointwise convergence.
\end{proof}
The finite construction also provides an exact control constant for the special indices: the smallest uniform norm of a continuous $H$ satisfying $S_{2^mq}H(x_i)=b_i$ is $\max_i|b_i|/L_q$. The upper estimate is the disjoint-cylinder construction, and the lower follows from $|S_nH(x_i)|\le L_n\norm H_\infty$.

This finite control and the countable Baire intersection prove the fixed-set assertion.

\subsection{A high-frequency block for weighted universality}
We give an expository block construction for the weighted mechanisms of \cite{GriEpi01,EpiWalsh99,EpiUniversal16}. The precise combined formulation in Theorem~\ref{th:weighteduniv}, including the prescribed envelope, is proved here from Lemma~\ref{lem:univblock}; it is not attributed verbatim to one of those papers and no priority claim is made for that combination. The distinction between equality on a large set and small unweighted norm is decisive.
\begin{lemma}\label{lem:univblock}
Let $P$ have spectrum in $[0,2^m)$, $\eps,\eta,\delta>0$, and let finitely many exponents $q_i>2$ and tolerances $b_i>0$ be prescribed. Above any starting frequency there is a polynomial $Q=\sum c_nw_n$ and a set $E$, $\mu(E)>1-\eps$, such that $Q=P$ on $E$,
\[
 \max_n|c_n|<\delta,\qquad \sum_n|c_n|^{q_i}<b_i,
\]
and for every measurable $e\subset E$ and every partial sum,
\begin{equation}\label{eq:univlocal}
 \int_e|S_NQ|\dd\le\int_e|P|\dd+\eta.
\end{equation}
\end{lemma}
\begin{proof}
Let $B_{L,s}=\prod_{j=L}^{L+s-1}(1+r_j)$. Choose a large $K$, let $s=\lceil\log_2(2K/\eps)\rceil$, and choose $K$ disjoint increasingly ordered coordinate groups $[L_i,L_i+s)$ above $m$ and the required threshold. Define
\[
 Q=P\Big(1-\frac1K\sum_{i=1}^KB_{L_i,s}\Big)
   =-\frac PK\sum_{i=1}^K(B_{L_i,s}-1).
\]
On $E=\bigcap_i\{B_{L_i,s}=0\}$, $Q=P$, and $\mu(E^c)\le K2^{-s}\le\eps/2$. The constant high-coordinate components cancel. As low and high digits are disjoint, each coefficient $-\wh P(a)/K$ is repeated exactly $K(2^s-1)$ times at distinct frequencies. Thus
\[
 \sum_n|c_n|^q=K^{1-q}(2^s-1)\sum_{a<2^m}|\wh P(a)|^q
 \le\frac4\eps K^{2-q}\sum_a|\wh P(a)|^q.
\]
This tends to zero for every prescribed $q>2$, and the largest coefficient also tends to zero.

In natural order a partial sum contains several completed coordinate groups and at most one unfinished group. On $E$ the completed part is $(\nu/K)P$, $0\le\nu\le K$. The unfinished part has $L^1$ norm at most
\[
 \frac{s+m+2}{K}\norm P_1.
\]
To verify this, enumerate high-digit frequencies first and the lower $m$ digits inside each. The two possible remainders are a multiple of $P(D_b^{[L_i]}-1)/K$ and a modulated $S_tP/K$. Independence gives the first norm bound $(s+1)\norm P_1/K$; the elementary Lebesgue bound gives $(m+1)\norm P_1/K$ for the second. If $N$ is between groups, the unfinished term is absent. Since $s=O(\log K)$, enlarge $K$ to make the remainder smaller than $\eta$. Integrating on any $e\subset E$ proves \eqref{eq:univlocal}.
\end{proof}
The factor $K^{2-q}$ identifies the threshold in this construction. It is not legitimate to set $q=2$. For $P=1$ the polynomial $Q$ has integral zero, so $\norm{Q-1}_1\ge1$ despite equality on most of the group.

\begin{theorem}[Constructive weighted representation]\label{th:weighteduniv}
Given a positive nondecreasing $\Lambda_N\to\infty$, there is a formal Walsh series $\sum c_nw_n$ with
\[
 c\in\bigcap_{q>2}\ell^q,\qquad
 \abs{\sum_{n<N}c_nw_n(x)}\le\Lambda_N\quad(x\in\D,N\ge1),
\]
such that, for every $\eps>0$, some weight $0<v_\eps\le1$ with $\mu\{v_\eps\ne1\}<\eps$ makes this series universal in $L^1(v_\eps\dd)$ both by subseries and by rearrangements. Explicitly, the order of choices is
\[
 \Lambda\ \longmapsto\ c\ \longmapsto\ \eps\ \longmapsto\ v_\eps\ \longmapsto\ f\ \longmapsto\ \text{selection or rearrangement}.
\]
The coefficient series depends on $\Lambda$ and is independent of $\eps$; the weight may depend on $\eps$, and is fixed before every target in its space.
\end{theorem}
\begin{proof}
Enumerate all rational complex Walsh polynomials as $P_j$, each infinitely often. Apply Lemma~\ref{lem:univblock} in successively disjoint frequency blocks, with $\mu(E_j^c)<2^{-j}$, error $2^{-j}$, largest coefficient below $2^{-j}$, and
\[
 \sum_{n\in\spec Q_j}|\wh Q_j(n)|^{2+1/j}<2^{-j}.
\]
After choosing a block, its coefficient $\ell^1$ norm is known and unaffected by raising its high coordinate positions. Raise them until $\Lambda_N$ exceeds the total absolute coefficient mass of this block and every preceding block. Fill all gaps with zero coefficients. This proves the pointwise envelope. For each fixed $q>2$, all sufficiently late exponents $2+1/j$ are at most $q$, and the coefficient moduli are at most one. Summation proves $c\in\ell^q$ and $c_n\to0$.

For the rest of the proof fix $\eps$ and abbreviate $v_\eps$ by $v$. Here is a direct common-weight step. Let $F_j$ be the maximum of the absolute values of $P_j$, $Q_j$, and all its distinct partial sums; it is bounded. Choose $J$ with $\sum_{j\ge J}2^{-j}<\eps$. On each $E_j^c$, $j\ge J$, impose
$v\le2^{-2j}/(1+\norm{F_j}_\infty)$, and set $v$ equal to the minimum of one and all applicable bounds. Borel--Cantelli shows that almost every point belongs to only finitely many $E_j^c$, so this minimum is positive almost everywhere. On the null exceptional set define $v$ arbitrarily in $(0,1]$. Then $v=1$ outside a set of measure less than $\eps$, and $\int_{E_j^c}vF_j\dd\le2^{-2j}$. For every $j\ge J$, equality $Q_j=P_j$ on $E_j$ gives
\[
 \norm{Q_j-P_j}_{L^1(v)}\le2^{1-2j},\qquad
 \sup_N\norm{S_NQ_j}_{L^1(v)}\le\norm{P_j}_{L^1(v)}+2^{1-j}.
\]
For the second estimate, layer cake gives the precise common-weight step:
\[
 \int_{E_j}v|S_NQ_j|\dd
 =\int_0^1\int_{E_j\cap\{v>t\}}|S_NQ_j|\dd\,dt
 \le\int_{E_j}v|P_j|\dd+2^{-j}.
\]
Add the controlled bad-set contribution. This explains why local control for every measurable subset was included in the block lemma.

Fix $f\in L^1(v\dd)$ and put $R_0=f$. Rational complex cylinder functions are dense for this finite measure and occur infinitely often in the enumeration. Choose $j_\nu>j_{\nu-1}$, $j_\nu\ge J$, so late that
\[
 \|P_{j_\nu}-R_{\nu-1}\|_{L^1(v)}<2^{-\nu-2},
 \qquad 2^{1-j_\nu}<2^{-\nu-2},
 \qquad R_\nu=R_{\nu-1}-Q_{j_\nu}.
\]
Then $\|R_\nu\|_{L^1(v)}<2^{-\nu-1}$. Moreover
$\sup_N\|S_NQ_{j_\nu}\|_{L^1(v)}\le\|R_{\nu-1}\|_{L^1(v)}+2^{-\nu-1}$.
Endpoint convergence and this estimate prove convergence of every intermediate sum of the selected subseries.

For a rearrangement, at stage $\nu$ first insert the least unused index $k_\nu$. Correct $R_{\nu-1}-c_{k_\nu}w_{k_\nu}$ by an entirely untouched sufficiently late block, chosen with the same two small errors. Only finitely many indices have been used, so infinitely many candidate blocks remain. The indices $k_\nu$ tend to infinity and $\|c_{k_\nu}w_{k_\nu}\|_{L^1(v)}\le|c_{k_\nu}|\to0$. Consequently the corrected residuals and all intermediate errors tend to zero. Taking the least unused index at every stage ensures that every term occurs exactly once.
\end{proof}
The weight precedes the arbitrary target. The coefficients are formal; the natural-order envelope need not survive a target-dependent rearrangement.

\subsection{Ordinary partial sums on common large sets}\label{sub:ordinary-vilenkin}
The preceding representation theorem changes the chosen terms. A different construction leaves every coefficient and its natural position unchanged and selects only the stopping indices. The bounded-Vilenkin result \cite{EpiVilenkin26} belongs to this second class. Let $\mathbf m=(m_1,m_2,\ldots)$, $2\le m_j\le B<\infty$, put $M_0=1$, $M_k=\prod_{j=1}^km_j$, and write $\psi_n$ for its naturally ordered characters. The same notation includes the binary and fixed-integer-base cases.

\begin{theorem}[Episkoposian \cite{EpiVilenkin26}, Theorem 3.3]\label{th:ordinary-vilenkin}
For each fixed $q>2$ there is a formal series $\sum_{n\ge0}c_n\psi_n$, with $c\in\ell^q$, such that
\[\begin{gathered}
 \forall\varepsilon>0\ \exists E_\varepsilon,\quad \mu(E_\varepsilon)<\varepsilon,\\
 \forall p\in[1,\infty)\ \forall g\in L^p\ \exists N_k\uparrow\infty:
 \left\|\sum_{n<N_k}c_n\psi_n-g\right\|_{L^p(E_\varepsilon^c)}\longrightarrow0.
\end{gathered}\]
The exceptional set precedes both $p$ and $g$. The subsequence may depend on them. No signs, deletion, or rearrangement are used.
\end{theorem}
The published article appeared on 24 September 2026. Its quantitative localization lemma reproduces an $\mathbf m$-adic step function outside a set of measure $\theta$, above any specified spectral threshold, with arbitrarily small $\ell^q$ coefficient mass. A finite-group multiplier estimate has the scale
\[
 \sum_{b\ne0}|\beta_b|^q
 <256 B\theta^{-1}\rho^{q-2}
       \log\!\frac{256}{\theta\rho^2},\qquad 0<\rho<1.
\]
Boundedness of the digit bases controls overshoot when a minimal digit block is selected. The constant is an upper estimate for that construction, not an optimal localization cost.

The useful induction targets a repeated dense sequence $h_j$: at stage $j$ localize the residual $h_j-S_{N_{j-1}}c$, not $h_j$ alone. Then $S_{N_j}c=h_j$ off a set $E_j$ with $\sum_j\mu(E_j)<\infty$. For a fixed $\varepsilon$, remove a sufficiently late union of these sets. Every subsequent endpoint is then exactly its assigned dense target on the remaining set. This explains both the ordinary-sum property and the target-independent exceptional set. It is different from summing independently constructed approximants to successive targets.

\begin{corollary}[A diagonal consequence of the localization lemma]\label{cor:allq-vilenkin}
Under the same bounded-digit hypothesis, one can require a single coefficient sequence in $\bigcap_{q>2}\ell^q$ in Theorem~\ref{th:ordinary-vilenkin}.
\end{corollary}
\begin{proof}
At stage $j$ use the source's localization lemma with $q_j=2+1/j$, coefficient maximum below $2^{-j}$, exceptional measure below $2^{-j}$, and $q_j$-mass below $2^{-j}$. The residual is a finite step function, so all lemma inputs are available before choosing its new block. Gaps are filled by zero coefficients. For a fixed $q>2$, eventually $q_j\le q$ and all coefficients have modulus below one, whence the $q$-mass of each late block is at most its $q_j$-mass. The finitely many early blocks do not affect membership in $\ell^q$. The endpoint equalities and common-set argument are unchanged.
\end{proof}
This is an expository strengthening of a fixed-$q$ formulation by its block lemma; no independent priority is claimed. It removes the mere intersection $\bigcap_{q>2}\ell^q$ from the list of unresolved requirements in this setting.

For the binary system, the explicit block in Lemma~\ref{lem:univblock} also explains a near-$\ell^2$ gauge condition without suggesting a false square-summable endpoint. If $P=\sum_a\alpha_aw_a$, its block uses coefficients $-\alpha_a/K$, each repeated $K(2^s-1)$ times, where $2^s\le4K/\theta$. For every nonnegative nondecreasing gauge $\omega$ with $\omega(t)\to0$ as $t\downarrow0$,
\begin{equation}\label{eq:gauge-localization}
 \sum_n|c_n|^2\omega(|c_n|)
 \le\frac4\theta\sum_a|\alpha_a|^2
       \omega\!\left(\frac{\max_a|\alpha_a|}{K}\right)
 \longrightarrow0.
\end{equation}
One may therefore impose a summable gauge budget at each binary residual step, together with the preceding $q_j$ budgets. This gives the same ordinary-sum common-set conclusion with $\sum_n|c_n|^2\omega(|c_n|)<\infty$. Gauge-controlled representation results already occur in \cite{Epi06,EpiUniversal16}; a gauge condition alone is not a novelty claim. A logarithmic restriction coming from a probabilistic estimate is not automatically an intrinsic Walsh boundary.

The actual endpoint obstruction is unchanged. If $c\in\ell^2$, Riesz--Fischer gives one $F\in L^2$ to which every ordinary partial-sum subsequence converges in $L^2$. Its restriction to a positive-measure complement cannot also converge to $F+1$. The source's Theorem 3.8 proves this even when the exceptional set is allowed to depend on the target. Separately, if a mean-zero multiplier equals one off a set of measure $\theta<1$, Cauchy--Schwarz gives $\|\Psi\|_2^2\ge(1-\theta)/\theta$. This is a local cost, whereas uniqueness of the $L^2$ limit is a global obstruction.

\begin{proposition}[A Fourier-origin obstruction for unmodified sums]\label{prop:ordinary-Fourier-obstruction}
For $F\in L^1(\D)$, the full sequence $S_NF$ converges to $F$ in measure. On every fixed measurable set $B$ of positive measure, any subsequential limit in measure of the restrictions $S_NF|_B$ is $F|_B$. In particular, ordinary-sum universality for arbitrary targets on such a set cannot have unweighted $L^1$ Fourier origin.
\end{proposition}
\begin{proof}
Given $t>0$, approximate $F$ by a Walsh polynomial $P$ in $L^1$. For $N$ above its spectrum, $S_NP=P$. Proposition~\ref{prop:weakindividual} and Markov's inequality give
\[
 \mu\{|S_NF-F|>t\}
 \le \frac{2(C+1)}t\|F-P\|_1,
\]
uniformly for those $N$. Polynomial density proves convergence in measure. Restriction preserves it and limits in measure are unique up to null sets. Norm convergence in any finite positive power implies convergence in measure on $B$, so it cannot restore universality there.
\end{proof}
This obstruction does not concern pointwise approximation on a fixed null set or a sign-modified Fourier series. It explains why the formal series in Theorem~\ref{th:ordinary-vilenkin} cannot simply be declared to be the Fourier series of an unweighted integrable Walsh function. A missing Fourier-origin argument in that ordinary-sum, positive-measure formulation cannot be promoted to an open existence problem.

\subsection{Classical representation and abstract universality}\label{sub:classical-universality}
The representation problem predates the Walsh constructions discussed here. Talalyan's survey \cite{Talalyan60} and Ul'yanov's problem survey \cite{Ulyanov64} place Menshov-type representation, rearrangement, summability, and uniqueness in a common historical setting. They are useful precisely because these operations solve different problems. Pogosyan \cite{Pogosyan75} proves a rearrangement-universal formal series for each uniformly bounded complete orthonormal system: one fixed series can be rearranged to represent each measurable target almost everywhere. That theorem does not impose Fourier origin from one fixed $L^1$ source, and does not identify ordinary unmodified partial sums with rearranged sums. It therefore does not answer Question~\ref{ques:permutations}.

There is a second, topological line. Grosse-Erdmann \cite{GrosseErdmann99} treats universal families and hypercyclicity; Kahane \cite{Kahane00} develops category arguments for trigonometric series; Bayart--Grosse-Erdmann--Nestoridis--Papadimitropoulos \cite{BayartUniversal08} provide an abstract theory of universal series. In the present setting, a category proof needs a specified complete coefficient or function space and continuity of each finite partial-sum map. A constructive block theorem instead supplies quantitative support, agreement, and coefficient budgets. Neither proof strategy supplies the other's additional conclusions for free.

For generalized Walsh Fourier sums, Episkoposian--M\"uller \cite{EpiMueller16} study pointwise universality and topological transitivity, complementing the earlier compact-set work \cite{EpiMuller14}. This identifies a Fourier-origin line of investigation distinct from weighted formal representation and from the ordinary-sum common-set theorem above. The referenced pointwise setting must be specified before comparing a result with convergence in a weighted norm. These historical and abstract contributions broaden the comparison without replacing its local Walsh arguments.

\subsection{Universalizing weights: elementary order structure and the remaining question}\label{sub:weight-order}
Consider weights universalizing one fixed formal Walsh series by subsequences of its ordinary sums. The distinction from sign-selected or rearranged series is essential. For a fixed $C=\sum_nc_nw_n$, write
\[
 \begin{aligned}
 \mathfrak U(C)=\{v:{}&0<v\le1\text{ a.e., the ordinary sums of }C\\
 &\text{are subsequentially universal in }L^1(v\,d\mu)\}.
 \end{aligned}
\]
Some of this structure follows from a general density argument and should not be posed as an open problem.

\begin{proposition}[Downward inheritance]\label{prop:weight-inheritance}
If $u\in\mathfrak U(C)$ and $0<v\le A u$ almost everywhere for a finite constant $A$, with $0<v\le1$, then $v\in\mathfrak U(C)$. Every $v\in\mathfrak U(C)$ has $\operatorname*{ess\,inf}v=0$.
\end{proposition}
\begin{proof}
The inclusion $L^1(u\,d\mu)\to L^1(v\,d\mu)$ is continuous. Its range is dense because bounded functions lie in both spaces and truncations approximate any $L^1(v)$ target. Every tail of the ordinary-sum orbit is dense in $L^1(u)$ by subsequential universality, hence also dense in $L^1(v)$. Approximating a fixed target with errors $2^{-j}$ using increasing indices proves the first assertion. If $v\ge a>0$, the weighted norm is equivalent to the unweighted $L^1$ norm. For every fixed frequency $k$, all sufficiently late ordinary sums have Fourier coefficient $c_k$, so continuity of this coefficient functional excludes arbitrary limits. This proves the second assertion.
\end{proof}
Downward inheritance does not say that an arbitrary weight with essential infimum zero is universalizing. It also does not give closure under a pointwise maximum: convergence in two weaker weighted norms along unrelated subsequences need not produce one subsequence in the stronger norm. A common approximation skeleton would be a sufficient mechanism for a compatible family: every prescribed target neighborhood would have to be reached at the same stopping index in the two norms. The meaningful structural task is to characterize additional compatibility or joining conditions for a fixed orbit, rather than rediscover the lower-order implication.

\paragraph{Local research boundary.}
The published bounded-digit theorem and the diagonal corollary settle ordinary-sum existence with all superquadratic coefficient exponents and a common exceptional set. Proposition~\ref{prop:ordinary-Fourier-obstruction} rules out an unweighted integrable Fourier origin in the binary positive-measure formulation. Remaining distinct questions concern removal of bounded digits, prescribed coefficient profiles with stricter modification alphabets, and compatibility of several universalizing weights. Section~\ref{sec:boundaries} distinguishes source-listed problems from proposals of this survey.

\subsection{Universal Fourier functions and monotone coefficients}
Gevorkyan--Navasardyan's theorem \cite{GevNava99} provides an important foundation: a prescribed positive monotone null sequence outside $\ell^2$ admits ternary Walsh representations of arbitrary finite almost everywhere measurable functions. Later work imposes Fourier origin, common sets, weighted topologies, or stricter selection rules \cite{Nava07,Gri20,GGS}.
\begin{theorem}[Grigoryan \cite{Gri24}, Theorem 3]\label{th:grigoryan24}
There exists $U\in L^1[0,1]$ with decreasing Walsh coefficients and Fourier series convergent both in $L^1$ and almost everywhere, such that every real measurable function finite almost everywhere is the almost everywhere sum of
\[
 \sum_{n\ge0}\delta_n\wh U(n)w_n,\qquad\delta_n\in\{-1,1\}.
\]
\end{theorem}
The original Fourier series converges to $U$; the signs depend on the desired target. There is no contradiction between those statements. Sargsyan's weighted rearrangement result \cite{Sargsyan25} likewise uses a deliberately signed series, not automatically the original positive-coefficient series. Universal pairs \cite{UniversalPair23}, conditional universality \cite{GriConditional26,Avetisyan26}, and weighted quasiuniversality \cite{GM26} place different restrictions on these choices. They should be compared by their definitions rather than arranged as if each were a stronger form of one fixed assertion.

\subsection{Coefficient envelopes and Fourier origin}\label{subsec:diagonal}
The published prescribed-coefficient theorem \cite{GevNava99} and the Fourier-origin sign theorem \cite{Gri24} impose different data. A useful elementary bridge concerns a fixed monotone coefficient envelope. The following calculation is included as a direct consequence of the Dirichlet kernel estimate, without a novelty claim.

\begin{proposition}[A monotone Fourier envelope]\label{prop:monotone-envelope}
Suppose $d_n\ge0$ is nonincreasing and $d_n\le C(n+1)^{-1/2}$. Then the natural-order Walsh series $\sum_{n\ge0}d_nw_n$ converges almost everywhere and in every $L^p$, $1\le p<2$, to one function $U$ with $\widehat U(n)=d_n$.
\end{proposition}
\begin{proof}
For $x\in I_k\setminus I_{k+1}$, the dyadic expansion of $D_N$ gives $\sup_N|D_N(x)|\le C_0 2^k$. Abel summation therefore bounds every tail starting at $M$ by $2C_0 2^k d_M$. Since $d_M\to0$, the series converges at each $x\ne0$. Split a partial sum at $M=2^k$. The initial part is bounded by $\sum_{n<M}d_n\le C_1C2^{k/2}$ and the remaining part by $C_2C2^{k/2}$. Thus all partial sums are dominated by $C_3C2^{k/2}$ on that cylinder shell. Its $p$th power is integrable for $p<2$, since $\sum_{k\ge0}2^{kp/2}2^{-k-1}<\infty$. Dominated convergence proves the stated norm convergence. In particular $L^1$ convergence permits testing each character and identifies all Fourier coefficients.
\end{proof}
The decreasing envelope $d_n=(n+1)^{-1/2}$ itself is not square summable. Consequently Fourier origin in every subquadratic $L^p$ does not force $\ell^2$ coefficients. Conversely, an $\ell^2$ coefficient sequence cannot be universally modified by a fixed finite alphabet to represent arbitrarily large constant targets almost everywhere: Parseval bounds every partial sum in $L^2$ and Fatou preserves this bound in the limit. The envelope proposition concerns the original unmodified series; a universal representation requires the separate published theorem and its precise modification rule.

For coefficient constraints it is useful merely to define, without asserting an existence classification,
\[
 M(X,\ell^2)=\{A:(A_nc_n)_n\in\ell^2\text{ for every }c\in X\}.
\]
If $X$ is a solid sequence space and $d=Ac$ must be universally representable by signs, selectors, or ternary factors, then $A\notin M(X,\ell^2)$ is necessary. Sufficiency with a prescribed profile, a fixed Fourier source, and a stricter modification alphabet is an additional question, formulated in Section~\ref{sec:boundaries}.

\subsection{Quantifiers and a general category mechanism}\label{sub:v2category}
\begin{center}\small
\begin{tabular}{p{.25\textwidth}p{.63\textwidth}}\toprule
Operation & Fixed data and target-dependent choices\\\midrule
Pointwise subsequences & A countable $E$ is fixed; a generic continuous $f$ is chosen; the target $h:E\to\R$ then chooses a subsequence.\\
Correction & A target $f$ is given; its modification and, in the basic theorem, its large agreement set may depend on $f$.\\
Common correction profile & The set and absolute coefficients precede every target; its signs determine the corrected function.\\
Weighted subseries & One coefficient series is fixed; a weight is chosen before all targets in its weighted space; the target chooses retained terms.\\
Signs or ternary factors & Absolute coefficient values are fixed; the permitted alphabet determines whether omission is allowed.\\
Rearrangements & Every term must appear exactly once; filling unused indices is part of the construction.\\\bottomrule
\end{tabular}\end{center}
These are different quantifier patterns even if the same polynomial blocks appear in their proofs.

\begin{lemma}[Category criterion]\label{v2:category}
Let $X$ be a complete metric space, $Y$ a separable metric space, and $T_n:X\to Y$ continuous. Suppose that for every nonempty open $U\subset X$, $V\subset Y$, and every $m$, some $n\ge m$ satisfies $T_n(U)\cap V\ne\varnothing$. Then a dense $G_\delta$ subset of $X$ consists of points whose $T_n$-orbit has every $y\in Y$ as a subsequential limit.
\end{lemma}
This is the Baire criterion applied to the open dense sets $\bigcup_{n\ge m}T_n^{-1}(V_j)$, with $(V_j)$ a countable base of $Y$. Requiring every $m$ makes the approximating indices strictly increase.
\begin{proposition}\label{v2:genericc0}
In $c_0$ with its supremum norm, a dense $G_\delta$ set of coefficient sequences has formal Walsh partial sums subsequentially universal for every finite measurable target in the topology of convergence in measure.
\end{proposition}
\begin{proof}
Use $Y=L^0(\D)$ with metric $d(f,g)=\int\min(1,|f-g|)\dd$. Rational Walsh polynomials are dense by truncation and cylinder approximation. The map $T_Nc=\sum_{n<N}c_nw_n$ is continuous, as it uses finitely many coordinates. To check the lemma, start with a ball about $a\in c_0$ of radius $\delta$. Retain a sufficiently long finite prefix $R$ so the discarded coefficients are smaller than $\delta/3$. Approximate the desired residual $g-R$ in measure by a polynomial $P$. Lemma~\ref{lem:univblock} gives a later block $Q$ with all coefficients below $\delta/3$ and $Q=P$ outside a sufficiently small set. Set $c$ equal to the retained prefix, the new block, and zeros elsewhere. It lies in the original ball, and the endpoint partial sum approximates $g$ in measure. Place the block above any given lower bound on $N$. The category criterion applies.
\end{proof}
This assertion concerns formal coefficients, not Fourier coefficients of an unweighted integrable function. Although $\ell^2$ is dense in $c_0$, it has no such universal elements, because its orthogonal series has a fixed $L^2$ limit. A dense exceptional linear subspace is consistent with a residual universal class.

\subsection{Translations, coefficient obstructions, and operator transfer}\label{sub:v2Wiener}
\begin{proposition}[Compact-group Wiener criterion]\label{v2:Wiener}
For $f\in L^1(\D)$, the closed linear span in $L^1$ of its translates is all of $L^1$ exactly when $\wh f(n)\ne0$ for every $n$.
\end{proposition}
The detailed derivation is recorded in Appendix~\ref{app:auxiliary-proofs}, preserving the normalization used here.
This compact character argument explains part of the broader dyadic Tauberian tradition \cite{GolubovWiener03}; a noncompact Wiener theorem has additional content. It is also a useful contrast with universal subseries: linear combinations of translates permit arbitrary scalar weights, while selecting a fixed Fourier term allows only a fixed finite alphabet.

If a continuous linear automorphism $T$ of the target space commutes with every permitted approximation operator $U_n$, then density of $\{U_nf\}$ is equivalent to density of $\{U_nTf\}$, because the latter is its image under a homeomorphism. A bounded multiplier alone need not have that property: zero and finite-rank symbols already obstruct it. For restrictions to a set $E$, an operator on the restriction space and a commuting diagram must be supplied; boundedness on the original $L^p$ does not define it automatically.

Recent published developments include uniformly greedy-universal functions for generalized Walsh systems \cite{ESG26}, correction to weighted quasiuniversal functions \cite{GM26}, and the distinct constructions \cite{EGS22,EpiTigran25,EpiGreedy26}. The target class, permitted operations, and meaning of quasiuniversality are source-specific.

\subsection{Fourier origin, coefficient obstructions, and critical profiles}\label{sub:universal-comparison}
The continuity obstruction can be stated without reference to a particular orthogonal system. Suppose $\ell_k$ is a continuous coefficient functional on a target space $Y$, with $\ell_k(w_n)=\delta_{kn}$. If all partial sums of a selected subseries converge in $Y$, the coefficient of its limit is either $0$ or $c_k$. Hence it cannot realize arbitrary values of $\ell_k$ on $Y$. For a sign-modified series there are at most two possibilities, and for a ternary-modified series at most three. A rearrangement using all terms retains coefficient $c_k$ because eventually its $k$th term has occurred and no other term contributes to this functional. None of these modes can be universal in unweighted $L^1$ for arbitrary targets.

For $L^1(v\dd)$ the functional $f\mapsto\int fw_k\dd$, initially defined on cylinder polynomials, extends continuously exactly when $1/v$ is essentially bounded, since $|w_k|=1$ and the dual of $L^1(v\dd)$ is represented by essentially bounded multipliers relative to $v\dd$. If $v$ is positive but has arbitrarily small values on positive-measure sets, this particular obstruction disappears. It does not follow that every such weight produces universality; the common-weight construction establishes a carefully chosen one. The distinction between absence of an obstruction and existence of a representation is essential.

The $\ell^2$ obstruction uses a different topology. Let $d\in\ell^2$, and suppose a ternary modification converges almost everywhere to $f$. Its partial sums have $L^2$ norms at most $\|d\|_2$. Fatou gives $\|f\|_2\le\|d\|_2$. In particular, constant targets of arbitrarily large magnitude cannot all be realized. This applies without continuity of coefficient functionals in almost everywhere convergence. It explains a square-summability obstruction independently of coefficient-functionals. The monotone prescribed-coefficient theorem supplies a representation mechanism only within its stated selection rule.

\paragraph{Power and logarithmic envelope tests.}
The square-summability boundary can be checked without a universality theorem. For
$d_n=(n+2)^{-\beta}(\log(n+2))^{-\gamma}$ with $\beta\ge0$,
\[
 d\in\ell^2\quad\Longleftrightarrow\quad
 \beta>\tfrac12\quad\hbox{or}\quad
 \beta=\tfrac12\text{ and }\gamma>\tfrac12,
\]
within the positive eventually monotone null profiles under discussion. At the critical power the integral test is $\int^\infty t^{-1}(\log t)^{-2\gamma}\,dt$. Finitely many initial terms do not affect it. Profiles with $\beta=1/2$ and $0\le\gamma\le1/2$ satisfy the envelope proposition while escaping the square-summability obstruction; this identifies admissible candidates, not a criterion for every representation mode. A parameter test must distinguish these two logical conclusions.

\paragraph{Comparison and research boundary.}
The universality results now have a fixed comparison record: coefficients, allowed operations, topology, and order of choices. The common-set ordinary-sum theorem is formal, the sign theorem has Fourier origin, and the prescribed-coefficient theorem permits ternary factors. Questions~\ref{ques:permutations}--\ref{ques:joining-weights} address specific remaining changes of these data.

\subsection{Three constructive costs and their analytic obstructions}\label{sub:representation-costs}
The correction and representation arguments change different data, and their quantitative costs explain the difference. The local corrector in \eqref{eq:localcorrector} reproduces a height on most of a cylinder, but its unweighted error from that height has fixed $L^1$ mass. The compensating part is concentrated on the exceptional set. Consequently a large agreement set provides local reproduction, while the correction theorem still needs the residual induction and the bound on every intermediate sum.

The block in Lemma~\ref{lem:univblock} spreads each original coefficient over $K(2^s-1)$ frequencies with size reduced by $K$. Since $2^s\asymp K/\eps$, the $q$-mass has the scale $\eps^{-1}K^{2-q}$. For $q>2$ one can simultaneously pay finitely many coefficient budgets and make the unfinished-group error $O((\log K)/K)$ small. At $q=2$ this calculation stops improving. The global obstruction is stronger than the failure of this calculation: square-summable coefficients have a fixed $L^2$ limit for ordinary sums, and bounded $L^2$ partial sums for finite-alphabet modifications. These exclude the respective universal target classes by uniqueness or Fatou.

The common weight pays for a different part of the construction. On each bad set it suppresses the already chosen bounded block and all its partial sums. On the good set, the local inequality for every measurable subset permits layer-cake integration against that same weight. This is why a single weight can precede every target, while the subseries or permutation is chosen afterwards. The prescribed natural-order envelope is obtained earlier, by moving the blocks farther out. It has a separate cost in frequency location and gives no analogous bound for the target-dependent permutations.

\begin{longtable}{@{}p{.20\textwidth}p{.28\textwidth}p{.24\textwidth}p{.15\textwidth}@{}}
\caption{What the constructive theorems change}\label{tab:constructive-costs}\\
\toprule Construction & Fixed object and target choice & Cost controlled by the proof & Obstruction retained\\\midrule\endfirsthead
\toprule Construction & Fixed object and target choice & Cost controlled by the proof & Obstruction retained\\\midrule\endhead
Correction & The target is given before its modified Fourier source. & Exceptional measure, coefficient order, intermediate partial sums. & The original thresholding operators remain unbounded for $p\ne2$.\\
Weighted representation & A formal coefficient series is fixed; $\eps$ chooses the weight; the target chooses terms or their order. & Small coefficient masses, a common weighted error, completion of every rearrangement. & Continuous unweighted coefficient functionals exclude arbitrary $L^1$ targets.\\
Ordinary-sum universality & The coefficient series and the large set precede the target; only endpoints are selected. & Exact residual matching on common good sets. & An integrable Fourier source would already have its unique limit in measure.\\\bottomrule
\end{longtable}
For a proposed extension, the proof should specify which of these costs is changed. Prescribing one coefficient envelope in all subquadratic $L^p$ spaces is a source-construction requirement. Restricting a ternary alphabet to signs changes the permitted corrections. Joining two weighted topologies requires common return indices. None of these additional requirements is supplied by merely repeating an existential theorem with a newly chosen source.

\section{Uniqueness, measures, and multiparameter geometry}\label{sec:uniqueness}
Coefficient recovery is the reverse side of representation: which exceptional-set or limiting geometry still determines a series? In several variables, rectangles, cubes, cones, and triangles create genuinely different limiting procedures. Finite product representations explain the algebra but do not replace the required multiparameter estimates.

Recovering coefficients from a convergent formal series depends on the convergence geometry and on control of its dyadic densities. A finite measure, an integrable function, and an arbitrary formal series have different recovery mechanisms. In several parameters, cubes, cones, and coupled triangular averages lead to different maximal and uniqueness questions.

\subsection{Coefficient recovery and exceptional sets}
For a fixed mode of convergence and coefficient class, a $U$-set is an exceptional set outside which convergence of a formal Walsh series to zero forces all coefficients to vanish. An $M$-set admits a nontrivial such series. Measure zero alone does not decide the classification. The classical foundations include Arutyunyan--Talalyan \cite{ArutyunyanTalalyan64} and Crittenden--Shapiro \cite{CrittendenShapiro65}; modern work continues, rather than initiates, that theory.

Two immediate tests should precede subtler uniqueness theorems. If all partial sums tend to zero in $L^1$, testing against each $w_n$ shows $c_n=0$. If $c\in\ell^2$ and all partial sums tend to zero almost everywhere, their $L^2$ limit has an almost everywhere convergent subsequence, hence is zero; Parseval again gives $c=0$. These arguments do not apply to arbitrary almost everywhere convergent formal series.

The measure formulation makes the missing condition explicit.
\begin{theorem}\label{th:measurerecovery}
For $F_m=\sum_{n<2^m}c_nw_n$, there is a finite complex Borel measure $\nu$ with Walsh coefficients $c_n$ exactly when $\sup_m\norm{F_m}_1<\infty$. It is unique. It is absolutely continuous with respect to $\mu$ exactly when $(F_m)$ is uniformly integrable.
\end{theorem}
\begin{proof}
If $\nu$ exists, $F_m$ is its cylinder density: on a rank-$m$ cylinder $I$, its value is $\nu(I)/\mu(I)$. Thus $\norm{F_m}_1\le\norm\nu_{\mathrm{TV}}$. Conversely, define $L(P)=\int PF_m\dd$ for $m$ above the degree of $P$. Martingale consistency makes this independent of $m$, and $|L(P)|\le C\norm P_\infty$. Extend by density to $C(\D)$ and apply the Riesz representation theorem. Testing characters proves the coefficients and uniqueness. If $\nu=f\mu$, then $F_m=E_mf\to f$ in $L^1$, so the family is uniformly integrable. Conversely, uniform integrability and martingale convergence give an $L^1$ limit with the same coefficients, identifying $\nu$ with its density.
\end{proof}
For $c_n=w_n(a)$, $F_m=D_{2^m}(\cdot\oplus a)$ represents the point mass at $a$. It tends to zero off $a$ but its $L^1$ norm is one. Boundedness is not uniform integrability.

\begin{example}[A singular Rajchman measure]\label{ex:riesz}
There exists a singular probability measure $\nu$ with $\wh\nu(n)\to0$ whose dyadic cylinder densities tend to zero $\mu$-almost everywhere.
\end{example}
\begin{proof}
Take independent coordinate signs of $\nu$-means $a_j=(j+2)^{-1/2}$. Then
\[
 F_m=\prod_{j<m}(1+a_jr_j),\qquad\wh\nu(n)=\prod_{j:n_j=1}a_j.
\]
The largest active digit tends to infinity with $n$, so these coefficients tend to zero. Set $A_M=\sum_{j<M}a_j^2\to\infty$ and $X_M=\sum_{j<M}a_jr_j$. Under $\mu$ its mean is zero and variance $A_M$; under $\nu$ its mean is $A_M$ and variance at most $A_M$. Choose $M_k$ with $A_{M_k}\ge k^2$. Chebyshev, Borel--Cantelli, and a countable sequence of positive tolerances give $X_{M_k}/A_{M_k}\to0$ $\mu$-almost everywhere and $\to1$ $\nu$-almost everywhere. These disjoint full-measure events prove singularity. Since $\sum a_j^3<\infty$,
\[
 \log F_M=X_M-\tfrac12A_M+O(1).
\]
Along $M_k$ the density tends to zero $\mu$-almost everywhere. The nonnegative martingale $F_m$ has a finite almost everywhere limit, which must therefore be zero.
\end{proof}
This is a classical Riesz-product construction. It gives a dyadically convergent null series, not automatically convergence of all ordinary partial sums. In one dimension, convergence of all partial sums even at one point forces $c_n\to0$ by subtracting consecutive sums; dyadic convergence alone has no such implication.

Modern uniqueness and coefficient recovery include Kozlovskaya \cite{Kozlovskaya20}, Plotnikov \cite{Plotnikov07,Plotnikov10,Plotnikov17}, Lukomskii \cite{Lukomskii21}. Generalized integration gives another recovery mechanism: Skvortsov--Tulone \cite{SkvortsovTulone15} use multidimensional dyadic Kurzweil--Henstock and Perron-type integrals. Their integral and exceptional-set hypotheses cannot be replaced by the finite-measure criterion above.

\subsection{Cylinder primitives, Haar coefficients, and generalized recovery}\label{sub:primitive-recovery}
A formal Walsh series determines more information than its almost everywhere dyadic limit. This can be recorded before any integrability hypothesis. For $F_m=\sum_{n<2^m}c_nw_n$ and a rank-$m$ cylinder $I$, define
\begin{equation}\label{eq:cylinder-primitive}
 \nu_c(I)=\int_I F_m\,d\mu=2^{-m}F_m|_I.
\end{equation}
For $l\ge m$, integration over $I$ annihilates the higher-frequency terms, so the same value is obtained with $F_l$. In particular, if $I^0,I^1$ are the two children of $I$, then
\[
 \nu_c(I)=\nu_c(I^0)+\nu_c(I^1).
\]
Finite additivity on the cylinder algebra is therefore automatic. It does not include a bound on variation or countable additivity.

The Haar coefficients of this primitive can also be read algebraically. Choose $h_I=2^{m/2}(\ind_{I^0}-\ind_{I^1})$. Testing any sufficiently late $F_l$ gives
\begin{equation}\label{eq:primitive-haar}
 b_I=2^{m/2}\big(\nu_c(I^0)-\nu_c(I^1)\big).
\end{equation}
Conversely, the total mass and the differences $b_I$ reconstruct the child masses recursively:
\[
 \nu_c(I^0)=\tfrac12\big(\nu_c(I)+2^{-m/2}b_I\big),\quad
 \nu_c(I^1)=\tfrac12\big(\nu_c(I)-2^{-m/2}b_I\big).
\]
For a fixed $n<2^M$, the character is constant on each rank-$M$ cylinder, whence
\[
 c_n=\sum_{I\text{ of rank }M}w_n|_I\,\nu_c(I).
\]
These finite identities prove that the formal Walsh coefficients, the consistent cylinder primitive, and its formal Haar data encode the same object. They are a coefficient calculation, not a new recovery theorem from a terminal function.

The analytic question is which condition allows that primitive to be recovered from a pointwise sum. Theorem~\ref{th:measurerecovery} answers the finite-measure question through bounded variation, and uniform integrability identifies an $L^1$ density. Without those conditions, a generalized integral requires an appropriate derivation basis and a uniqueness principle for its primitives. The literature of generalized recovery, including \cite{Plotnikov10,SkvortsovTulone15}, addresses precisely this extra step.

\begin{theorem}[Plotnikov \cite{Plotnikov10}, Theorem 17]\label{th:countable-recovery}
Let a real formal Walsh series on $\D^d$ have cubical partial sums
\[
 S_N^\square(x)=\sum_{0\le n_1,\ldots,n_d<N}c_{\mathbf n}w_{\mathbf n}(x),\qquad N=1,2,\ldots.
\]
Suppose these sums converge to a finite function $f$ at every point outside an at most countable set. If $f\in L^1(\D^d)$, then its formal coefficients are the ordinary Fourier--Walsh coefficients of $f$:
\[
 c_{\mathbf n}=\int_{\D^d}f(x)w_{\mathbf n}(x)\,d\mu(x).
\]
The analogous complex assertion follows by applying the real result to the real and imaginary parts.
\end{theorem}
The source's Theorems 7 and 13 give the intervening recovery by its Perron-type $(PW)$ integral and agreement with Lebesgue integration. Skvortsov--Tulone \cite{SkvortsovTulone15} compare multidimensional dyadic Kurzweil--Henstock and Perron constructions for rectangular convergence. Their derivation bases and continuity requirements remain part of the external results.

A point mass illustrates why the mode cannot be weakened silently. For $c_n=w_n(a)$, the primitive in \eqref{eq:cylinder-primitive} is the restriction of $\delta_a$, while $F_m(x)\to0$ for every $x\ne a$. Nevertheless each term $c_nw_n(x)$ has modulus one, so the full ordinary series fails the term test at every point. The countable exceptional set is small enough for Theorem~\ref{th:countable-recovery}, but dyadic convergence is not its hypothesis. This separates failure of recovery from a dyadic limit and recovery under all integer cubes without contradicting either statement.

In one dimension, bounded variation of the cylinder primitive and absolute continuity are therefore different questions from uniqueness outside a specified set. In several dimensions a further distinction appears between rectangles and cubes: the convergence basis changes the derivative used to reconstruct the primitive. A complete recovery statement must fix the coefficient class, the limiting geometry, the exceptional set, and the integral. This is the common structure behind the measure criterion, generalized integration, and the $U$-set results discussed next.

\subsection{Cubes, cones, and the exceptional-set distinction}
On $\D^d$, write $w_n(x)=\prod_{j=1}^dw_{n_j}(x_j)$. Cubical convergence uses $S_{N,\ldots,N}$ for every integer $N\to\infty$. For fixed $\lambda>1$, restricted rectangular convergence uses $\min N_j\to\infty$ with $\max_{j,k}N_j/N_k\le\lambda$. Dyadic cubes, all integer cubes, and unrestricted rectangles are three different modes.

Plotnikov's countable-set uniqueness result \cite{Plotnikov10} states that every at most countable subset of the full compact product group is a $U$-set for cubical convergence. The group formulation retains all binary representatives and is relevant when translating interval statements. The statement is tied to all integer cubes; restricting to dyadic cubes changes the recovery problem. The distinction is already visible in the singular measure example above, whose cylinder densities vanish almost everywhere although the measure and some coefficients are nonzero. A null set, or a set cut out by character equations, does not by itself imply a uniqueness theorem. Any extension to a positive-dimensional exceptional set must specify the ordinary partial sums, the coefficient class, and the behavior outside that set. None of these hypotheses can be supplied merely by its Haar measure.

\subsection{Product maxima and triangular summation}
The rectangular martingale maximal operator in dimension $d$ is
\[
 M_{\mathrm{rect}}f=\sup_{m_1,\ldots,m_d}|E_{m_1}^{(1)}\cdots E_{m_d}^{(d)}f|.
\]
Positivity and iteration of Doob's inequality give $\|M_{\mathrm{rect}}f\|_p\le(p')^d\|f\|_p$ for $p>1$. For $d\ge2$ it is not of weak type $(1,1)$; the endpoint involves $L(\log L)^{d-1}$ instead.

One model demonstrates the logarithm. In two dimensions take $f_M=2^{2M}\ind_{I_M\times I_M}$, of $L^1$ norm one. At a point whose first nonzero coordinates are $k,l<M$, the average over $I_k\times I_l$ is $2^{k+l}$. For $k+l\ge M+1$, this exceeds $2^M$, while the corresponding shells have total measure at least $cM2^{-M}$. Thus $2^M\mu\{M_{\mathrm{rect}}f_M>2^M\}\gtrsim M$, ruling out weak $(1,1)$. Products of one-dimensional weak estimates cannot fix this obstruction. Multiparameter Hardy and summability theory \cite{Weisz02,MoriczSchippWade,MOProduct25} uses product structure explicitly.

Triangular sections are different again. Put
\[
 S_n^\triangle f=\sum_{i+j<n}\wh f(i,j)w_iw_j,
 \qquad\tau_n^\triangle f=\frac1n\sum_{l=0}^{n-1}S_l^\triangle f.
\]
A finite coefficient count gives
\[
 \tau_n^\triangle f=\frac1n\sum_{k=0}^{n-1}S_{k,n-k}f.
\]
Both sides assign coefficient $(i,j)$ the multiplier $(n-i-j-1)_+/n$.
\begin{theorem}[G\'at \cite{GatTriangle18}, Theorems 1.1 and 2.10]\label{th:triangle}
The maximal triangular Fej\'er operator has weak type $(1,1)$ and strong type $(p,p)$ for $1<p\le\infty$. In particular $\tau_n^\triangle f\to f$ almost everywhere for every $f\in L^1(\D^2)$.
\end{theorem}
The rectangular sizes in the displayed identity are coupled by $k+(n-k)=n$. Its kernel is not a product of one-dimensional Fej\'er kernels. The difficult step is the off-square localization estimate; density then gives convergence. Restricted two-dimensional N\"orlund convergence \cite{GogNagy25Restricted} introduces further geometric restrictions and should not be identified with either unrestricted products or triangles.

\subsection{Vilenkin groups and representative product systems}
For bases $b_j\ge2$, a Vilenkin group is $\prod_j\Z_{b_j}$ with place values $M_0=1$, $M_{j+1}=b_jM_j$. Its character blocks satisfy $D_{M_j}=M_j\ind_{I_j}$. If $\sup_jb_j<\infty$, martingale regularity is uniform and many atomic arguments have controlled constants. For unbounded bases, a child's average can exceed its parent's by an unbounded factor; a proof using a fixed branching bound must be reconsidered. G\'at's work \cite{GatUnbounded03,GatUnbounded07} identifies substantive intermediate classes, while interval estimates for unbounded systems \cite{TselishchevUnbounded24} require a different argument from merely replacing binary digits.

Representative product systems \cite{GatToledo24} enlarge the algebraic setting. On a nonabelian finite factor, irreducible representations have matrix entries and dimensions, and orthogonality carries dimension weights. Scalar multiplication by a Walsh character gives an exact XOR shift; multiplication of general representation coefficients decomposes through tensor products rather than one translated index. This is what fails in an unqualified scalar transfer. The Fourier and martingale structures still interact, but the binary identity and the normalization of kernels must be rebuilt in the representation-theoretic model.

\subsection{A quantitative product comparison and a recent divergence theorem}\label{sub:v2products}

There is a useful new connection between multidimensional Lebesgue functions and strong divergence. Getsadze \cite{Getsadze26}, Theorem 5, proves that for a uniformly bounded orthonormal system and finite spectral sets $A_n$, with $N_n=|A_n|\to\infty$, the condition
\[
 L_{A_n}(x)/\log N_n\longrightarrow\infty\quad\text{almost everywhere}
\]
produces an integrable $h$ and a subsequence for which $|S_{A_{n_k}}h(x)|\to\infty$ almost everywhere. The theorem is stated on a unit cube; its character-model specialization applies after the measure-preserving dyadic identification. For two-dimensional Walsh squares, put $A_q=[0,v_q)^2$ and $N_q=v_q^2$. Its Lebesgue function is constant in the spatial variables and equals $L_{v_q}^2$. The recursion gives $L_{v_q}\asymp q$, whereas $\log N_q\asymp q$. Consequently
\[
 \frac{L_{A_q}}{\log N_q}=\frac{L_{v_q}^2}{\log(v_q^2)}\asymp q\longrightarrow\infty.
\] This is a multidimensional pointwise mechanism, stronger than merely applying uniform boundedness to diverging operator norms. It does not contradict dyadic-cube martingale convergence.

Tephnadze's 2025 work \cite{TephStrong25} studies sharp strong convergence for two-dimensional Walsh partial sums in $H_p^\square$, $0<p<1$, without the earlier restriction that the two summation indices have comparable size. This is a further reason to state the geometry of the limiting indices explicitly. We record the direction of that published advance without reproducing its summability formula here.

\subsection{Finite representation theory behind representative products}\label{sub:v2representations}
Let $K$ be a finite group and $\pi$ a unitary irreducible representation of dimension $d_\pi$. Fix the convention
\[
 \wh f(\pi)=\frac1{|K|}\sum_{x\in K}f(x)\pi(x),\quad
 (f*g)(x)=\frac1{|K|}\sum_yf(y)g(y^{-1}x).
\]
Then
\begin{align*}
 f(x)&=\sum_{\pi\in\wh K}d_\pi\operatorname{Tr}(\wh f(\pi)\pi(x)^*),\\
 \frac1{|K|}\sum_x|f(x)|^2&=\sum_\pi d_\pi\norm{\wh f(\pi)}_{\mathrm{HS}}^2,\\
 \wh{f*g}(\pi)&=\wh f(\pi)\wh g(\pi).
\end{align*}
The matrix order in the last formula follows directly by writing $x=yz$ and using $\pi(yz)=\pi(y)\pi(z)$. It changes with the Fourier convention.

To verify orthogonality, average $\pi(x)A\rho(x)^*$ over $K$. The result intertwines $\rho$ and $\pi$. Schur's lemma makes it zero for inequivalent irreducibles and $\operatorname{Tr}(A)I/d_\pi$ when they coincide. Apply this to matrix units. Decomposing the regular representation supplies completeness and $\sum_\pi d_\pi^2=|K|$, giving inversion and Parseval. These are the exact dimension factors absent from scalar Walsh formulas.

On $S_3$, the standard two-dimensional character has values $2,0,-1$ on the three conjugacy classes. Its square is $1+\operatorname{sgn}+\chi_{\mathrm{std}}$, proving
$\mathrm{std}\otimes\mathrm{std}\simeq1\oplus\operatorname{sgn}\oplus\mathrm{std}$.
Products of matrix coefficients therefore split into several representations instead of translating one binary frequency. Infinite representative product systems \cite{GatToledo24} additionally require an ordering and uniform analytic control; the finite calculation explains the obstruction without claiming the infinite maximal theorems from algebra alone.

For completeness, hyperbolic scale sets have another geometry. The number of $d$-tuples $(k_1,\ldots,k_d)\in\N_0^d$ with $\sum_jk_j=M$ is
\[
 \binom{M+d-1}{d-1}\asymp_d(M+1)^{d-1}.
\]
The stars-and-bars count proves the identity. Thus summing dyadic blocks over a hyperbolic boundary can produce polynomial-in-$M$ factors even if each block has the same exponential scale. Such factors are the source of logarithmic terms when $M\asymp\log N$. This counting fact is not by itself a convergence theorem for arbitrary hyperbolic partial sums; kernel or square-function bounds still have to be supplied.

\paragraph{Comparison and research boundary.}
Product characters explain finite tensor identities, whereas uniqueness and divergence depend on the limiting geometry. Their application to numerical integration in the next section uses finite annihilators and a function norm, not a direct transfer of a multidimensional almost-everywhere convergence theorem.

\section{Digital nets and quantitative integration}\label{sec:qmc}
Digital integration asks for a bound on a linear functional over a prescribed function class. Character cancellation identifies the error frequencies exactly; smoothness, weights, and the chosen norm decide their total contribution. We retain the distinction between a finite identity, a worst-case error estimate, and a constructive rate as dimension grows.

Digital cubature turns a finite character average into exact cancellation outside an annihilator. Integration rates then follow only after the function space gives suitable coefficient or kernel decay. We develop the algebraic identity, a smooth Sobolev theorem with its parameter dependence, and a finite merit criterion in their respective normalizations.

\subsection{Annihilators and the exact error}
Digital integration gives a direct passage from character algebra to approximation rates. Let $P$ be a finite subgroup of $\D^s$ and let
\[
 P^\perp=\{k\in\N_0^s:w_k(x)=1\text{ for all }x\in P\}.
\]
Character orthogonality on $P$ gives
\begin{equation}\label{eq:qmcerror}
 \frac1{|P|}\sum_{x\in P}f(x)-\int f
 =\sum_{k\in P^\perp\setminus\{0\}}\wh f(k),
\end{equation}
initially for polynomials and then whenever the series is absolutely convergent. For a digital shift $\sigma$, the summands acquire factors $w_k(\sigma)$. Hence
\[
 \mathbb E_\sigma Q_{P\oplus\sigma}f=\int f,\qquad
 \Var_\sigma Q_{P\oplus\sigma}f
 =\sum_{k\in P^\perp\setminus\{0\}}|\wh f(k)|^2.
\]
For $L^2$ functions the identities can be obtained by polynomial approximation in the shift variable, avoiding evaluation of a fixed equivalence class at an exceptional deterministic point. Statements about deterministic cubature on a function space instead need well-defined evaluation functionals.

For the Walsh-diagonal RKHS with norm $\sum_k|\widehat f(k)|^2/r(k)$, the squared worst-case error is $\sum_{P^\perp\setminus\{0\}}r(k)$; the exact proof is in Proposition~\ref{v2:RKHSequality}. Smooth Sobolev kernels need not be Walsh diagonal. Averaging digital shifts removes their off-diagonal coefficients, connecting coefficient decay to smooth real functions \cite{Dick08,Dick09,DickPillichshammer,DKS13,GodaSuzuki20}.

\subsection{A precise high-order Sobolev rate}
Fix an integer $\alpha\ge2$. The space $H_\alpha$ consists of functions with absolutely continuous derivatives through order $\alpha-1$ and $f^{(\alpha)}\in L^2(0,1)$, with continuous representatives. Its unanchored inner product is
\[
 \langle f,g\rangle_{H_\alpha}=
 \sum_{r=0}^{\alpha-1}\Big(\int_0^1f^{(r)}\Big)
 \overline{\Big(\int_0^1g^{(r)}\Big)}
 +\int_0^1f^{(\alpha)}\overline{g^{(\alpha)}}.
\]
Its kernel, using Bernoulli polynomials, is
\[
 K_\alpha(x,y)=\sum_{r=0}^\alpha\frac{B_r(x)B_r(y)}{(r!)^2}
 +\frac{(-1)^{\alpha+1}B_{2\alpha}(|x-y|)}{(2\alpha)!}.
\]
For nonnegative coordinate-subset weights $\gamma_u$, define the mixed-smoothness space $H_{\alpha,\gamma}$ by the kernel
\[
 K_{\alpha,\gamma}(x,y)=\sum_{u\subseteq\{1,\ldots,s\}}
 \gamma_u\prod_{j\in u}(K_\alpha(x_j,y_j)-1).
\]
These are weights on groups of variables, not Muckenhoupt weights.

For an integer $\beta\ge1$ and an integer $k$, let $\mu_\beta(k)$ be the sum of the positions of its largest $\min(\beta,v)$ nonzero binary digits, with the lowest position numbered one. Extend by addition to vectors. An order-$\beta$ digital $(t,m,s)$ net satisfies the dual condition $\mu_\beta(k)>\beta m-t$ for every nonzero $k\in P^\perp$.
\begin{theorem}[Goda--Suzuki--Yoshiki \cite{GSY17}, Theorem 1 and Corollary 1]\label{th:qmcrate}
Let $\alpha,\beta,m$ be integers with $\alpha\ge2$, $m\ge1$, $\beta\ge2\alpha$, and let $0\le t\le\beta m$. Put $N=2^m$ and let $P$ be an order-$\beta$ digital $(t,m,s)$ net over $\F_2$, with $N$ nodes counted with the digital-net convention. Then
\[
 \Big(\mathbb E_\sigma e_{\mathrm{wor}}(H_{\alpha,\gamma};P\oplus\sigma)^2\Big)^{1/2}
 \le N^{-\alpha}\sum_{\varnothing\ne u\subseteq\{1,\ldots,s\}}
 \sqrt{\gamma_u}\,C_{\alpha,\beta,2,t,u}(\log N)^{(|u|-1)/2}.
\]
The constants are independent of $m$ for fixed $t$. In particular, at least one shift has this worst-case error bound.
\end{theorem}
For fixed dimension and positive product weights, the order $N^{-\alpha}(\log N)^{(s-1)/2}$ is optimal in the corresponding worst-case setting. The dependence on dimension and weights is not absent; dimension-independent tractability needs further summability assumptions on the weights.

The following is a proof mechanism, not a reproduction of the digital counting estimates in \cite{GSY17}. The first reduction is
\[
 \mathbb E_\sigma e_{\mathrm{wor}}(H;P\oplus\sigma)^2
 =\sum_{k\in P^\perp\setminus\{0\}}\wh K(k,k),
\]
under the kernel convergence conditions verified for this space. Character orthogonality cancels every off-diagonal term. The difficult part estimates the remaining diagonal sum using smoothness, digital weights, and the count of admissible dual frequencies. The explicit optimal construction \cite{GSYExplicit17} is stronger algorithmic information than existence of a good shift. The median constructions \cite{GodaSuzukiMatsumoto24} provide probabilistic robustness in a different model; worst-case, mean-square worst-case, and high-probability errors are not interchangeable.

\subsection{Walsh figure of merit and finite discretization}
Let $V=\F_2^{s\times n}$ with pairing $A\cdot B=\sum_{i,j}a_{ij}b_{ij}$ modulo two, and let $P\subset V$ be linear. Put $\mu(A)=\sum_{i=1}^s\sum_{j=1}^nj a_{ij}$. A finite Walsh figure of merit is
\begin{align}\label{eq:wafom}
 \operatorname{WAFOM}(P)&=\sum_{A\in P^\perp\setminus\{0\}}2^{-\mu(A)}\\
 &=\frac1{|P|}\sum_{B\in P}\prod_{i=1}^s\prod_{j=1}^n
       (1+(-1)^{b_{ij}}2^{-j})-1.\nonumber
\end{align}
Expand the product and average each character over $P$ to prove the second formula. This changes the cost from enumerating a large annihilator to evaluating a product on the net. For the finite normalized Fourier transform,
\[
 |Q_Pf-\mathbb E_Vf|
 \le\Big(\sup_{A\ne0}2^{\mu(A)}|\wh f(A)|\Big)\operatorname{WAFOM}(P).
\]
The coefficient assumption is essential. WAFOM alone is not an error bound for an unrestricted integrand \cite{MatsumotoOhori15}. Passing from an $n$-digit discretization to the continuum also incurs a discretization error; for a function Lipschitz in the $\ell^\infty$ metric, replacing each point by a representative of its digit cell adds at most $L2^{-n}$. These explicit hypotheses make the finite criterion useful without conflating it with every higher-order net theorem.

\subsection{Generating matrices, reproducing kernels, and equality cases}\label{sub:v2digital}
In a finite digital model, let $C_j\in\F_2^{n\times m}$ and define
$P=\{(C_1z,\ldots,C_sz):z\in\F_2^m\}$, assuming the combined map is injective when $|P|=2^m$ is required. Character averaging gives the exact dual equations
\[
 P^\perp=\{(u_1,\ldots,u_s):\sum_{j=1}^sC_j^Tu_j=0\}.
\]
The equations remain valid without injectivity for the image subgroup; repeated parametrizations have equal multiplicities. This small rank issue matters when counting the number of cubature nodes.

For positive summable weights $r(k)$ on the compact dual, define
\[
 \norm f_{\mathcal H_r}^2=\sum_k\frac{|\wh f(k)|^2}{r(k)},\qquad
 \mathcal K_r(x,y)=\sum_kr(k)w_k(x)\overline{w_k(y)}.
\]
Cauchy--Schwarz and $\sum r(k)<\infty$ imply absolute uniform convergence of every member's Fourier series, so point evaluations and cubature are defined.
\begin{proposition}\label{v2:RKHSequality}
For a finite subgroup $P$, the squared worst-case error on the unit ball is exactly
\[
 e_{\mathrm{wor}}(\mathcal H_r;P)^2
 =\sum_{k\in P^\perp\setminus\{0\}}r(k).
\]
\end{proposition}
\begin{proof}
Use \eqref{eq:qmcerror} and Cauchy--Schwarz with weights $r(k)$. If the right side is $A>0$, choose coefficients $r(k)/\sqrt A$ on the indicated dual set and zero elsewhere. This gives unit norm and error $\sqrt A$. If $A=0$, the functional vanishes. Uniform convergence justifies the evaluations and exchanges of sums.
\end{proof}
Coordinate weights, digit-position decay, and mixed smoothness are different ways to specify $r$ or the more general non-diagonal kernel. The error formula itself does not make constants independent of dimension.

For shifted nodes, polynomial character orthogonality proves unbiasedness and the variance formula in Section~\ref{sec:qmc}. For $f\in L^2$, each translated evaluation is an $L^2$ function of the random shift. Approximation by polynomials and the bound
$\norm{\mathcal Q_{P+\cdot}g}_{L^2(\D^s)}\le\norm g_2$
pass the identities to the limit. This establishes the random-shift formulas without pretending that deterministic point evaluation is continuous on $L^2$.

\subsection{Median amplification and the distinction between three guarantees}\label{sub:v2median}
\begin{lemma}\label{v2:median}
Suppose independent real-valued estimates $Q_1,\ldots,Q_L$, with $L$ odd, each satisfy
$\mathbb P\{|Q_j-I|>\eps\}\le q<1/2$. Their median satisfies
\[
 \mathbb P\{|\operatorname{med}(Q_j)-I|>\eps\}
 \le\exp(-2L(1/2-q)^2).
\]
\end{lemma}
A bad median requires at least half the independent estimates to be bad; Hoeffding's bound gives the display. The target and tolerance are fixed before this probability calculation.
A QMC application needs an integration-specific single-rule bound and independent randomizations, supplied in \cite{GodaSuzukiMatsumoto24}. The displayed probability fixes the target function before the randomization; it is not a probability bound for a supremum over the function ball.

There are two computational costs in the WAFOM expression \eqref{eq:wafom}. Enumerating $P^\perp$ can require $2^{sn}/|P|$ terms, while the product formula uses $O(sn|P|)$ elementary factors. Neither cost is independent of digit precision. The deterministic inequality weights coefficients by their full digit-position sum; the order-$\beta$ Sobolev theorem instead uses the first $\beta$ nonzero positions and a smoothness-specific kernel. The two weights are related in applications but are not interchangeable definitions of worst-case error.

\subsection{Exact tests and the recovery obstruction}\label{sub:qmc-bridge}
\paragraph{A one-dimensional exact test.}
Let $P_m\subset\D$ consist of the $2^m$ points whose coordinates vanish from position $m$ onward. Then
\[
 P_m^\perp=\{2^m\ell:\ell\ge0\}.
\]
The average of $w_n$ on $P_m$ is zero if one of its first $m$ binary frequency digits is nonzero, and one otherwise. In particular, $w_{2^m}$ has integral zero but cubature average one. The rule therefore cannot approximate the integral uniformly on the unit ball of unrestricted $L^\infty$. In the diagonal reproducing-kernel space with $r(k)=(k+1)^{-2\alpha}$, $\alpha>1/2$, the exact worst-case formula instead gives
\[
 e_{\rm wor}(H_r;P_m)^2
 =\sum_{\ell\ge1}(2^m\ell+1)^{-2\alpha}
 \asymp_\alpha2^{-2\alpha m}.
\]
The two-sided estimate follows by comparison with $\sum_{\ell\ge1}\ell^{-2\alpha}$. The rate is a consequence of this explicit Walsh-diagonal space and subgroup; it is not a claim that these points achieve arbitrary smooth Sobolev order. The latter is the separate content of Theorem~\ref{th:qmcrate}.

\paragraph{A finite uncertainty obstruction to joint design.}
On a finite group $V$, let $H$ be a subgroup. The indicator $1_H$ has Fourier support $H^\perp$. If a measurement design omits all rows in $H^\perp$, it annihilates this signal. The sparsity in the physical and spectral descriptions is linked by $|H||H^\perp|=|V|$. Thus the same annihilator geometry that gives exact integration cancellation can create a recovery nullspace for a particular signal class. A useful joint design must state which signals are excluded, which measurement rows are protected, and which restricted lower bound is sought. There is no general implication from a small digital integration error to a stable sparse-recovery matrix.

\paragraph{Comparison and research boundary.}
Digital cancellation and the RKHS identity establish an integration guarantee on a specified class. The recovery obstruction shows why the same nodes do not automatically identify a signal. Programme B in Section~\ref{sec:boundaries} asks for both guarantees with one finite design and explicitly retains their separate hypotheses.

\paragraph{Integration, recovery, and the meaning of a small spectral quantity.}
The digital error formula and sparse recovery examine different directions of the same measurement map. Integration controls one functional on a smoothness ball; recovery controls all differences of admissible signals. For the diagonal space $\mathcal H_r$, the first quantity is exactly the weighted mass of the unobserved dual characters. For a signal class $\Sigma$, a stable recovery requirement instead has the form
\[
 \inf_{z\in(\Sigma-\Sigma)\setminus\{0\}}
       \frac{\|Mz\|_2}{\|z\|_2}>0.
\]
A character in $P^\perp$ gives a direct integration error test, whereas a signal with Fourier support confined to omitted measurement rows gives a recovery null vector. The two tests share annihilator algebra but apply it to different admissible sets. The next section makes the signal side concrete through Walsh--Haar coherence, Boolean coefficient concentration, and query-based support discovery. This supplies an analytic connection between the application chapters without treating a cubature rate as a reconstruction theorem.

\subsection{An error dictionary for three Walsh integration models}\label{sub:qmc-dictionary}
The previous formulas describe three different ways to obtain an error estimate. In a Walsh-diagonal RKHS, the norm itself specifies the coefficient budget, and the worst-case error is an exact annihilator sum. In a smooth real-variable Sobolev space, the kernel is generally not diagonal in the Walsh basis; a digital shift removes cross terms in the mean-square calculation. WAFOM instead controls a finite digit model under a full digit-position coefficient decay condition. The statements become comparable only after the norm, randomization, and digit precision have been fixed.

\begin{longtable}{@{}p{.23\textwidth}p{.29\textwidth}p{.35\textwidth}@{}}
\caption{Integration error: function class, spectral budget, and guarantee}\label{tab:qmc-dictionary}\\
\toprule Function class & Spectral quantity & Established conclusion\\\midrule\endfirsthead
\toprule Function class & Spectral quantity & Established conclusion\\\midrule\endhead
Diagonal $\mathcal H_r$ & $\sum|\widehat f(k)|^2/r(k)$ & Exact worst-case error squared $\sum_{k\in P^\perp\setminus\{0\}}r(k)$; Proposition~\ref{v2:RKHSequality}.\\
One-dimensional $r(k)=(k+1)^{-2\alpha}$, $\alpha>1/2$ & Dual frequencies $2^m\ell$ & Error $\asymp_\alpha N^{-\alpha}$ for the explicit subgroup rule, $N=2^m$.\\
Mixed Sobolev $H_{\alpha,\gamma}$ & Diagonal Walsh coefficients of the shifted kernel & RMS worst-case order $N^{-\alpha}(\log N)^{(s-1)/2}$ under Theorem~\ref{th:qmcrate}.\\
Finite $n$-digit class & Full digit-position coefficient decay & Deterministic WAFOM bound plus a separate discretization error; Section~\ref{sec:qmc}.\\\bottomrule
\end{longtable}
The first two rows are exact model calculations. The third uses a substantial smoothness and digital-counting theorem. The last row does not imply the third for every Sobolev integrand: its coefficient condition must be verified at the chosen precision. Conversely, a rate for a random digital shift does not certify every deterministic shift. These comparisons give a concrete test for a proposed application: identify which row's assumptions hold before quoting its rate.

\section{Finite transforms, Boolean spectra, and mathematical applications}\label{sec:finite}
Finite transforms connect the survey with data and algorithms, but different applications ask for different guarantees. We begin with exact character algebra, then identify the extra structure needed for recovery, cryptanalytic screening, finite-state reversibility, reference scoring, or statistical calibration. A complete orthogonal change of coordinates is treated as such before any application-specific advantage is considered.

A finite Walsh transform is the character table of a binary vector space. This gives exact normalization, coding, covariance, and noise identities. Algorithmic and statistical guarantees enter after additional operations such as sampling or energy aggregation; those operations are treated explicitly rather than inferred from orthogonality.

\subsection{From an exact transform to an application guarantee}\label{sub:application-map}
The finite applications answer a common question: which additional structure turns character algebra into a conclusion about the scientific task? A fast orthogonal transform supplies a representation. The error criterion, signal class, probability model, or finite-state map supplies the remaining mathematics.
\begin{longtable}{@{}p{.17\textwidth}p{.23\textwidth}p{.23\textwidth}p{.23\textwidth}@{}}
\caption{Applications compared by the assumptions that create the guarantee}\label{tab:application-map}\\
\toprule Task & Exact Walsh input & Additional structure & Result and boundary\\\midrule\endfirsthead
\toprule Task & Exact Walsh input & Additional structure & Result and boundary\\\midrule\endhead
Digital integration & Character cancellation on a finite subgroup. & Smoothness or an RKHS norm and weighted coordinate assumptions. & Worst-case integration error; cancellation alone does not give recovery.\\
Sparse recovery & Character measurements and Walsh--Haar coherence. & Restricted injectivity and a sparsity model. & Stable identification only on the specified signal class.\\
S-box analysis & Walsh correlations and mask transformations. & A nonlinear Boolean map and a specified equivalence relation. & A spectral invariant or a correlation criterion, not a full-cipher security proof.\\
Hybrid encryption & An invertible transform or an explicitly defined finite-state layer. & Exact encoding and a protective cryptographic construction. & Reversibility and specified diffusion properties; secrecy requires its own model.\\
Quadratic signal models & Orthogonal conjugation and quadratic forms. & A positive covariance operator and a stated structural approximation. & Coordinate-invariant full scores; band constraints introduce a different model.\\
Covariance approximation & Walsh diagonalization of XOR convolution. & Translation invariance or an explicitly fitted structured approximation. & An exact diagonalization or an orthogonal projection onto the structured class.\\
Generalized bases & A unitary digit matrix with a constant first row. & Tensor products and the associated digit filtration. & Orthogonal coordinates and complete-block conditional expectations.\\\bottomrule
\end{longtable}
This map also identifies the appropriate baseline. A complete orthogonal rotation should be compared with the same invariant full-covariance score in the original coordinates; a structured covariance estimator should be compared with other estimators under the same data split. A reversible preprocessing scheme should be compared with its protective cipher alone, including serialization cost. A spectral invariant should be compared with established invariants on the same equivalence classes. The same mathematical baseline applies to every proposed implementation.

\subsection{Normalization, sampling, and randomized rotations}
On $V=\F_2^m$, $N=2^m$, write
\[
 H_N[u,x]=(-1)^{u\cdot x},\qquad Q_N=N^{-1/2}H_N.
\]
Then $H_N^TH_N=NI$, whereas $Q_N$ is orthogonal. We distinguish the unnormalized transform $H_Nf$, the average-normalized coefficients $N^{-1}H_Nf$, and the unitary coefficients $Q_Nf$. With average-normalized coefficients, inversion is $f(x)=\sum_u\wh f(u)(-1)^{u\cdot x}$ and Parseval is $N^{-1}\sum_x|f(x)|^2=\sum_u|\wh f(u)|^2$.

The recursion $H_{2N}=\left(\begin{smallmatrix}H_N&H_N\\H_N&-H_N\end{smallmatrix}\right)$ gives $N\log_2N$ additions and subtractions. Euclidean perturbations preserve their norm under $Q_N$. Keeping the largest coefficients minimizes the Euclidean approximation error, but Section~\ref{sec:greedy} explains why this fact is norm-specific.

Character algebra gives a finite uncertainty bound. For $f\not\equiv0$, if $A=\supp f$ and $B=\supp\wh f$, then
\[
 \norm f_\infty\le\sqrt{|B|}\norm{\wh f}_2,
 \qquad \norm f_{L^2(V)}\le\sqrt{|A|/N}\norm f_\infty,
\]
and Parseval implies $|A||B|\ge N$. Equality occurs for an indicator of a subgroup and its annihilator. Consequently a Euclidean fractal uncertainty principle cannot simply be transferred to Walsh analysis; the specific failure is studied by Demeter \cite{Demeter22}.

For a fixed $x\in\R^N$, $0<\eta<1$, and independent random signs $D=\operatorname{diag}(\eps_i)$,
\begin{equation}\label{eq:rotation}
 \mathbb P\left\{\norm{Q_NDx}_\infty>
 \norm x_2\sqrt{\frac{2\log(2N/\eta)}N}\right\}\le\eta.
\end{equation}
Each coordinate is a signed sum with moment-generating function at most $\exp(t^2\norm x_2^2/(2N))$; the subgaussian tail and a union bound prove the formula. For a general orthogonal $Q$, it remains true with a factor $\kappa$ if $\max_{i,j}|Q_{ij}|\le\kappa/\sqrt N$. Orthogonality alone is insufficient: $Q=I$, $x=e_1$ leaves maximum coordinate one with probability one.

Subsampled randomized Hadamard transforms \cite{Tropp11} combine this spreading with sampling. Sparse recovery requires assumptions about both the sparsifying basis and the measurement distribution. A constant vector has only its DC Walsh coefficient, so a sampling scheme omitting that row loses it entirely. Local coherence and structured sparsity resolve such mismatches in \cite{AdcockEtal}. Sparse transform algorithms \cite{CheraghchiIndyk16} target another problem: locating heavy Walsh coefficients with fewer operations or queries. An algorithm for a sparse transform is not automatically a recovery theorem for every measurement design.

Modern machine-learning implementations use Hadamard rotations to redistribute coordinate magnitudes, for example in quantization \cite{QuaRot,QuIP}. Equation \eqref{eq:rotation} explains the fixed-vector probabilistic mechanism. Accuracy of a quantized trained model needs additional error and data assumptions; it is not a corollary of orthogonality alone.

For a matrix $W$, orthogonal $U,V$ satisfy $\|UWV\|_F=\|W\|_F$; a unitary discrete Fourier transform has the same energy-preservation property. The implementation distinction is that a real Hadamard rotation has entries of modulus $N^{-1/2}$ and an addition--subtraction recursion. Random signs then provide the coordinate-spreading estimate \eqref{eq:rotation}. An arbitrary fixed Hadamard transform does not flatten every vector: $x=Q_N^Te_j$ is transformed to $e_j$. Thus quantization error, coordinate outliers, and end-task accuracy require their own assumptions and comparisons; no energy advantage over the unitary DFT is asserted.

\subsection{Boolean functions, coding, and differential identities}
For $f:V\to\F_2$, let $F=(-1)^f$ and $W_f(u)=\sum_x(-1)^{f(x)+u\cdot x}=N\wh F(u)$. The distance from $f$ to the affine function $u\cdot x+b$ is $(N-(-1)^bW_f(u))/2$. Therefore
\begin{equation}\label{eq:nonlinearity}
 \mathrm{NL}(f)=N/2-\tfrac12\max_u|W_f(u)|.
\end{equation}
Parseval gives $\sum_uW_f(u)^2=N^2$, so $\max|W_f|\ge\sqrt N$. Bent functions attain equality at every frequency, possible only when $m$ is even. Rothaus \cite{Rothaus} and Carlet \cite{Carlet} supply the foundational and modern coding-cryptography setting.

A concrete family is $f(x,y)=x\cdot\pi(y)+g(y)$ on $\F_2^r\times\F_2^r$, with $\pi$ a permutation. Summing first over $x$ gives
\[
 W_f(a,b)=2^r(-1)^{g(\pi^{-1}(a))+b\cdot\pi^{-1}(a)},
\]
so every coefficient has modulus $2^r$. The construction proves a spectral property, not security of an entire cryptographic protocol. Correlation immunity likewise has an exact spectral meaning: independence of the output from every set of at most $t$ input coordinates is equivalent to vanishing $W_f(u)$ for $1\le|u|\le t$, with balance handled separately. Reed--Muller codes encode bounded algebraic degree; their Hamming distance questions connect to \eqref{eq:nonlinearity}, but algebraic degree and Fourier Hamming support are distinct.

For a vectorial map $S:V\to\F_2^r$, put $R=2^r$ for the codomain cardinality and define
\[
 \mathrm{DDT}_S(a,b)=|\{x:S(x+a)+S(x)=b\}|,
\quad W_S(u,v)=\sum_x(-1)^{u\cdot x+v\cdot S(x)}.
\]
A useful fourth-moment identity is
\begin{equation}\label{eq:ddt}
 \sum_{a,b}\mathrm{DDT}_S(a,b)^2
 =\frac1{NR}\sum_{u,v}W_S(u,v)^4.
\end{equation}
To prove it, apply Parseval in $b$ to the DDT row; its transform is the autocorrelation $\sum_x(-1)^{v\cdot(S(x+a)+S(x))}$. Then apply Parseval in $a$, using that this autocorrelation has Fourier coefficients $W_S(u,v)^2/N$. The normalization in \eqref{eq:ddt} follows. Differential and linear statistics are connected by this aggregate identity; it does not identify their worst individual entries.

\paragraph{Classical input and modern use.}
The finite identities above are background tools for the 2000--2026 developments in sparse Walsh computation, structured measurements, and quantization. The older learning, coding, and quantum examples below explain the distinct computational models in which the same spectrum is queried. Their dates identify their historical role; they are included as prerequisites rather than as recent advances. The operative comparison is between a known coefficient support, a support that must be found, and access to individual function queries.

\subsection{Noise, influence, and low-degree approximation}
For $0\le\rho\le1$, define Boolean noise by independently flipping each bit with probability $(1-\rho)/2$. Then
\[
 T_\rho F=\sum_u\rho^{|u|}\wh F(u)w_u,\qquad
 \mathbb E[F(X)F(Y)]=\sum_u\rho^{|u|}|\wh F(u)|^2
\]
for a real function and its noisy input. This is Hamming geometry, unlike the numerical multiplier $e^{-tn}$. For Boolean-valued signs, the influence of coordinate $j$ is $\sum_{u:u_j=1}|\wh F(u)|^2$, obtained by Parseval on the discrete difference. Hypercontractivity controls low-degree polynomials; low-influence invariance results such as \cite{MOO10} link the Boolean model to Gaussian analysis. The systematic reference is \cite{ODonnell}.

If $\sum_{u\notin\Lambda}|\wh F(u)|^2\le\eps$, the polynomial $P_\Lambda F$ has mean-square error at most $\eps$. For signs $F\in\{-1,1\}$, classification by the sign of its real-valued approximation has error probability at most $\eps$, since a sign mistake forces pointwise squared error at least one. Estimating coefficients, locating the unknown heavy support, and controlling generalization are additional tasks. Goldreich--Levin \cite{GoldreichLevin} and sparse transform methods supply algorithmic connections; the approximation statement alone is not a complete learning guarantee.

\paragraph{Fourier learning and Boolean structure.}
The distinction between known support and finding support leads to two classical lines: Kushilevitz--Mansour \cite{KM93} give a membership-query learning algorithm under the uniform distribution, whereas Linial--Mansour--Nisan \cite{LMN93} exploit low-degree Fourier concentration of constant-depth circuits. These are model-specific learning guarantees, rather than consequences of Parseval alone. Two structural theorems explain different spectral restrictions. Kahn--Kalai--Linial \cite{KKL88} imply
\[
 \max_j\operatorname{Inf}_j(F)\ge c\operatorname{Var}(F)\frac{\log m}{m},
 \qquad F:\F_2^m\to\{-1,1\},\quad m\ge2,
\]
with an absolute $c>0$. Friedgut--Kalai--Naor \cite{FKN02} show that Boolean functions with almost all Fourier mass at degrees zero and one are close to a constant or a signed coordinate. These are sensitivity and stability statements, not unrestricted sparsity assertions; see also \cite{ODonnell} for conventions and proofs.

\subsection{Hadamard sampling, codes, and sequency}\label{sub:hadamard-sampling}
The same finite character table has a quantum interpretation. Write $\mathsf H=2^{-1/2}\left(\begin{smallmatrix}1&1\\1&-1\end{smallmatrix}\right)$, so $\mathsf H^{\otimes m}=Q_N$. After uniform preparation and a phase oracle for $f:V\to\F_2$, the state is $N^{-1/2}\sum_x(-1)^{f(x)}|x\rangle$. Applying $Q_N$ gives amplitude
\[
 \frac1N\sum_x(-1)^{f(x)+u\cdot x}=\widehat{(-1)^f}(u)
\]
at $|u\rangle$; hence measurement samples squared Walsh coefficients. For the linear promise $f(x)=a\cdot x$, the final state is exactly $|a\rangle$. This is the character calculation underlying the Bernstein--Vazirani oracle example \cite{BV97}. It is a statement about oracle queries and a specified quantum circuit, not a classical speedup from performing an FWHT on an explicitly supplied table.

For synchronous finite-channel coding, distinct Walsh rows satisfy $N^{-1}\langle H_N[u,\cdot],H_N[v,\cdot]\rangle=\delta_{uv}$. An aligned noiseless superposition is therefore separated by matched Walsh correlations. A time shift, multipath, or asynchronous sampling changes these cross-correlations and needs a channel model; code orthogonality at zero lag alone supplies no asynchronous guarantee. Harmuth's engineering treatment \cite{Harmuth72} is a historical source for orthogonal signaling and sequency. Sequency reorders a fixed finite family and leaves full-block orthogonality intact, while changing prefixes, bandwidth interpretation, and truncation. The Gray and bit-order dictionary in Section~\ref{sec:foundations} must therefore accompany an implementation.

\subsection{Exact Walsh--Haar coherence and missing-scale information}\label{sub:v2coherence}
For a rank-$j$ cylinder $I$ with digits $a_0,\ldots,a_{j-1}$, choose $h_I=2^{j/2}r_j\ind_I$.
\begin{proposition}\label{v2:coherence}
If $n\notin[2^j,2^{j+1})$, then $\langle w_n,h_I\rangle=0$. If $n=2^j+u$, $u<2^j$, then
\[
 \langle w_n,h_I\rangle=2^{-j/2}(-1)^{\sum_{k<j}u_ka_k}.
\]
\end{proposition}
\begin{proof}
A nonzero higher digit of $n$ averages to zero on the free coordinates of $I$. In their absence, a missing digit at $j$ leaves the mean-zero $r_j$ factor. In the remaining case $r_j$ cancels and all lower digits are constant on $I$; the integral is $2^{j/2}\mu(I)$ times the indicated sign.
\end{proof}
The Haar scale-$j$ subspace and Walsh block $[2^j,2^{j+1})$ coincide. Their change-of-basis matrix is a normalized $2^j$-dimensional Hadamard matrix. Global coherence is one because of the lowest scales, although local coherence decreases with $j$. If a measurement set misses this whole Walsh block, every Haar function of that scale has zero measurements. Even a one-coefficient signal is then indistinguishable from zero; no reconstruction algorithm can repair the information loss. Full measurement of that block, conversely, recovers its subspace exactly. Intermediate guarantees need a level-dependent sparsity and sampling model, as in \cite{AdcockEtal}.

\subsection{XOR stationarity and spectral covariance}\label{sub:v2stationarity}
\begin{theorem}\label{v2:covariance}
For a finite binary group $V$, $|V|=N$, the matrix $[C(x\oplus y)]_{x,y}$ is positive semidefinite exactly when all normalized Walsh coefficients $\wh C(u)$ are real and nonnegative. Its eigenvalues are $N\wh C(u)$. A centered field with this covariance has transformed coordinates $Z_u=N^{-1/2}\sum_xX_xw_u(x)$ satisfying
\[
 \mathbb E Z_u\overline{Z_v}=N\wh C(u)\ind_{\{u=v\}}.
\]
\end{theorem}
\begin{proof}
Changing variables $z=x\oplus y$ gives
$\sum_yC(x\oplus y)w_u(y)=N\wh C(u)w_u(x)$.
The characters form a complete orthogonal eigenbasis, proving the positivity criterion. One may realize the covariance by
$X_x=\sum_u\sqrt{\wh C(u)}\xi_uw_u(x)$ with independent centered unit-variance $\xi_u$. Transforming this formula, or conjugating the covariance matrix by $Q_N$, gives the last assertion.
\end{proof}
Uncorrelated transformed coordinates are independent under joint real Gaussianity, but not in general. Ordinary time stationarity uses addition with carries and need not have this XOR covariance. The expected periodogram $N^{-1}|Z_u|^2$ equals $\wh C(u)$; consistency of an estimator under growing sample size requires a separate asymptotic model.

\subsection{Coding identities, hypercontractivity, and the three scales}\label{sub:v2coding}
The distances in \eqref{eq:nonlinearity} give nearest-word decoding in $RM(1,m)$ by maximizing $|W_f(u)|$; the sign determines the affine constant. This uses $O(N\log N)$ operations in the truth-table model, while heavy-coefficient oracle search is measured in $m=\log_2N$.

The Boolean hypercontractive theorem \cite{ODonnell} states
\[
 \norm{T_\rho F}_q\le\norm F_p,
 \qquad1<p\le q<\infty,\quad0\le\rho\le\sqrt{\frac{p-1}{q-1}}.
\]
For a polynomial of Hamming degree at most $d$, take $p=2$, $\rho=(q-1)^{-1/2}$, and apply the theorem to $\sum\rho^{-|u|}\wh F(u)w_u$. Parseval gives
\[
 \norm F_q\le(q-1)^{d/2}\norm F_2,\qquad q\ge2.
\]
The deep input is the tensorized two-point inequality, not the last coefficient calculation. This is a concrete dimension-free consequence rather than a generic statement that noise ``improves integrability.''

There are one-sided comparisons between the three scales even though they are not equivalent. For a sign-valued Boolean function, Parseval gives
\[
 \sum_j\operatorname{Inf}_j(F)=\sum_nh(n)|\wh F(n)|^2
 \le\sum_nn|\wh F(n)|^2\le\sum_nn^2|\wh F(n)|^2.
\]
Thus a small numerical $L^2$ derivative cannot coexist with arbitrarily large total influence in this normalization. The converse comparison fails: $w_{2^j}$ has Hamming weight one while its numerical eigenvalue tends to infinity. This clarifies which comparisons are false and which inequalities remain valid.

\begin{proposition}[Autocorrelation and the bent criterion]\label{v2:autocorrelation}
For $F=(-1)^f$, put $R_f(a)=\sum_xF(x)F(x\oplus a)$. Then
\[
 \sum_aR_f(a)w_u(a)=W_f(u)^2,\qquad
 R_f(a)=N^{-1}\sum_uW_f(u)^2w_u(a).
\]
Consequently $f$ is bent exactly when $R_f(a)=0$ for every $a\ne0$.
\end{proposition}
\begin{proof}
In the double sum substitute $y=x\oplus a$ and separate the character factors. Fourier inversion gives the second formula. Constant squared spectral modulus $N$ is equivalent to autocorrelation $N\ind_{\{0\}}$, proving the criterion in both directions.
\end{proof}
Balance of the Boolean derivative $f(x\oplus a)+f(x)$ is exactly $R_f(a)=0$. The strict avalanche criterion restricts this to unit-weight $a$, and $PC(k)$ to $1\le|a|\le k$; these are not the same as vanishing first-order Walsh coefficients. For correlation immunity, condition on coordinates in $J$. The conditional expectation of $F$ is
\[
 \mathbb E(F\mid X_J)=\sum_{\supp u\subseteq J}\wh F(u)w_u(X_J).
\]
It is constant exactly when its nonzero frequencies vanish. Since $F$ is sign-valued, this is equivalent to independence of its output distribution from $X_J$. Taking every $|J|\le t$ proves the immunity criterion; the separate zero-frequency condition gives balance.

\subsection{What spectral aggregation preserves and loses}\label{sub:finite-guarantees}
Let $P_b$ be mutually orthogonal real projections with positive ranks $d_b$ and $\sum_bP_b=I$. The energies $e_b(x)=\|P_bx\|_2^2$ and normalized profile $p_b(x)=e_b(x)/\|x\|_2^2$ illustrate the cost of passing from coefficients to features. For nonzero $x,y$,
\[
 \|p(x)-p(y)\|_1
 \le 2\left\|\frac{x}{\|x\|_2}-\frac{y}{\|y\|_2}\right\|_2.
\]
Indeed, writing the two normalized vectors as $u,v$, factor each difference of squared block norms and apply Cauchy--Schwarz across the blocks. The square sum of the block differences is at most $\|u-v\|_2^2$, and the square sum of the block sums is at most four.

Stability does not imply recovery. Every $U=\bigoplus_bU_b$ with orthogonal $U_b$ preserves the profile. Two nonzero blocks already permit independent sign changes; a block of dimension at least two permits a continuous rotation orbit. A recovery statement must therefore target this quotient or impose an additional signal model. Dyadic translations preserve spectral-group energies because they multiply each Walsh coefficient by a sign. Ordinary cyclic shifts need not do so.

The statistical conclusion also depends on the noise model. If $Y=x+\sigma Z$ with $\sigma>0$, deterministic $x$, and $Z\sim N(0,I_N)$, then the block energies are independent and
\[
 e_b(Y)/\sigma^2\sim\chi^2_{d_b}(\|P_bx\|_2^2/\sigma^2).
\]
For $x=0$, independent gamma normalization gives
\[
 (p_b(Y))_b\sim\operatorname{Dir}((d_b/2)_b),
 \qquad \mathbb Ep_b=d_b/N.
\]
For $x\ne0$ the normalized law is noncentral and is not this Dirichlet law. Thus neither calibration nor discrimination follows from transform orthogonality alone. These elementary boundaries delimit the analytic content of applications without importing empirical claims into the survey.

\subsection{Cryptanalytic spectra: concentration and affine changes}\label{sub:applied-crypto}

There are two distinct uses of Walsh analysis in information protection. One computes the correlations of a nonlinear Boolean component; the other transforms a plaintext or signal before a separate encryption stage. The former has a direct interpretation in linear cryptanalysis. The latter requires a complete finite-state encoding and a security argument for the resulting cryptosystem. Carlet's monograph \cite{Carlet} supplies the Boolean-function framework; Lu--Desmedt \cite{LuDesmedt16} illustrate the use of Walsh expansions to compute biases when the input bits are constrained or dependent. Multiplication of individual biases without checking their joint model is not justified by the existence of a fast transform.

For a permutation $F:\F_2^m\to\F_2^m$, put $N=2^m$ and use the unnormalized table $W_F(a,b)$ from \eqref{eq:ddt}. The fourth-power sum and the maximum of the same table are
\[
 I_4(F)=\sum_{b\ne0}\sum_a|W_F(a,b)|^4,
 \qquad \Lambda(F)=\max_{a,\,b\ne0}|W_F(a,b)|.
\]
Every nonzero component of a permutation is balanced, so $W_F(0,b)=0$. Parseval gives $\sum_{b\ne0,a}|W_F(a,b)|^2=(N-1)N^2$. Consequently
\begin{equation}\label{eq:applied-fourth}
 N^4\le I_4(F)\le N^4(N-1),\qquad
 \frac{I_4(F)}{(N-1)N^2}\le\Lambda(F)^2\le I_4(F)^{1/2}.
\end{equation}
For the lower bound apply Cauchy--Schwarz to the at most $(N-1)^2$ nonzero-mask entries; for the upper bound use $|W_F|\le N$. The remaining inequalities follow by comparing the fourth sum with $\Lambda^2$ times the second sum and with its largest summand. These are concentration bounds, not a determination of the minimum fourth moment in each dimension. An affine permutation attains the upper bound. Under invertible affine changes of input and output, the Walsh magnitudes are permuted, so $I_4$ is invariant.

Identity~\eqref{eq:ddt} also gives the exact redundancy relation
\begin{equation}\label{eq:applied-moment-ddt}
 I_4(F)=N^2\sum_{a,b}\mathrm{DDT}_F(a,b)^2-N^4.
\end{equation}
Indeed, the omitted output mask $b=0$ contributes only $W_F(0,0)^4=N^4$. Thus $I_4$ can distinguish permutations having the same maximum Walsh magnitude, but it is already determined by the full DDT multiset. It complements nonlinearity; it cannot be presented as information independent of all differential statistics. These identities follow directly from the established Walsh/DDT relation; no additional cryptanalytic invariant is claimed.

Under invertible affine input and output changes, $G(x)=AF(Bx+c)+d$, direct substitution gives
\[
 W_G(a,b)=(-1)^{(B^{-T}a)\cdot c+b\cdot d}
 W_F(B^{-T}a,A^Tb).
\]
Thus moments of the full magnitude table are invariant under these changes, whereas Hamming weight of an individual mask generally is not. For a permutation, $W_{F^{-1}}(b,a)=W_F(a,b)$. These are changes of character coordinates; they do not establish an attack bound for a cipher containing the component. The published Boolean-function theory \cite{Carlet} distinguishes the linear and differential criteria that use these tables.

\subsection{Generalized Walsh bases and finite encoding}\label{sub:applied-encryption}
Dutkay--Picioroaga \cite{DutkayPicioroaga14} develop generalized Walsh bases from unitary matrices and discuss encoding applications. The pertinent analytic distinction is between complete digit blocks and arbitrary partial sums: a block of size $b^q$ reproduces a conditional expectation for the associated $b$-adic filtration. An unrestricted almost-everywhere theorem for all $L^1$ partial sums is not inferred from this block identity.

Here is the finite-dimensional mechanism in explicit normalization. Let $b\ge2$, let $A=(a_{ij})_{0\le i,j<b}$ be unitary, and suppose its first row is constant $b^{-1/2}$. On $[0,1)$ set
\[
 m_i(x)=\sqrt b\sum_{j=0}^{b-1}a_{ij}\ind_{[j/b,(j+1)/b)}(x),
 \qquad R_bx=bx\pmod1.
\]
If $n=\sum_{\ell=0}^{q-1}i_\ell b^\ell<b^q$, put
\[
 W^A_n(x)=\prod_{\ell=0}^{q-1}m_{i_\ell}(R_b^\ell x).
\]
Leading zero digits add no factor, since $m_0=1$. On a rank-$q$ digit interval with digits $j_0,\ldots,j_{q-1}$ the value is
$b^{q/2}\prod_\ell a_{i_\ell j_\ell}$. Summing the product inner products over all digit intervals gives
\[
 \langle W^A_n,W^A_k\rangle
 =\prod_{\ell=0}^{q-1}\sum_{j=0}^{b-1}a_{i_\ell j}\overline{a_{k_\ell j}}
 =\delta_{nk}.
\]
There are $b^q$ such functions and the rank-$q$ step space has dimension $b^q$. Thus they form its orthonormal basis. The union of these step spaces is dense in $L^2[0,1)$, proving completeness. In particular,
\[
 \sum_{n<b^q}\langle f,W^A_n\rangle W^A_n
 =\mathbb E(f\mid\mathcal F_q),
\]
where $\mathcal F_q$ is the $b$-adic digit sigma field. The equality first follows in $L^2$ and extends to $L^1$ because the finite projection and the conditional expectation have the same bounded step kernel. Consequently these complete-block sums converge in $L^1$ and almost everywhere for $f\in L^1$. The martingale conclusion concerns the block indices $b^q$.

For a step signal $x=(x_j)_{j<b^q}$, choose compatible tensor ordering and define $U=\overline A^{\otimes q}$. Its coefficient vector is
\[
 \widehat x=b^{-q/2}Ux,\qquad
 \sum_n|\widehat x_n|^2=b^{-q}\sum_j|x_j|^2.
\]
The unitary coordinate vector is $Ux$, not $\widehat x$. For $b=2$ and the real Hadamard digit matrix this is exactly the normalization used earlier. A product of small unitary digit transforms therefore supplies an invertible encoding with exact Euclidean stability. For integer messages an exact real transform still requires a specified output representation; reducing its entries modulo two gives a different operator. In particular the real Sylvester matrix reduces to the rank-one all-one matrix for size at least two.

\paragraph{Encoding and information preservation.}
For a public invertible map $R$ on a finite message space, exact representation preserves full-vector Shannon entropy: $H(R(X))=H(X)$. Orthogonality additionally gives $\|Q^Te\|_2=\|e\|_2$ for a perturbation in unitary coordinates. These are information-preservation and numerical-stability statements. They do not supply secrecy. A pooled histogram of coordinates is a different statistic: its plug-in entropy on a sample of length $L$ over an alphabet of size $a$ is at most $\log_2\min(a,L)$, because at most that many bins are occupied. Neither a uniform histogram nor a large change in selected coefficients is a security definition. If transform preprocessing is followed by a protective cryptographic operation, the latter and its finite encoding determine the security problem.

\paragraph{A basis choice is a model choice.}
Different digit matrices preserve the complete-block Hilbert geometry but change the individual coordinate functions. Similarly, a Walsh-packet or wavelet-packet best-basis procedure \cite{CoifmanWickerhauser} selects a representation by a specified cost. The selection can improve that cost on a signal class; preservation of energy alone cannot prove it improves every cost or every signal. Comparisons must keep fixed the number of retained coefficients, the sampling rule, and the reconstruction metric. This returns to the distinction between a complete unitary representation and nonlinear support selection in Section~\ref{sec:greedy}.

\paragraph{Structured covariance as a finite model.}
The XOR-stationary covariance theorem above identifies the exact statistical model diagonalized by the WHT. Appendix~\ref{app:quadratic-models} records the supplementary quadratic-score identities and the Frobenius projection onto this model. A full positive covariance transforms by orthogonal conjugation; a diagonal or sequency-band restriction selects a smaller family of models. The distinction connects finite spectral structure with estimation while leaving the general coordinate identities outside the main analytic progression.

\section{Noncompact dyadic transforms and affine analysis}\label{sec:noncompact}
The carry-free field model below introduces two-sided scales. Character integration on balls gives inversion on a dense test space. Completeness of its affine Haar system then requires two limits: fine averages recover the function, and coarse averages vanish. Both are proved in the stated Haar normalization.

The compact group $\D$ is not the additive group of $\mathbb Q_2$. Every element of $\D$ has additive order at most two, whereas $\mathbb Q_2$ is torsion-free. A noncompact local field compatible with binary carry-free addition is instead
\[
 G=\F_2((t)),\qquad H=\F_2[[t]],\qquad\mu(H)=1.
\]
Multiplication by $t^m$ scales Haar measure by $2^{-m}$. Define $\chi(x)=(-1)^{[t^{-1}]x}$, where $[t^{-1}]x$ denotes the coefficient of $t^{-1}$. The pairing $\chi(x\xi)$ identifies $G$ with its dual, and $(t^mH)^\perp=t^{-m}H$.

\begin{proposition}\label{prop:noncompact}
For the transform $\mathcal F_Gf(\xi)=\int_Gf(x)\chi(x\xi)\dd(x)$,
\[
 \mathcal F_G\ind_{a+t^mH}(\xi)=2^{-m}\chi(a\xi)\ind_{t^{-m}H}(\xi).
\]
On locally constant compactly supported functions, $\mathcal F_G^2f=f$ and $\norm{\mathcal F_Gf}_2=\norm f_2$. The transform extends unitarily to $L^2(G)$.
\end{proposition}
\begin{proof}
A character integrates to zero on a compact subgroup unless it is trivial there, in which case its integral is the subgroup's measure. This proves the indicator identity. Apply it twice to get inversion; reflection disappears because the characteristic is two. A common support and constancy scale reduces Parseval to finite Fourier orthogonality on $t^{-M}H/t^mH$. Such functions are dense in $L^2(G)$ by support truncation and local averaging, giving the unitary extension.
\end{proof}
For $x\in H$ and $\xi_n=\sum_jn_jt^{-j-1}$, $\chi(x\xi_n)=w_n(x)$. Thus the dual of $H$ is the discrete quotient $G/H$, while $G$ itself is noncompact and self-dual. These statements are compatible, but describe different groups.

The real half-line model replaces $\sum_ka_kt^k$ by $\sum_ka_k2^{-k-1}$ with a fixed binary convention. It preserves measure up to the chosen normalization; addition becomes digitwise XOR, not ordinary real addition. This is the setting of the Walsh--Fourier transform on $\R_+$ in \cite{GES,DyadicBook15}. Fractional calculus or pseudodifferential results on $\mathbb Q_2$ cannot be cited as identical binary Walsh theorems merely because both are nonarchimedean.

The comparison has three layers. Fine's generalized-character framework, the Gibbs--Butzer differentiation programme, and its developments by Schipp, Wade, Su and collaborators are documented in \cite{GES,ButzerWagner,DyadicBook15}. The half-line refinable-function theory of Protasov--Farkov \cite{ProtasovFarkov06} adds translation and dilation, rather than just a reindexing of compact Fourier sums. In our field model, $\mu(t^mA)=2^{-m}\mu(A)$; under the half-line digit map this becomes a real dilation with the same measure factor but with XOR translation. Consequently orthogonality is local character algebra, whereas completeness requires both limiting directions of the conditional expectations, proved below. The local-field statement identifies the exact model behind the analogy.

Let $E_mf$ be averaging over $x+t^mH$, for every integer $m$. Proposition~\ref{prop:noncompact} gives $\mathcal F_G(E_mf)=\ind_{t^{-m}H}\mathcal F_Gf$. Put $\psi=\ind_{tH}-\ind_{1+tH}$. The functions
\[
 \psi_{m,a}(x)=2^{m/2}\psi(t^{-m}x-a),\qquad m\in\Z,\ a\in G/H,
\]
use a fixed set of representatives of $G/H$. Changing representatives may change signs of individual functions; fixing them makes the notation unambiguous. The following proposition proves the orthogonal decomposition and completeness. Its role is to connect the Fourier projections with affine synthesis \cite{ProtasovFarkov06,Terekhin09,AstashkinTerekhin18}.

\subsection{Completeness of the noncompact Haar expansion}\label{sub:v2noncompact}
\begin{proposition}\label{v2:noncompactHaar}
With $\psi_{m,a}$ as above, these functions are an orthonormal basis of $L^2(G)$ and
\[
 (E_{m+1}-E_m)f=\sum_{a\in G/H}\langle f,\psi_{m,a}\rangle\psi_{m,a},
 \qquad f=\sum_{m\in\Z}(E_{m+1}-E_m)f
\]
in $L^2$.
\end{proposition}
\begin{proof}
On each coset of $t^mH$, the mean-zero subspace of functions constant on its two child cosets has dimension one and is spanned by the corresponding normalized $\psi_{m,a}$. This proves the projection identity. Equal-level supports are disjoint; different-level orthogonality follows since a coarse function is constant on the support of a fine mean-zero function.

As $m\to+\infty$, local averaging converges to a locally constant compactly supported function exactly once its constancy scale is reached, and hence to every $L^2$ function by density and contractivity. For the other limit first take such a compactly supported $f$. When $m$ is sufficiently negative, its support is in one coset of $t^mH$, on which the average equals $2^m\int f$. Its squared $L^2$ norm is $2^m|\int f|^2\to0$. Density and contractivity give $E_mf\to0$ for general $f\in L^2$. Telescoping and orthogonality now prove completeness and the asserted two-sided expansion.
\end{proof}

\subsection{Exact local-field tests and the choice of geometry}\label{sub:localfield-tests}

An exact test is a ball indicator. For $f=1_{t^kH}$ and any integer $m$,
\[
 E_mf=
 \begin{cases}
  1_{t^kH},&m\ge k,\\
  2^{m-k}1_{t^mH},&m<k.
 \end{cases}
\]
For $m<k$, its squared norm is $2^{m-2k}$, tending to zero as $m\to-\infty$. The Fourier transform is $2^{-k}1_{t^{-k}H}$. Both identities use the normalization $\mu(H)=1$ and can be checked by integration over a coset. Translation by $a$ introduces the character factor $\chi(a\xi)$ in the transform, rather than an ordinary real exponential.

The finite-quotient proof above extends to $L^2(G)$ by density and isometry. For an $L^1$ function the transform is bounded and continuous, but inversion needs additional hypotheses or a limiting interpretation.

\paragraph{Comparison and research boundary.}
The local-field model preserves character integration and affine synthesis with a two-sided scale parameter. It also fixes the boundary of the analogy: ordinary real addition and the additive group of $\mathbb Q_2$ are different models. The final section uses the same discipline to distinguish genuine extensions from changes of terminology.

\subsection{Wavelet eigenfunctions and the choice of a spectral generator}\label{sub:wavelet-generators}
The noncompact Haar expansion above is an orthogonal decomposition by scale. Turning it into a differential calculus requires a chosen symbol. On $G=\F_2((t))$, the radial Fourier multiplier $|\xi|^\alpha$, where $\alpha>0$ and $|t|=1/2$, is constant with value $2^{\alpha(m+1)}$ on
\[
 t^{-(m+1)}H\setminus t^{-m}H.
\]
Since $E_{m+1}-E_m$ projects onto this annulus, its Haar functions are eigenvectors with that eigenvalue. On $L^2(G)$ the operator is defined by the closed multiplier domain
\[
 \int_G|\xi|^{2\alpha}|\mathcal F_Gf(\xi)|^2\,d\mu(\xi)<\infty.
\]
The projection identity proves the spectral assertion directly. It uses a radial scale symbol, rather than the integer symbol $n$ of the compact numerical-frequency generator.

There is a related but distinct $p$-adic theorem. Normalize Haar measure on $\mathbb Q_p$ by $\mu(\mathbb Z_p)=1$ and choose the standard additive character $\chi_p$ trivial on $\mathbb Z_p$. Kozyrev \cite{Kozyrev02}, Theorem 2, gives the orthonormal basis
\[
 \psi_{\gamma j a}(x)=p^{-\gamma/2}
 \chi_p(p^{\gamma-1}jx)\ind_{\{|p^\gamma x-a|_p\le1\}},
\]
where $\gamma\in\Z$, $a\in\mathbb Q_p/\mathbb Z_p$ has the fixed finite negative-digit representative, and $j=1,\ldots,p-1$. For the Vladimirov multiplier $\mathcal F_p^{-1}|\xi|_p^\alpha\mathcal F_p$, $\alpha>0$,
\[
 D^\alpha\psi_{\gamma j a}=p^{\alpha(1-\gamma)}\psi_{\gamma j a}.
\]
This is an external spectral theorem. The normalization is visible already from the support measure $p^\gamma$, which makes each displayed function have norm one. The theorem concerns $\mathbb Q_p$, not the compact binary group.

\begin{longtable}{@{}p{.24\textwidth}p{.30\textwidth}p{.33\textwidth}@{}}
\caption{Three spectral generators with different frequency geometry}\label{tab:generator-models}\\
\toprule Model & Spectral parameter & Consequence of the chosen geometry\\\midrule\endfirsthead
\toprule Model & Spectral parameter & Consequence of the chosen geometry\\\midrule\endhead
Compact numerical calculus & $n$ on $w_n$ & Finite-dimensional frequency cutoffs and compact resolvent.\\
Radial calculus on $\F_2((t))$ & $2^{\alpha(m+1)}$ on the Haar scale $m$ & Frequency annuli determine scale eigenspaces.\\
Vladimirov calculus on $\mathbb Q_p$ & $p^{\alpha(1-\gamma)}$ on Kozyrev wavelets & A $p$-adic Fourier symbol determines the wavelet eigenvalues.\\\bottomrule
\end{longtable}
The compact product and the characteristic-two field have carry-free addition, while $\mathbb Q_p$ has characteristic zero. A measure-preserving coding can identify Hilbert spaces or selected bases without identifying the additive groups. The operator must consequently be transported as well as the functions. This is why a result about a Vladimirov operator does not follow merely by renaming the compact numerical derivative.

Group wavelets also allow filters beyond the Haar choice. Lang \cite{Lang96} constructs multiresolution analyses and orthogonal wavelets on the Cantor dyadic group. Protasov--Farkov \cite{ProtasovFarkov06} study dyadic refinable functions on the half-line, and the Banach-space synthesis questions in \cite{Terekhin09,AstashkinTerekhin18} impose additional basis conditions. In the Haar case, one mean-zero direction on every binary parent supplies the exact projection formula proved above. With a different filter, orthogonality, completeness, regularity, and boundedness of the synthesis map must each be established. They are separate additions to the finite character algebra.

\section{Research boundaries and quantitative questions}\label{sec:boundaries}
The research questions below are organized by the implication that is missing from an established result. An endpoint question, a proposed refinement, and a model-design programme have different status. Universality questions receive explicit quantifiers because changes in coefficient origin, alphabet, order, or target topology can either settle a question or make it impossible.

The preceding results separate three sources of strength: localization in a fixed order, freedom to select frequencies, and freedom to construct a representation for the target. Binary shells give sharp Hardy thresholds; interval trees supply control beyond complete blocks; correction and universality change the quantifiers of the approximation problem. The finite applications use exact character cancellation followed by additional smoothness, sampling, or probabilistic assumptions. The questions below concern the remaining analytic or constructive step in each setting.

\subsection{Precise problems, established boundaries, and proposed extensions}\label{sub:precise-problems}
We separate the classical endpoint question from proposed refinements and model-design programmes. The following reductions identify conclusions already supplied by the results above.

\paragraph{Problem 1: a specified full-sequence endpoint.}
For the Walsh--Paley maximum $S^*$, consider the proposed weak Luxemburg-norm estimate
\[
 \|S^*f\|_{1,\infty}\le C\|f\|_{L\log L},
 \qquad \Phi_1(t)=t\log(e+t).
\]
The reference point is the verified sufficient Antonov-type class $\Phi_A$ in Section~\ref{sub:endpoint-refinement}, together with the quantitative weak-$L^p$ estimate. The estimate displayed here is a question, not a consequence of either result. A sufficient proof would give almost everywhere convergence by polynomial density in the indicated Orlicz norm and the maximal approximation argument. A negative result would have to address this specified estimate or construct divergence in the specified class; divergence in arbitrary $L^1$ does not do so.

A smaller sufficient Young-function growth than $\Phi_A$, an improved distribution estimate, or divergence in the specified class would address this endpoint. A trigonometric result requires a separate Walsh argument.

\paragraph{Programme A: geometry beyond the settled endpoint-count growth.}
For $\Lambda\subset\N$, write $A_s(\Lambda)$ for the endpoints of its maximal binary strings in scale $[2^s,2^{s+1})$, as in Section~\ref{sec:hardy}, and define
\[
 \mathfrak C(\Lambda)=\sup_{F\ne0}
 \frac{\|\sup_{n\in\Lambda}|\sigma_nF|\|_{1/2}}{\|F\|_{H^{1/2}}}.
\]
The weighted theorem of \cite{ABPTWeighted23}, Theorem 3.1(a), gives
\[
 \mathfrak C(\Lambda)\lesssim K^2
 \quad\hbox{if}\quad\sup_s|A_s(\Lambda)|\le K.
\]
This worst-case order is already sharp. Indeed, for
$\Lambda_M=\{2^M+2^l:1\le l<M\}$ and the unit $H^{1/2}$ atom from Theorem~\ref{th:fejerhardy}, its endpoint calculation gives
$\int(\sup_{n\in\Lambda_M}|\sigma_na_M|)^{1/2}\,d\mu\gtrsim M$.
Hence $\mathfrak C(\Lambda_M)\gtrsim M^2$, while $|A_M(\Lambda_M)|\asymp M$. Thus
\[
 \sup_{\Lambda:\,\sup_s|A_s(\Lambda)|\le K}\mathfrak C(\Lambda)
 \asymp K^2 \qquad(K\ge K_0).
\]
This is a consequence of the cited theorem and the existing atom calculation, not an open growth problem. A finer programme is to distinguish families with the same endpoint count by their endpoint spacing and cross-scale organization, or to determine exact constants for specified finite configurations. A proposed invariant should first be tested against $\Lambda_M$ and the restricted-maximal results \cite{BPTRestricted23,ABPTWeighted23}.

\paragraph{Proposed question: same-space compactness outside the regular class.}
For fixed $1<p_0<\infty$, $p_0\ne2$, suppose $a\in c_0$ and $T_a$ is bounded on $L^{p_0}$. A precise question raised here is whether boundedness at this exponent together with $a\in c_0$ implies operator-tail decay at the same exponent. A positive result beyond the regular class would identify an explicit hypothesis weaker than \eqref{eq:marcink}. The Marcinkiewicz case is already settled by Proposition~\ref{prop:marccompact}; a vanishing-variation assumption would be stronger than necessary. Boundedness at an exponent beyond $p_0$ toward the same side of two also settles the question by interpolation. A useful result must consequently identify a genuine boundary class, prove an endpoint tail estimate, or construct a bounded null-symbol obstruction. This is an operator-tail programme; it is not settled by scalar decay or a formal coefficient estimate.

\paragraph{A representation constraint leading to Question~\ref{ques:strict-diagonal}.}
Fix a solid rearrangement-stable sequence space $X\subset c_0$, a positive nonincreasing profile $A$, and natural Walsh order. Question~\ref{ques:strict-diagonal} formulates a representation problem with these data. For sign-only representation one must use every coefficient magnitude; for subseries representation one cannot reverse a coefficient's sign. Determine a criterion only after choosing one of these modes and a target class. The square-summability obstruction remains necessary in all three modes, but necessity alone proves no existence statement.

Fourier origin should be tracked as an additional conclusion. A useful result may prescribe an $L^p$ class for one source function, a common exceptional set, or a quantitative coefficient envelope, but it must exhibit the simultaneous construction. Combining separately existential statements with different coefficient sequences is not enough. Comparison with \cite{GevNava99,Gri20,Gri24} is needed before any claim of novelty.

\paragraph{Programme B: integration and recovery in one finite design.}
Fix a digital integration space $H_r$, a class of node subgroups or generating matrices, and a sparse signal model $\Sigma$. Specify a measurement matrix $M$ built from the same finite character table. Seek a design for which both
\[
 \sum_{k\in P^\perp\setminus\{0\}}r(k)\le\varepsilon^2,
 \qquad \|Mz\|_2\ge\kappa\|z\|_2
       \quad(z\in\Sigma-\Sigma)
\]
hold with stated costs and $\kappa>0$. This is a proposed joint design problem, not an equivalence between two known theories. The subgroup-indicator obstruction in Section~\ref{sub:qmc-bridge} must be excluded or accommodated. If $\Sigma$ is a fixed subspace, the second condition reduces to a smallest-singular-value calculation; for a union of sparse subspaces it is a uniform restricted estimate. These are different levels of difficulty.

\paragraph{Two analytic directions with a fixed statement.}
The optimal universal lacunary Walsh class remains the source-listed question in \cite{GogOniani21}; Table~\ref{tab:orlicz-gap} records the distinct upper result for each lacunary sequence and the lower obstruction for particular nested spectra. Removing the factor $\log_4L$ is therefore a research direction, not a theorem supplied by the trigonometric classification.

For quantitative weights, one must fix the exact class of weights and operator before asking for optimal powers. For example, with $S^*$ in Paley order and dyadic $A_p^d$, define the proposed extremal profile
\[
 \mathcal C_p(A)=\sup_{\substack{w\in A_p^d\\ {[w]_{A_p^d}\le A}}}
       \|S^*\|_{L^p(w)\to L^p(w)},\qquad A\ge1.
\]
Determine its growth up to specified powers or logarithms, comparing any new bound with the precise weighted theorem used. An upper exponent quoted in Section~\ref{sec:timefreq} is not a matching lower bound. This profile is a proposed way of stating the quantitative task, rather than an assertion that every parameter range lacks an answer. For N\"orlund means, the corresponding task must first fix the coefficient class, the $L^\infty$ hypothesis, and whether the endpoint is strong, weak, or normalized. A conjecture about all such means without these data would conflate different problems.

\subsection{Universality questions with specified data}\label{sub:universality-questions}
The following statements separate an open question identified in a source from a proposed extension of the present comparison. An incomplete argument is neither evidence for a theorem nor evidence that its conclusion is open in the literature. The adjacent established results specify what a solution would add.

\begin{question}[Source-listed: Fourier origin and permutations]\label{ques:permutations}
Does there exist $U\in L^1(\D)$ such that, for every real measurable $h$ finite almost everywhere, there is a permutation $\pi_h$ of $\N_0$ for which
\[
 \sum_{j=0}^{N-1}\widehat U(\pi_h(j))w_{\pi_h(j)}(x)
 \longrightarrow h(x)\quad\text{for almost every }x?
\]
Every coefficient is used exactly once, without changing its sign, and the source function is fixed before $h$.
\end{question}
This is the finite-target part of Question 1 in \cite{Gri24}. The sign theorem there settles a different operation. The weighted result \cite{Sargsyan25} allows a prior sign modification and a selected weight. Neither is an answer to the displayed question. A ternary representation theorem permits a different operation and likewise does not answer it. The source posed the problem in 2024, and the nearby results identified here do not resolve it.

\begin{question}[Source-motivated: unbounded digit sequences]\label{ques:unbounded-digits}
Let $m_j\ge2$ be an arbitrary sequence of integers, possibly unbounded. For a fixed $q>2$, under what necessary or sufficient condition on $\mathbf m$ does there exist a formal $\ell^q$ Vilenkin series satisfying the full common-exceptional-set conclusion of Theorem~\ref{th:ordinary-vilenkin}, with unmodified ordinary sums and all $1\le p<\infty$?
\end{question}
Remark 3.12 of \cite{EpiVilenkin26} points to unbounded generating sequences; the displayed question fixes the desired global conclusion. The known theorem covers $\sup_jm_j<\infty$. At the proof level the unresolved step is replacement of the bounded overshoot in the finite digit-block selection by a sufficient estimate for the given sequence. An uncontrolled factor depending on the terminal digit cannot be treated as a uniform constant. A concrete first test is a generating sequence with $m_j\le C(j+2)^A$ for fixed $C,A$: one must bound the last-digit overshoot inside the localization budget, rather than hide it in a constant. This is a proposed tractable model, not an assertion that it is unresolved independently of the literature. A proof for a specified unbounded class would have to be compared with earlier unbounded-system constructions. The diagonal $\bigcap_{q>2}\ell^q$ improvement in the bounded case is not part of this unresolved step.

\begin{question}[Proposed here: a stricter diagonal representation mode]\label{ques:strict-diagonal}
Fix a solid rearrangement-stable real sequence space $X\subset c_0(\mathbb R)$ and a positive nonincreasing real profile $A$. Is there a necessary and sufficient criterion for the existence of $c\in X$, $d=Ac=\widehat U$, with one real $U\in\bigcap_{1\le p<2}L^p$, such that every real measurable target finite almost everywhere is represented in natural Walsh order using only signs $\{-1,1\}$? Pose separately the same question with selectors $\{0,1\}$.
\end{question}
The condition $A\notin M(X,\ell^2)$ is necessary for either mode by the Parseval--Fatou obstruction. No sufficient classification for this full combination of constraints is asserted here. The sign-only existence result \cite{Gri24} does not prescribe this factor/profile/Fourier-integrability combination. A solution must retain the same profile before the target; a construction choosing a new profile for each target answers a different question. As a first test, take $X=\ell^s(\mathbb R)$, $s>2$, and $A_n=(n+1)^{-\beta}$. If $\beta>1/2-1/s$, then $A\in\ell^r$ for $1/r=1/2-1/s$, so H\"older forces $Ac\in\ell^2$ and rules out the requested universality. The remaining parameter range requires an existence argument with the chosen alphabet and source integrability. This gives a fixed model on which to test a proposed criterion; no claim is made that every such specialization is a new open problem.

\begin{question}[Proposed here: joining universalizing weights]\label{ques:joining-weights}
For one fixed formal series $C$, let $u,v\in\mathfrak U(C)$. Which verifiable compatibility condition on $C,u,v$ ensures $\max\{u,v\}\in\mathfrak U(C)$? Equivalently, when does one obtain, for every target in the common weighted space, one increasing sequence of ordinary stopping indices converging in both weighted norms?
\end{question}
The equivalence uses
\[
 \max\{\|h\|_{L^1(u)},\|h\|_{L^1(v)}\}
 \le\|h\|_{L^1(\max\{u,v\})}
 \le\|h\|_{L^1(u)}+\|h\|_{L^1(v)}.
\]
Comparable weights are already settled by Proposition~\ref{prop:weight-inheritance}. Compatible common approximation skeletons give a sufficient mechanism. Separate universality in the two norms supplies no immediate simultaneous stopping indices. The question concerns the remaining compatibility structure of a fixed orbit; it does not assert that an arbitrary positive weight with essential infimum zero is universalizing.

There is also an elementary necessary condition for joining. If $\max\{u,v\}\ge a>0$ almost everywhere, its norm is equivalent to the unweighted $L^1$ norm, so it cannot universalize one ordinary-sum orbit. Thus
\[
 \operatorname*{ess\,inf}\max\{u,v\}=0
\]
is necessary for a positive joining conclusion. The separate conditions $\operatorname*{ess\,inf}u=\operatorname*{ess\,inf}v=0$ do not imply it: $u$ and $v$ can be bounded below on complementary sets. This observation supplies an obstruction test, not an example asserting separate universality for those weights. A proposed compatibility condition should first survive this complementary-set test and then produce common stopping indices.

\subsection{Which apparent gaps are already closed or obstructed}\label{sub:excluded-questions}
Several tempting questions disappear once the quantifiers and the local construction are fixed.

\paragraph{One coefficient sequence for all $q>2$.}
For the bounded-digit ordinary-sum theorem this is Corollary~\ref{cor:allq-vilenkin}. The analogous diagonal device already occurs in the proof of Theorem~\ref{th:weighteduniv}. It is not an open endpoint problem.

\paragraph{A vanishing near-square-summability gauge.}
For the binary block, \eqref{eq:gauge-localization} allows any nondecreasing $\omega$ tending to zero, after the target step function and exceptional measure are fixed. Such a gauge does not imply $\ell^2$. A logarithmic restriction in a different construction should first be tested against this deterministic example before it is called necessary.

\paragraph{Integrable Fourier origin for ordinary-sum universality on a large set.}
Proposition~\ref{prop:ordinary-Fourier-obstruction} excludes this binary formulation. The sign-modified and formal-series theorems evade the obstruction by different mechanisms. Leaving Fourier origin unproved in an ordinary-sum argument does not leave an open theorem of this form.

\paragraph{One continuous function universal on every perfect null set.}
This unrestricted proposal is impossible. A continuous function belongs to $L^2$, so its full Walsh sums converge to it almost everywhere. For every point $x$ of the compact binary group there is a perfect Haar-null compact set containing $x$: fix infinitely many coordinates to those of $x$ and leave infinitely many others free. Universality for all continuous targets on every such set would, by choosing the constant target $f(x)+1$, force failure of convergence at every point. This contradicts the almost-everywhere theorem. Fixing a countable set first, as in Theorem~\ref{th:muller}, is consistent. Generic compact-set conclusions in \cite{EpiMuller14} also have a different quantifier order.

\paragraph{Lower order ideals of weights.}
The inheritance $0<v\le A u$ from one universalizing weight is Proposition~\ref{prop:weight-inheritance}, using continuity and density. The joining question above is the stronger step. Separate existence of several weights does not by itself create one common subsequence.

These deductions record boundaries of plausible research programmes. A proposed extension must specify an independently meaningful conclusion and its relation to known results; a failure of one proof strategy does not establish that the conclusion is open.

\subsection{A problem register and the role of applications}\label{sub:problem-register}
\begin{longtable}{@{}p{.23\textwidth}p{.27\textwidth}p{.37\textwidth}@{}}
\caption{Questions distinguished by status and the missing implication}\label{tab:problem-register}\\
\toprule Question or programme & Status in this survey & Known input and additional step\\\midrule\endfirsthead
\toprule Question or programme & Status in this survey & Known input and additional step\\\midrule\endhead
Full-sequence $L\log L$ maximum & Specified endpoint question above. & Near-endpoint sufficient classes do not prove the proposed weak Luxemburg estimate.\\
Restricted Fej\'er geometry & Proposed quantitative refinement. & Worst-case endpoint-count growth $K^2$ is settled; spacing must supply additional information.\\
Same-exponent multiplier tails & Proposed operator question above. & Marcinkiewicz interpolation settles the regular class; boundedness at only one exponent leaves a separate tail issue.\\
Fourier-origin permutations & Source-listed, Question~\ref{ques:permutations}. & Sign-only and weighted rearrangement theorems use different data.\\
Unbounded digits & Source-motivated, Question~\ref{ques:unbounded-digits}. & The finite digit-block estimate has to survive uncontrolled overshoot.\\
Strict diagonal modes & Proposed, Question~\ref{ques:strict-diagonal}. & The square-summability obstruction alone gives no sufficient criterion under a prescribed profile.\\
Joining weights & Proposed, Question~\ref{ques:joining-weights}. & Downward inheritance is settled; simultaneous return times need compatibility.\\
Integration and recovery & Proposed joint design programme above. & Annihilator cancellation and restricted injectivity are different requirements.\\\bottomrule
\end{longtable}

The application questions are model-design questions rather than unresolved character identities. A Walsh covariance approximation should specify a population class and compare its estimation risk with a full or alternative structured estimator; a spectral cryptanalytic invariant should identify an equivalence relation and a class of attacks it informs. For generalized digit bases, a block convergence argument must not be mistaken for a maximal theorem over all indices. These are concrete ways to extend the mathematical mechanisms in Section~\ref{sec:finite}, but no uniform application superiority is conjectured.

The resulting research programme has three distinct levels: sharpen a known analytic boundary; change a prescribed representation constraint while retaining its quantifiers; or add a statistical/finite-state model to an exact transform identity. Stating the level makes both the usefulness of an extension and the evidence needed to establish it assessable.

\section*{Declaration of AI-assisted preparation}
AI-assisted tools were used for editorial restructuring, literature retrieval, preparation of applied exposition, and consistency checking. Responsibility for the mathematical statements, source attributions, and final submitted text rests with the author.

\appendix

\section{Orthogonal and independent-character tests}\label{app:orthogonal}
\subsection{What orthogonality gives under a permutation}\label{sub:v2RM}
\begin{theorem}[Rademacher--Menshov]\label{v2:RMfull}
For an orthonormal sequence $(\phi_n)$ on a probability space,
\[
 \sum_{n\ge1}|a_n|^2\log^2(n+2)<\infty
\]
implies almost everywhere convergence of $\sum_na_n\phi_n$ in that order.
\end{theorem}
We use this classical theorem as an external input; see \cite{GES,SWS}. Its finite mechanism is the decomposition of a prefix into at most one interval of each dyadic length, which gives an $L^2$ maximal loss $O(\log N)$. It applies after a fixed permutation, but the weighted coefficient hypothesis must then be checked in that order. Natural-order Carleson convergence does not supply a bound uniform over all permutations. Modern Weyl-multiplier questions are discussed in \cite{KaragulyanLaceyNavoyan25}.

\subsection{The independent-character endpoint}
Walsh characters are independent random signs exactly when their frequency vectors are linearly independent over $\F_2$. Surjectivity of the finite character-value map proves sufficiency; a nonempty character product equal to one proves necessity. For real $a_j,t$, Khintchine and $\prod_j\cosh(ta_j)\le\exp(t^2\sum_ja_j^2/2)$ give uniform $\exp L^2$ control when $\sum_j a_j^2=1$. For complex coefficients with $\sum_j|a_j|^2=1$, apply the real estimate separately to the real and imaginary parts; the union bound and $|z|^2=(\Re z)^2+(\Im z)^2$ give the same Orlicz bound up to an absolute constant. By Orlicz duality,
\[
 \left(\sum_j|\widehat f(n_j)|^2\right)^{1/2}
 \lesssim\|f\|_{L(\log L)^{1/2}}
\]
for an independent family; see \cite{SagherZhou,Carrillo}. Numerical lacunarity must still be checked against the actual Walsh algebra.

The logarithmic exponent has an elementary sharp test. For $f_M=D_{2^M}=2^M\ind_{I_M}$, the first $M$ Rademacher coefficients are one, giving coefficient norm $\sqrt M$. For $\Phi_\beta(t)=t\log^\beta(e+t)$, $0\le\beta<1/2$, its Luxemburg norm is comparable to $M^\beta$, by inserting the indicator into $2^{-M}\Phi_\beta(2^M/\lambda)=1$. The ratio diverges. This identifies the coefficient geometry behind the logarithm in endpoint multi-frequency estimates, rather than deriving such estimates from orthogonality alone.

The same independent-character test rules out a uniform synthesis inequality $\ell^{p'}\to L^p$ for $1<p<2$. Coefficients equal to one at $M$ Rademacher frequencies have synthesis norm $\asymp_p M^{1/2}$ and coefficient norm $M^{1/p'}$, with $1/p'<1/2$. This rules out that coefficient shortcut without deciding a same-space compactness statement under an independent boundedness hypothesis.

\clearpage
\section{Published-source scope and normalization register}\label{app:source-register}
The following pointers record hypotheses and conventions that matter when comparing the published results. An open author copy provides access to an identified publication; the publication is the bibliographic record. A source pointer identifies an external input, whereas the proofs supplied in this survey remain separately visible.
{\small
\begin{longtable}{@{}p{.22\textwidth}p{.30\textwidth}p{.37\textwidth}@{}}
\toprule Result & Published source & Scope retained\\\midrule\endhead
N\"orlund criterion & \cite{GogNorlund23}, Theorem 3 & $1/2<p\le1$; positive decreasing coefficients and the stated $L^\infty$ maximal bound.\\
Matrix kernel bound & \cite{BlahotaGat25}, Theorems 1--2 & Nonnegative normalized rows; the square-sum condition controls kernel $L^1$ norms. Lacunarity is additional in Theorem 4.\\
Separated subsequences & \cite{Konyagin93}, Theorem 1 & Divisibility between successive binary supports; unbounded Lebesgue constants can coexist with a.e. convergence for every $L^1$ input.\\
Variable exponent maximum & \cite{HaoLiWeisz26}, Theorem 1.4 & $p_->1/2$, local oscillation condition, $1/p_--1/p_+<1$, and $q\in(0,\infty]$.\\
Interval synthesis & \cite{Osipov16}, Theorem 1 & Finite disjoint families; uniform family constant.\\
Weighted variation & \cite{DoLaceyWeighted}, Theorem 1.2 & $w\in A_a^d$, $1\le a<p$, $r>2a$, $1/r<1/a-1/p$; translate the source's inclusive partial-sum index.\\
Weighted Walsh maximum & \cite{DiPlinioLerner14}, Theorem 1.4(ii) & Exponent $\max\{p',2/(p-1)\}$ for Walsh; keep separate from the logarithmic factor in the Fourier bound.\\
Weighted interval bound & \cite{DFPR25}, Theorem D, \S5.3 & The proof supplies the $L^p(w)$ input norm used here; dyadic characteristics and the stated operator are fixed.\\
Compact interpolation & \cite{Krasnoselskii60} & Qualitative compactness by interpolation; Proposition~\ref{prop:marccompact} also gives an explicit norm-tail estimate.\\
Quartile endpoint & \cite{DemeterDiPlinio13} & Mixed endpoint at output exponent $2/3$; retain its Lorentz input.\\
Critical restricted Fej\'er maximum & \cite{ABPTWeighted23}, Theorem 3.1(a) & Endpoint-count normalization; journal online version used, with volume and pages unassigned.\\
Correction block & \cite{Epi12} & Binary localization; the assembly argument is proved here.\\
Common correction set & \cite{GriNava16}, Theorem 3 & The set and strictly decreasing positive coefficient profile precede the arbitrary $L^1$ target.\\
Fourier-origin sign universality & \cite{Gri24}, Theorem 3, Question 1 & The original Fourier series converges; target-dependent signs universalize it. Permutations are a separate source question.\\
Countable-set coefficient recovery & \cite{Plotnikov10}, Theorems 7, 10, 13, 17 & All integer cubes; exact recovery for an integrable finite sum outside a countable set.\\
Bounded-digit universality & \cite{EpiVilenkin26}, Theorem 3.3, Lemma 2.12 & Fixed $q>2$; one exceptional set for all $p$ and targets. The numbered statements were also compared with the author's text.\\
Triangular maximum & \cite{GatTriangle18}, Theorems 1.1, 2.10 & Triangular averaging and a one-index maximal operator; distinct from unrestricted rectangular maxima.\\
Lebesgue-function divergence & \cite{Getsadze26}, Theorem 5 & Uniformly bounded orthonormal system; Lebesgue function divided by $\log|A_n|$ tends to infinity a.e.\\
Digital rate & \cite{GSY17}, Theorem 1, Corollary 1 & RMS worst-case error with fixed dimension, order, weights, and the required $\beta\ge2\alpha$.\\
Wavelet spectral theorem & \cite{Kozyrev02}, Theorem 2 & $\mathbb Q_p$ and the Vladimirov multiplier; the displayed representative and Haar normalization are fixed.\\
Generalized bases & \cite{DutkayPicioroaga14} & Unitary digit matrices; conditional expectations at complete digit blocks.\\\bottomrule
\end{longtable}
}

\clearpage
\section{Proofs of auxiliary statements}\label{app:auxiliary-proofs}
The following are proofs of known results or elementary consequences, included to reproduce the normalizations used in the main text. They contain no claim of new results. The principal threshold and constructive proofs remain in their respective sections.

\subsection{Derivation for Proposition~\ref{v2:absolute}}
\begin{proof}
Cauchy--Schwarz on a block gives
\[
 \sum_{2^j\le n<2^{j+1}}|\wh f(n)|\le2^{j/2}\norm{\Delta_jf}_2.
\]
Sum for the first assertion. For $s>1/2$, apply sequence H\"older with the geometrically decaying factor $2^{-j(s-1/2)}$. For failure at $s=1/2$, choose a positive $a\in\ell^q\setminus\ell^1$ and put $\wh f(n)=a_j2^{-j}$ throughout block $j$. For finite $q>1$, $a_j=(j+1)^{-\theta}$ with $1/q<\theta\le1$ works; for $q=\infty$ take $a_j=1$. The squared coefficient sum is $\sum a_j^22^{-j}<\infty$, defining an $L^2$ function. Its block Besov sequence is exactly $a_j$, while its absolute coefficient sum is $\sum a_j=\infty$.
\end{proof}

\subsection{Derivation for Proposition~\ref{prop:compacttail}}
\begin{proof}
If the tails vanish, $E_mT_a$ gives finite-rank norm approximants. Conversely, $E_m\to I$ strongly, with uniformly bounded norms. This convergence is uniform on compact subsets: cover a compact set by a finite $\eps$-net and use the operator bound on the remainder. Apply this to the compact closure of $T_a$ of the unit ball.
\end{proof}

\subsection{Derivation for Proposition~\ref{v2:BV}}
\begin{proof}
Write $a_k-L=\sum_{n\ge k}(a_n-a_{n+1})$. Since the natural partial-sum operators have a common $L^p$ norm, the operator series
\[
 T_a-LI=\sum_{n\ge0}(a_n-a_{n+1})S_{n+1}
\]
converges absolutely in operator norm. Each summand is finite rank, so its limit is compact. If $L\ne0$, subtracting a compact operator would make $LI$ compact were $T_a$ compact, impossible in infinite-dimensional $L^p$.
\end{proof}

\subsection{Derivation for Proposition~\ref{v2:Wiener}}
\begin{proof}
A zero coefficient is preserved by translations up to its character factor, so its continuous functional annihilates the span. Conversely, let $M$ be the closed span. The map $h\mapsto T_hf$ is norm-continuous in $L^1$. Its Bochner integral against $w_n(h)$ is a limit of finite linear combinations of translates and equals $\wh f(n)w_n$. Thus every $w_n$ belongs to $M$ when all coefficients are nonzero. Polynomial density proves $M=L^1$.
\end{proof}

\clearpage
\section{Supplementary quadratic models and XOR covariance}\label{app:quadratic-models}
The formulas below are elementary consequences of orthogonal coordinates and the finite character theorem in Section~\ref{sub:v2stationarity}. They serve as model checks; their inclusion carries no claim of a new detection or clinical theorem.
\subsection{Orthogonal coordinates and quadratic signal models}\label{sub:applied-ecg}

A full quadratic score gives an elementary boundary test for any claim that a complete orthogonal transform improves discrimination. Let a signal window be $x\in\R^N$ and let $Q$ be any orthogonal matrix, including a sequency permutation of $Q_N$. The calculation below concerns only the chosen reference mean and covariance; it makes no empirical or application-specific claim.

Let $\widehat\mu$ and $\widehat\Sigma$ be the time-domain empirical mean and unbiased covariance of reference windows. For $\delta>0$, set
\[
 R_\delta=\widehat\Sigma+
 \delta\frac{\operatorname{tr}\widehat\Sigma}{N}I.
\]
If $\widehat\Sigma\ne0$, then $R_\delta$ is positive definite. In the zero-covariance case a separate positive absolute ridge is needed. Transforming the same reference sample gives $\widehat\mu_Q=Q\widehat\mu$ and $R_{Q,\delta}=QR_\delta Q^T$, because trace and isotropic regularization commute with orthogonal conjugation. Therefore
\begin{equation}\label{eq:applied-ecg-invariance}
 (Qx-Q\widehat\mu)^TR_{Q,\delta}^{-1}(Qx-Q\widehat\mu)
 =(x-\widehat\mu)^TR_\delta^{-1}(x-\widehat\mu).
\end{equation}
The equality is exact for the empirical model, not merely a population approximation. Full-covariance Mahalanobis scoring with this regularization is consequently unchanged by a complete WHT.

Symmetric whitening gives a canonical back-projected residual
\begin{equation}\label{eq:applied-ecg-decomposition}
 r=Q^TR_{Q,\delta}^{-1/2}Q(x-\widehat\mu)
   =R_\delta^{-1/2}(x-\widehat\mu),\qquad
 d_\delta^2(x)=\sum_{n=1}^Nr_n^2.
\end{equation}
The second equality follows from the spectral functional calculus under orthogonal conjugation. Thus not only the scalar score but this symmetric full-covariance attribution is coordinate invariant. Replacing the symmetric square root by a Cholesky factor preserves the scalar quadratic form but need not preserve the same coordinate attribution. The terms $r_n^2$ decompose the chosen quadratic form; their interpretation requires an additional signal model.

A specifically Walsh model appears when covariance structure is imposed after the transform. For mutually orthogonal sequency-band projections $\Pi_b$, consider
\[
 \widetilde R_Q=\sum_b\Pi_b\widehat\Sigma_Q\Pi_b+\lambda I,
 \qquad\lambda>0.
\]
This retains within-band covariance and removes cross-band covariance. Its pullback generally differs from the time-domain full covariance, so the model now depends on the bands and their ordering. Choosing diagonal covariance, truncating coordinates, or weighting bands likewise introduces a substantive spectral assumption. Computational convenience is one motivation; predictive benefit requires evidence under a specified signal model and validation design.

The exact invariant is the full quadratic form. A band approximation is a new estimator, whose bias and variance depend on the population class and the available reference sample. The distinction is mathematical: a unitary change of coordinates and a structural restriction on covariance are different operations. It is the latter, together with sampling assumptions, that can change statistical performance.

\subsection{Projection onto XOR covariance and concentrated contrasts}\label{sub:applied-xor-fit}

The covariance theorem in Section~\ref{sub:v2stationarity} gives an elementary orthogonal projection formula. For a matrix $A$ indexed by $\F_2^m$, put
\begin{equation}\label{eq:applied-xor-proj}
 c_A(u)=\frac1N\sum_xA_{x,x\oplus u},\qquad
 [\mathcal P_{\rm xor}A]_{x,y}=c_A(x\oplus y).
\end{equation}
The matrices supported on each fixed XOR diagonal have disjoint supports and squared Frobenius norm $N$. Averaging on these supports is therefore the Frobenius-orthogonal projection onto XOR-stationary matrices. With the normalized Sylvester matrix,
\[
 Q_N\mathcal P_{\rm xor}(A)Q_N^T
 =\operatorname{diag}\bigl(\operatorname{diag}(Q_NAQ_N^T)\bigr),
\]
and
\[
 \|A-\mathcal P_{\rm xor}(A)\|_F
 =\|\operatorname{offdiag}(Q_NAQ_N^T)\|_F.
\]
For a positive semidefinite covariance the projection remains positive semidefinite, since the retained diagonal entries are nonnegative. This provides a model-distance diagnostic, not an automatic justification for discarding correlations. A learned feature permutation changes which raw feature pairs share an XOR diagonal, so searching an ordering is itself model fitting and belongs in the independent training stage.

There is also a direct connection to nonlinear approximation. Let a transformed mean shift have standardized entries $v_j=d_j/\sigma_j$. If its largest $k$ entries retain at least $(1-\eta)\|v\|_2^2$, then
\[
 \|v\|_\infty\ge\sqrt{\frac{1-\eta}{k}}\,\|v\|_2.
\]
This follows because the energy of any $k$ coordinates is at most $k\|v\|_\infty^2$. Concentration of standardized contrasts therefore gives an exact comparison between maximum-coordinate and total-energy criteria. A detection probability requires a separate noise model and threshold; the displayed inequality alone is deterministic.

\paragraph{Comparison and research boundary.}
The finite applications distinguish exact character identities from the assumptions that create a task guarantee. Digital cancellation needs a smoothness class for an integration error bound; sparse recovery needs restricted injectivity; covariance diagonalization needs XOR stationarity; a selected Boolean correlation has a component-level interpretation. The complete-block and covariance identities apply to their specified models. A change of ordering or covariance class changes the estimation problem.

\end{document}